\documentclass[11pt]{article}
\usepackage[utf8]{inputenc}

\usepackage{geometry, amsmath, amsthm, latexsym, amssymb, graphicx, verbatim, enumerate, mathrsfs,  tikz-cd, bbm, upgreek, accents}
\usepackage{hyperref}
\hypersetup{colorlinks=true}
\usepackage{todonotes}
\usepackage{stackengine}
\usepackage{scalerel}
\usepackage{standalone}

\usepackage{ytableau}
\usepackage{youngtab}

\usepackage{wrapfig}

\usepackage{float}

\usepackage{caption, subcaption}
\newcommand\explainalot[1]{}

\newcommand\justthesis[1]{}

\newtheorem{innercustomthm}{Theorem}
\newenvironment{customthm}[1]
  {\renewcommand\theinnercustomthm{#1}\innercustomthm}
  {\endinnercustomthm}
  
\newtheorem{innercustomlemma}{Lemma}
\newenvironment{customlemma}[1]
  {\renewcommand\theinnercustomlemma{#1}\innercustomlemma}
  {\endinnercustomlemma}
  
\newtheorem{theorem}{Theorem}
\newtheorem*{theorem*}{Theorem}
\newtheorem{lemma}{Lemma}
\newtheorem{proposition}{Proposition}
\newtheorem{corollary}{Corollary}
\theoremstyle{definition}
\newtheorem{definition}{Definition}
\newtheorem{remark}{Remark}

\usepackage{calc} 
\usepackage{adjustbox} 

\newlength{\myheight}

\usepackage{filecontents}
\usepackage[english]{babel}

\usepackage[
backend=biber,
style=alphabetic,
citestyle=alphabetic,
doi=false,isbn=false,url=false,eprint=true
]{biblatex}

\newcommand{\R}{\mathbb{R}}

\newcommand{\C}{\mathbb{C}}

\newcommand{\Z}{\mathbb{Z}}

\newcommand{\ra}{\rightarrow}

\newcommand{\Ical}{\mathcal{I}}

\newcommand{\Ocal}{\mathcal{O}}
\newcommand{\Pcal}{\mathcal{P}}

\newcommand{\Hfrak}{\mathfrak{H}}

\newcommand{\bfe}{\mathbf{e}}

\newcommand{\wt}{\text{wt}}

\newcommand{\rbox}{\todo[inline,color=red!20]}

\makeatletter
\DeclareRobustCommand{\cev}[1]{%
  \mathpalette\do@cev{#1}%
}
\newcommand{\do@cev}[2]{%
  \fix@cev{#1}{+}%
  \reflectbox{$\m@th#1\vec{\reflectbox{$\fix@cev{#1}{-}\m@th#1#2\fix@cev{#1}{+}$}}$}%
  \fix@cev{#1}{-}%
}
\newcommand{\fix@cev}[2]{%
  \ifx#1\displaystyle
    \mkern#23mu
  \else
    \ifx#1\textstyle
      \mkern#23mu
    \else
      \ifx#1\scriptstyle
        \mkern#22mu
      \else
        \mkern#22mu
      \fi
    \fi
  \fi
}
\makeatother

\def\IntKern{}
\newcommand\MySymbolint[2][0pt]{%
  \mathchoice
  {\def\IntKern{\dimexpr-.6\wd0+#1\relax}\Symbolint\displaystyle\displaystyle{\scalebox{.8}{$#2$}}}%
  {\def\IntKern{\dimexpr-.62\wd0+#1\relax}\Symbolint\textstyle\textstyle{\scalebox{.6}{$#2$}}}%
  {\def\IntKern{\dimexpr-.642\wd0+#1\relax}\Symbolint\scriptstyle\scriptscriptstyle{\scalebox{.45}{$#2$}}}%
 {\def\IntKern{\dimexpr-.67\wd0+#1\relax}\Symbolint\scriptscriptstyle\scriptscriptstyle{\scalebox{.35}{$#2$}}}%
  \!\int}
\newcommand\Symbolint[3]{%
  {\setbox0=\hbox{$#1{#2#3}{\int}$ }%
  \vcenter{\hbox{$#2#3$ }}\kern\IntKern}
}
\newcommand{\Kint}{\MySymbolint{\textnormal{K}}}

\newcommand{\Gr}{\textnormal{Gr}}

\newcommand{\GL}{\textnormal{GL}}
\newcommand{\Fl}{F\ell}

\DeclareMathOperator{\codim}{codim}

\DeclareMathOperator{\sort}{sort}

\newcommand{\pt}{\text{pt}}

\newcommand{\zerostr}{\mathtt{0}}
\newcommand{\onestr}{\mathtt{1}}
\newcommand{\tenstr}{\mathtt{1\! 0}}
\newcommand{\twostr}{\mathtt{2}}
\newcommand{\threestr}{\mathtt{3}}
\newcommand{\dstr}{\mathtt{d}}
\newcommand{\mstr}{\mathtt{m}}

\newcommand{\pad}{^!}

\newcommand{\flag}[1]{#1_\bullet}

\newcommand{\Ftilde}{\Tilde{F}}

\newlength\lthk
\makeatletter
\newcommand\xleftrightarrow[2][]{%
  \ext@arrow 9999{\longleftrightarrowfill@}{#1}{#2}}
\newcommand\longleftrightarrowfill@{%
  \arrowfill@\leftarrow\relbar\rightarrow}
\makeatother

\newcommand{\inlinetri}{
\tikzset{every picture/.style={line width=0.75pt}} 
\begin{tikzpicture}[x=0.75pt,y=0.75pt,yscale=-1,xscale=1,scale=.65]
\raisebox{-.1ex}
{
\draw  [line width=0.75]  (18.6,6.68) -- (28.6,24) -- (8.6,24) -- cycle ;
}
\end{tikzpicture}
}

\newcommand{\inlinetrapside}{
\tikzset{every picture/.style={line width=0.75pt}} 
\begin{tikzpicture}[x=0.75pt,y=0.75pt,yscale=-1,xscale=1,scale=.75]
\raisebox{-.35ex}
{
\draw   (22.8,24) -- (12.8,24) -- (7.8,15.34) -- (12.8,6.68) -- cycle ;
}
\end{tikzpicture}
}

\newcommand{\inlinetrapsiderotate}{
\tikzset{every picture/.style={line width=0.75pt}} 
\begin{tikzpicture}[x=0.75pt,y=0.75pt,yscale=-1,xscale=1,scale=.75]
\raisebox{-.5ex}
{
\draw   (11.8,10.68) -- (21.8,10.68) -- (26.8,19.34) -- (21.8,28) -- cycle ;
}
\end{tikzpicture}}

\newcommand{\inlinetrapsiderotatetwo}{
\tikzset{every picture/.style={line width=0.75pt}} 
\begin{tikzpicture}[x=0.75pt,y=0.75pt,yscale=-1,xscale=1,scale=.75]
\raisebox{-.35ex}
{
\draw   (28.05,9.35) -- (33.05,18.01) -- (28.05,26.67) -- (18.05,26.67) -- cycle ;
}
\end{tikzpicture}
}

\newcommand{\inlinepar}{
\tikzset{every picture/.style={line width=0.75pt}} 
\begin{tikzpicture}[x=0.75pt,y=0.75pt,yscale=-1,xscale=1,scale=.65]
\raisebox{-.7ex}
{
\draw  [line width=0.75]  (14.4,11.02) -- (24.4,28.34) -- (19.4,37) -- (9.4,19.68) -- cycle ;
}
\end{tikzpicture}
}

\newcommand{\inlinerhombus}{
\tikzset{every picture/.style={line width=0.75pt}} 
\begin{tikzpicture}[x=0.75pt,y=0.75pt,yscale=-1,xscale=1,scale=.5]
\raisebox{-.7ex}
{
\draw   (23.33,8.96) -- (33.33,26.28) -- (23.33,43.6) -- (13.33,26.28) -- cycle ;
}
\end{tikzpicture}
}

\newcommand{\inlinepent}{
\tikzset{every picture/.style={line width=0.75pt}} 
\begin{tikzpicture}[x=0.75pt,y=0.75pt,yscale=-1,xscale=1,scale=.55]
\raisebox{-.5ex}
{
\draw   (22,8.02) -- (32,25.34) -- (27,34) -- (17,34) -- (12,25.34) -- cycle ;
}
\end{tikzpicture}}

\newcommand{\inlinepentside}{
\tikzset{every picture/.style={line width=0.75pt}} 
\begin{tikzpicture}[x=0.75pt,y=0.75pt,yscale=-1,xscale=1,scale=.6]
\raisebox{-.3ex}
{
\draw   (32.75,35.5) -- (42.75,18.18) -- (52.75,18.18) -- (57.75,26.84) -- (52.75,35.5) -- cycle ;
}
\end{tikzpicture}}

\newcommand{\inlinehex}{
\tikzset{every picture/.style={line width=0.75pt}} 
\begin{tikzpicture}[x=0.75pt,y=0.75pt,yscale=-1,xscale=1,scale=.65]
\raisebox{-.3ex}
{
\draw  [line width=0.75]  (24,9.68) -- (29,18.34) -- (24,27) -- (14,27) -- (9,18.34) -- (14,9.68) -- cycle ;
}
\end{tikzpicture}
}

\title{Commutative Properties of Schubert Puzzles with Convex Polygonal Boundary Shapes}
\author{Portia X. Anderson}
\date{}

\begin{document}

\maketitle

\begin{abstract}

We generalize classical triangular Schubert puzzles to puzzles with convex polygonal boundary. We give these puzzles a geometric Schubert calculus interpretation and derive novel combinatorial commutativity statements, using purely geometric arguments, for puzzles with four, five, and six sides, having various types of symmetry in their boundary conditions. We also present formulas for the associated structure constants in terms of Littlewood-Richardson numbers, and we prove an analogue of commutativity for parallelogram-shaped equivariant puzzles.

\end{abstract}

\setcounter{tocdepth}{2}
\tableofcontents

\section{Introduction}

Schubert puzzles were first introduced in 1999 in \cite{KTW}, and there has been continuous research activity on them since then.
However, rather little work has been done exploring non-triangular puzzles up until now. (As far as we are aware, only rhombus-shaped puzzles have received some attention in the literature, e.g. in \cite{KTpuzzles} and \cite{Purbhoo}.) In this paper, we make an effort to give a comprehensive survey of the commutative properties of generalized polygonal puzzles of the type that compute structure constants in the ordinary cohomology of the Grassmannian. These properties will generalize the fundamental commutative property of classical triangular puzzles.

Our commutativity results were originally discovered and are proved here purely through geometric means, though they are primarily combinatorial statements. We describe a way to give a classical Schubert calculus interpretation to a puzzle with convex polygonal boundary, and we perform our computations in a geometric and cohomological context. As another result of this process, we will present formulas for the associated structure constants in terms of Littlewood-Richardson numbers.

The paper will proceed by proving statements for simpler boundary shapes first, and then it will progressively build upon those and specialize them further as we increase the number of sides of our boundaries and reach more symmetric cases.

We prove an analogue of commutativity for parallelogram-shaped equivariant puzzles, and we include a number of remarks about how our results extend to puzzles used for K-theory, and to 2-step and 3-step puzzles. In the future, we would like to expand our understanding to encompass other boundary shapes for equivariant puzzles, as well as polygonal puzzles using other extended puzzle piece sets purposed for Schubert calculus in other cohomology theories.

\subsection{Definition of puzzles}

\subsubsection{Classical triangular puzzles}
Classically, a \textbf{Schubert puzzle} is a tiling of an equilateral triangular region, whose boundary we denote $\inlinetri$, using a set of labeled unit triangles 

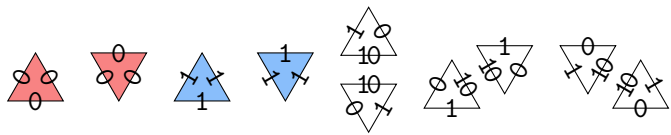
\begin{center}

\begin{tikzpicture}[x=1pt,y=1pt,yscale=1,xscale=1,scale=.7]
\draw [fill={rgb, 255:red, 248; green, 131; blue, 131 }] (0.0,25.98) -- (15.0,51.96) -- (30.0,25.98) -- cycle;
\draw (7.5,38.97) node  [rotate=60]  {$\mathtt{0}$};
\draw (22.5,38.97) node  [rotate=-60]  {$\mathtt{0}$};
\draw (15.0,25.98) node  [rotate=0]  {$\mathtt{0}$};
\draw [fill={rgb, 255:red, 248; green, 131; blue, 131 }] (75.0,51.96) -- (60.0,25.98) -- (45.0,51.96) -- cycle;
\draw (67.5,38.97) node  [rotate=60]  {$\mathtt{0}$};
\draw (52.5,38.97) node  [rotate=-60]  {$\mathtt{0}$};
\draw (60.0,51.96) node  [rotate=0]  {$\mathtt{0}$};
\draw [fill={rgb, 255:red, 137; green, 190; blue, 251 }] (90.0,25.98) -- (105.0,51.96) -- (120.0,25.98) -- cycle;
\draw (97.5,38.97) node  [rotate=60]  {$\mathtt{1}$};
\draw (112.5,38.97) node  [rotate=-60]  {$\mathtt{1}$};
\draw (105.0,25.98) node  [rotate=0]  {$\mathtt{1}$};
\draw [fill={rgb, 255:red, 137; green, 190; blue, 251 }] (165.0,51.96) -- (150.0,25.98) -- (135.0,51.96) -- cycle;
\draw (157.5,38.97) node  [rotate=60]  {$\mathtt{1}$};
\draw (142.5,38.97) node  [rotate=-60]  {$\mathtt{1}$};
\draw (150.0,51.96) node  [rotate=0]  {$\mathtt{1}$};
\draw [fill={rgb, 255:red, 255; green, 255; blue, 255 }] (180.0,50.66100000000001) -- (195.0,76.641) -- (210.0,50.66100000000001) -- cycle;
\draw (187.5,63.65100000000001) node  [rotate=60]  {$\mathtt{1}$};
\draw (202.5,63.65100000000001) node  [rotate=-60]  {$\mathtt{0}$};
\draw (195.0,50.66100000000001) node  [rotate=0]  {$\mathtt{1\! 0}$};
\draw [fill={rgb, 255:red, 255; green, 255; blue, 255 }] (210.0,35.073) -- (195.0,9.093) -- (180.0,35.073) -- cycle;
\draw (202.5,22.083) node  [rotate=60]  {$\mathtt{1}$};
\draw (187.5,22.083) node  [rotate=-60]  {$\mathtt{0}$};
\draw (195.0,35.073) node  [rotate=0]  {$\mathtt{1\! 0}$};
\draw [fill={rgb, 255:red, 255; green, 255; blue, 255 }] (225.0,22.083000000000002) -- (240.0,48.063) -- (255.0,22.083000000000002) -- cycle;
\draw (232.5,35.073) node  [rotate=60]  {$\mathtt{0}$};
\draw (247.5,35.073) node  [rotate=-60]  {$\mathtt{1\! 0}$};
\draw (240.0,22.083000000000002) node  [rotate=0]  {$\mathtt{1}$};
\draw [fill={rgb, 255:red, 255; green, 255; blue, 255 }] (283.49960399419183,55.857) -- (268.49960399419183,29.877) -- (253.49960399419183,55.857) -- cycle;
\draw (275.99960399419183,42.867) node  [rotate=60]  {$\mathtt{0}$};
\draw (260.99960399419183,42.867) node  [rotate=-60]  {$\mathtt{1\! 0}$};
\draw (268.49960399419183,55.857) node  [rotate=0]  {$\mathtt{1}$};
\draw [fill={rgb, 255:red, 255; green, 255; blue, 255 }] (326.99920798838366,22.083000000000002) -- (341.99920798838366,48.063) -- (356.99920798838366,22.083000000000002) -- cycle;
\draw (334.49920798838366,35.073) node  [rotate=60]  {$\mathtt{1\! 0}$};
\draw (349.49920798838366,35.073) node  [rotate=-60]  {$\mathtt{1}$};
\draw (341.99920798838366,22.083000000000002) node  [rotate=0]  {$\mathtt{0}$};
\draw [fill={rgb, 255:red, 255; green, 255; blue, 255 }] (328.49960399419183,55.857) -- (313.49960399419183,29.877) -- (298.49960399419183,55.857) -- cycle;
\draw (320.99960399419183,42.867) node  [rotate=60]  {$\mathtt{1\! 0}$};
\draw (305.99960399419183,42.867) node  [rotate=-60]  {$\mathtt{1}$};
\draw (313.49960399419183,55.857) node  [rotate=0]  {$\mathtt{0}$};
\end{tikzpicture}

\end{center}

called \textbf{puzzle pieces}, 
so that any two glued edges have the same label, and only $\zerostr$ and $\onestr$ labels appear along the outer boundary $\inlinetri$, not $\tenstr$s.

Throughout this paper, we will let $\Hfrak$ denote this base set of puzzle pieces. To be more explicit, $\Hfrak$ consists of 

\begin{itemize}
    \item an upward-pointing unit triangle whose edges are all labeled $\zerostr$,
    \item an upward-pointing unit triangle whose edges are all labeled $\onestr$,
    \item an upward-pointing unit triangle whose NW, NE, and South edges are labeled $\onestr$, $\zerostr$, and $\tenstr$, respectively,
\end{itemize}
along with all their rotations in $60^\circ$ increments. 

\vspace{.5\baselineskip}
Often, we label the boundary ahead of time with binary strings of $\zerostr$s and $\onestr$s with the intent of finding all puzzles that agree with that fixed boundary labeling. An equilateral triangular boundary whose NW, NE, and South sides are labeled with strings $\lambda$, $\mu$, and $\nu$ in the orientations shown in Figure \ref{fig:lambdamunuboundary} will be denoted $\inlinetri_{\lambda,\mu,\nu}$, and a puzzle with this boundary labeling will be called a \textbf{$\inlinetri_{\lambda,\mu,\nu}$-puzzle}. See Figure \ref{fig:1010_0101_0011} for an example.

\begin{figure}[h]
\centering
    \begin{subfigure}[t]{0.4\textwidth}
        \centering
        \input{examples/lambdamunuboundary}
        \caption{The labeled boundary $\protect\inlinetri_{\lambda,\mu,\nu}$}
        \label{fig:lambdamunuboundary}
    \end{subfigure}
    ~
    \begin{subfigure}[t]{0.4\textwidth}
        \centering
        \input{examples/1010_0101_0011}
        
        \caption{A $\protect\inlinetri_{\onestr\zerostr\onestr\zerostr,\zerostr\onestr\zerostr\onestr,\zerostr\zerostr\onestr\onestr}$-puzzle}
        \label{fig:1010_0101_0011}
    \end{subfigure}
    \caption{}
    \label{fig:triangularboundaryandpuzzle}
\end{figure}

\subsubsection{Puzzles with convex polygonal boundary}

Now in this paper, we will generalize the definition of ``puzzles'' to include puzzle piece tilings of convex polygonal shapes, where
\begin{enumerate}
    \item the angle of each puzzle piece edge relative to the $x$-axis is a multiple of $60^\circ$, and 
    \item only $\zerostr$ and $\onestr$ labels are allowed to appear along the outer boundary, not $\tenstr$s.
\end{enumerate}

This allows us to have puzzles with trapezoidal, parallelogram-shaped, pentagonal, and hexagonal boundary as well. In fact, a hexagon generalizes all these shapes; a boundary with three, four, or five sides is just a hexagon where some sides have length 0.

\begin{figure}[h]
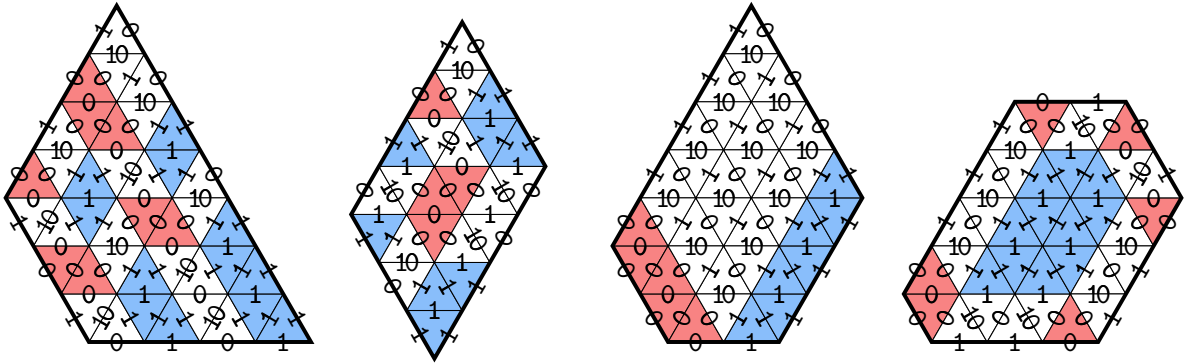

\centering
\begin{subfigure}[t]{0.25\textwidth}
\centering
\input{examples/trapezoidexample}
\caption{An example of a $ \protect\inlinetrapside_{\onestr\zerostr\onestr,\zerostr\onestr\zerostr\onestr,\zerostr\zerostr\onestr\zerostr\onestr\onestr\onestr,\onestr\zerostr\onestr\zerostr}$-puzzle}
\label{fig:trapezoidexample}
\end{subfigure}
\begin{subfigure}[t]{0.215\textwidth}
\centering
\input{examples/parallelogramexample}
\caption{An example of a $ \protect\inlinepar_{\onestr\zerostr\onestr,\zerostr\onestr\zerostr\onestr,\zerostr\onestr\onestr,\zerostr\zerostr\onestr\onestr}$-puzzle}
\label{fig:parallelogramexample}
\end{subfigure}
\begin{subfigure}[t]{0.235\textwidth}
\centering
\input{examples/pentagonexample}
\caption{An example of a $ \protect\inlinepent_{\zerostr\zerostr,\zerostr\onestr\onestr\onestr\onestr,\zerostr\zerostr\zerostr\onestr,\onestr\onestr\onestr,\onestr\zerostr}$-puzzle}
\label{fig:pentagonexample}
\end{subfigure}
\begin{subfigure}[t]{0.24\textwidth}
\centering
\input{examples/hexagonexample}
\caption{An example of a $ \protect\inlinehex_{\zerostr,\zerostr\onestr\onestr\onestr,\zerostr\onestr,\zerostr\onestr,\zerostr\onestr\onestr,\zerostr\onestr\onestr}$-puzzle}
\label{fig:hexagonexample}
\end{subfigure}

\centering
\caption{Examples of puzzles with different boundary shapes. The boundary labels are read clockwise starting from the SW side.}
\label{fig:boundaryshapespuzzles}
\end{figure}

Similarly to $\inlinetri_{\lambda,\mu,\nu}$, we will use a shape symbol (where only the number of sides and angles matter, not side lengths) with a subscript sequence of strings (which we usually represent with lowercase Greek letters) to denote a labeled boundary, for example  $\inlinetrapside_{\lambda,\mu,\nu,\xi}$, $\inlinepar_{\lambda,\mu,\nu,\xi}$, $\inlinepent_{\lambda,\mu,\nu,\xi,o}$, or $\inlinehex_{\lambda,\mu,\nu,\xi,o,\pi}$. 
For any shape, we will always begin at the SW side of the boundary and read the sequence of labels, assigning one to each side of the boundary, in a clockwise orientation. 
We will also use the suffix \textbf{-puzzle} to refer to a puzzle with the indicated boundary (see Figure \ref{fig:boundaryshapespuzzles}).

\subsection{Schubert calculus}

Classically, \textbf{Schubert calculus} refers to the activity of computing the structure constants in the Schubert bases for the cohomology ring of the Grasmmannian, but it also generalizes to other cohomology theories, such as $T$-equivariant cohomology or K-theory, and to $d$-step flag manifolds. Our reference for Schubert calculus in ordinary cohomology will be \cite{Fulton}.

In this paper, we will focus mainly on ordinary cohomology of the Grassmannian, though with one theorem about $T$-equivariant cohomology, as well as a number of remarks about K-theory and $d$-step Schubert calculus.

\subsubsection{Schubert varieties}
\label{section:schubertvarieties}
Let $\Gr(k;\C^n)$ denote the \textbf{Grassmannian} of $k$-dimensional subspaces of $\C^n$. It is a smooth complex manifold and projective algebraic variety of dimension $k(n-k)$.

We will write $\lambda\in\binom{[n]}{k}$ to mean $\lambda$ is a \textbf{binary string} with $k$ $\onestr$s and $n-k$ $\zerostr$s, and we will use $\lambda_i$ to denote the $i$th digit of $\lambda$. Let $\lambda=\lambda_1\lambda_2\cdots\lambda_n\in\binom{[n]}{k}$, and let $F_\bullet$ be a complete flag $\{\mathbf{0}\}=F_0\subset F_1\subset\cdots\subset F_n=\C^n$ of nested subspaces whose dimensions are indicated by their subscripts. Then $X_\lambda(F_\bullet)$ is the \textbf{Schubert variety} in $\Gr(k;\C^n)$ defined by 
$$X_\lambda(F_\bullet):= \{V\in\Gr(k;\C^n) : \dim(V\cap F_i)\geq \lambda_1+\lambda_2+\cdots+\lambda_i\}.$$

We use the shorter notation $X_\lambda:=X_\lambda(\flag{F})$ specifically when $\flag{F}$ is the \textit{standard complete flag} $\{\mathbf{0}\}\subset \langle \bfe_1\rangle\subset \langle \bfe_1,\bfe_2\rangle \subset \cdots\subset \langle \bfe_1,\bfe_2,\ldots,\bfe_n\rangle=\C^n$, where $\bfe_1,\bfe_2,\ldots,\bfe_n$ are the standard basis vectors. We will use the notation $^\vee$ throughout the paper to denote the operation of reversing a string. Then we use the notation $X^\lambda:=X_{\lambda^\vee}(\flag{\Ftilde})$ specifically when $\flag{\Ftilde}$ is the \textit{anti-standard complete flag} $\{\mathbf{0}\}\subset \langle \bfe_n\rangle\subset \langle \bfe_n,\bfe_{n-1}\rangle \subset \cdots\subset \langle \bfe_n,\bfe_{n-1},\ldots,\bfe_1\rangle=\C^n$. $X^\lambda$ is called an \textbf{opposite Schubert variety}, and this is indicated by the superscript instead of a subscript.

An \textbf{inversion} in $\lambda$ is a pair $i<j$ such that $\onestr=\lambda_i>\lambda_j=\zerostr$. Let $\ell(\lambda)$ denote the \textbf{length} of $\lambda$, which is defined as the number of inversions in $\lambda$. Then $\codim X_\lambda(F_\bullet)=\ell(\lambda)$. 

Each Schubert variety $X_\lambda(\flag{F})$ defines a homology class $[X_\lambda(\flag{F})]\in H_{2(k(n-k)-\ell(\lambda))}(\Gr(k;\C^n))$. $\GL_n(\C)$ acts transitively on the complete flag manifold $\Fl(\C^n)$. So for any pair $F_\bullet,G_\bullet\in\Fl(\C^n)$, there exists some $g\in\GL_n(\C)$ such that $F_\bullet=gG_\bullet$. $\GL_n(\C)$ is path-connected, so the inclusion $X_\lambda(F_\bullet)\hookrightarrow \Gr(k;\C^n)$ is homotopic to inclusion followed by the left multiplication by $g$. Thus we have $[X_\lambda(F_\bullet)]=g_*[X_\lambda(F_\bullet)]=g\cdot [X_\lambda(F_\bullet)] = [g\cdot X_\lambda(F_\bullet)]=[X_\lambda(g\cdot F_\bullet)]= [X_\lambda(G_\bullet)]$. So in fact, $[X_\lambda(F_\bullet)]=[X_\lambda]$, regardless of what flag $\flag{F}$ is.

Since $\Gr(k;\C^n)$ is an oriented closed manifold, it has Poincar\'e duality. From now on in this paper, we will use the notation $[X_\lambda(F_\bullet)]$ to refer interchangeably to the homology class or to its Poincar\'e dual class in the cohomology $H^*(\Gr(k;\C^n))$, and for the most part we will be talking about cohomology. 

The Schubert classes $\{[X_\lambda] : \lambda \in \binom{[n]}{k}\}$ form a $\Z$-basis for $H^*(\Gr(k;\C^n))$. So we can expand a product $[X_\lambda][X_\mu]$ uniquely as a linear combination
$$[X_\lambda][X_\mu]=\sum_{\nu\in \binom{[n]}{k}}c_{\lambda,\mu}^\nu [X_\nu].$$

The structure constants $c_{\lambda,\mu}^\nu$ are positive, and in fact they are the \textit{Littlewood-Richardson numbers} that arise in symmetric function theory.

The opposite Schubert classes $\{[X^\lambda] : \lambda\in\binom{[n]}{k}\}$ form the dual basis under the perfect pairing 
$$\underset{\Gr(k;\C^n)}{\int} [X_\lambda][X^\mu] = \langle \lambda,\mu\rangle :=\begin{cases}1, &\text{ if } \lambda=\mu\\0, &\text{ if }\lambda\neq\mu \end{cases}.$$

\subsubsection{Intersection theory}
\label{section:intersectiontheory}

Using the dual basis, we can find the structure constants as 
$$c_{\lambda,\mu}^\nu= \underset{\Gr(k;\C^n)}{\int} [X_\lambda][X_\mu][X^\nu].$$
The integral can be thought of as ``pushforward to a point.''

A fact that we use often is that the cup product in cohomology is Poincar\'e dual to intersection in homology. If two subvarieties $Y$ and $Z$ intersect transversely, then we have $[Y][Z]=[Y\cap Z]$. Schubert varieties intersect transversely if the flags they are defined with respect to are in general position relative to each other. This is a consequence of Kleiman's transversality theorem in \cite{Kleiman}, which implies that general translates of subvarieties in a homogeneous space intersect transversely. 
\explainalot{By choosing flags $\flag{F},\flag{G},\flag{H}$ in general position, we get that $$c_{\lambda,\mu}^\nu= \underset{\Gr(k;\C^n)}{\int} [X_\lambda(\flag{F})][X_\mu(\flag{G})][X_\nu(\flag{H})] =\underset{\Gr(k;\C^n)}{\int} [X_\lambda(\flag{F})\cap X_\mu(\flag{G})\cap X_\nu(\flag{H})]$$ equals a count of the intersection points when it is finite and 0 otherwise.
\rbox{only intersect properly if the subvarieties aren't smooth. do I need Buch's thing bc schubert varieties aren't smooth so this isn't enough?} 
This is worth mentioning, but it will not be used in this paper.}
The standard and anti-standard flags are in general position relative to each other, so $X_\lambda$ and $X^\mu$ always intersect transversely, and $[X_\lambda][X^\mu]=[X_\lambda\cap X^\mu]$. 

When considering the conditions encoded by $\lambda$ when defining $X_\lambda$, we call the condition $\dim(V\cap F_i)\geq \lambda_1+\lambda_2+\cdots+\lambda_i$ an \textbf{essential condition} whenever $\lambda_i= \onestr$ and $\lambda_{i+1}=\zerostr$. The set of essential conditions completely determines $X_\lambda$. 

Throughout the paper, we will let \textbf{sort} denote the operation of moving all the $\zerostr$s in a binary string ahead of all the $\onestr$s. So for $\lambda\in\binom{[n]}{k}$, we have $\sort(\lambda)=\zerostr^{n-k}\onestr^{k}$. (Note that $\sort(\lambda)$ indicates the \textbf{content} of $\lambda$, i.e. the number of $\zerostr$s and $\onestr$s it contains.) Then, for example, we have $X_{\sort(\lambda)}= \Gr(k;\C^n)$ because the string $\zerostr^{n-k}\onestr^{k}$ encodes no essential conditions, and we have $[X_{\sort(\lambda)}]=[\Gr(k;\C^n)]=1\in H^*(\Gr(k;\C^n)$.



\subsubsection{Equivariant cohomology}
Our reference for this section is \cite{AndersonFulton}. Let $T=(\C^\times)^n$, an $n$-dimensional torus, and let $ET$ be a contractible space with a free $T$-action. The \textbf{$T$-equivariant cohomology} of the Grassmannian is defined as $H^*_T(\Gr(k;\C^n)):=H^*(ET\times^T\Gr(k;\C^n))$, and it is an algebra over $H^*_T(\pt)=\Z[y_1,\ldots,y_n]$. Essential facts for this paper are that the Schubert varieties $X_\lambda$ and $X^\lambda$ are $T$-invariant, and $\{[X_\lambda] : \lambda \in \binom{[n]}{k}\}$ and $\{[X^\lambda] : \lambda\in \binom{[n]}{k}\}$ form dual bases for $H^*_T(\Gr(k;\C^n))$. (We will use the same notation $[X_\lambda]$ to denote either the ordinary or equivariant cohomology class, depending on the context.) There is also a surjective map $H^*_T(\Gr(k;\C^n)) \twoheadrightarrow H^*(\Gr(k;\C^n))$ sending $y_1,\ldots,y_n\mapsto 0$ and $[X_\lambda]\mapsto [X_\lambda]$, hence the structure constants for $H^*_T(\Gr(k;\C^n))$ (which we will denote $(c_T)_{\lambda,\mu}^\nu$) agree with those for $H^*(\Gr(k;\C^n))$ in degree 0.

 Be warned that in equivariant cohomology, the classes $[X_\lambda(\flag{F})]$ are \textit{not} independent of the choice of flag, and the Schubert varieties cannot simply be perturbed to lie in general position relative to each other while remaining $T$-invariant. However, we still have $[X_\lambda][X^\mu]=[X_\lambda\cap X^\mu]$, which our proofs involving equivariant cohomology will employ.

\subsubsection{K-theory}

Our discussion of K-theory in this paper will be rather limited and will mostly be confined to several remarks. The only necessary background here is that the classes of the structure sheaves of the Schubert varieties form a $\Z$-basis $\{[\Ocal_\lambda]: \lambda\in\binom{[n]}{k}\}$ for the K-theory of the Grassmannian $K(\Gr(k;\C^n))$, and dual to that basis is what is referred to as the ideal sheaf basis and denoted $\{[\Ical^\lambda]: \lambda\in\binom{[n]}{k}\}$.

\subsubsection{Schubert calculus on $d$-step flag manifolds}
\label{section:dstepschubertcalculus}


Let $\Fl(k_1,k_2,\ldots,k_d;\C^{n})$ denote the space (a smooth complex manifold and projective variety) of $d$-step flags $\{\mathbf{0}\} \subset V_1\subset V_2\subset \cdots\subset V_d\subset \C^n$, where $\dim(V_j)=k_j$ for all $1\leq j\leq d$. 

Now the set of strings $\lambda$ with the content of $\zerostr^{n-k_d}\onestr^{k_d-k_{d-1}}\twostr^{k_{d-1}-k_{d-2}}\cdots\dstr^{k_1}$ index the Schubert varieties $X_\lambda(F_\bullet)$ in $\Fl(k_1,k_2,\ldots,k_d;\C^{n})$, which are defined by 
$$X_\lambda(F_\bullet):= \left\{\flag{V} \in \Fl(k_1,k_2,\ldots,k_d;\C^{n}) : \sum_{j=1}^d \dim(V_j\cap F_i)\geq \lambda_1+\lambda_2+\cdots+\lambda_i\right\},$$
with respect to a complete flag $\flag{F}$.
The rest of the discussion of Schubert calculus on the Grassmannian (a 1-step flag manifold), following the definition Schubert varieties, extends almost identically to the more general case of $\Fl(k_1,k_2,\ldots,k_d;\C^{n})$. 

There are a couple of other general facts that we will need for this paper. Define the map $$p^i:\Fl(k_1,k_2,\ldots,k_d;\C^{n})\twoheadrightarrow \Gr(k_i;\C^{n}),\quad  \{\mathbf{0}\} \subset V_1\subset V_2\subset \cdots\subset V_d\subset \C^n\mapsto V_i,$$ 
i.e. projection onto the $i$th component.
Then $\flag{V}\in X_\lambda$ if and only if $p^i(\flag{V})\in p^i(X_\lambda)$ for all $1\leq i\leq d$. We also have that $p^i(X_\lambda)=X_{\lambda^i}\subseteq \Gr(k_i;\C^n)$, where $\lambda^i$ is the binary string obtained from $\lambda$ by changing each occurrence of $\mstr$ to a $\onestr$ for all $m > d-i$ and changing each occurrence of $\mstr$ to a $\zerostr$ for all $m\leq d-i$.




\subsection{Puzzles compute Schubert calculus}


\subsubsection{Puzzles using the base set of pieces}

It is a result of Knutson--Tao--Woodward (see \cite[Theorem 1]{KTW}) that puzzles using the base set of pieces $\Hfrak$ compute the structure constants of $H^*(\Gr(k;\C^n))$ in the Schubert basis:
$$\underset{\Gr(k;\C^n)}{\int} [X_\lambda][X_\mu][X^\nu] = \#\{\inlinetri_{\lambda,\mu,\nu^\vee}\text{-puzzles}\}.$$


\subsubsection{Extended puzzle piece sets for other cohomology theories}
\label{section:puzzlepiecesets}

There are other puzzle piece sets, comprised of the base set $\Hfrak$ plus additional special pieces, which compute structure constants in other cohomology theories for the Grassmannian in the Schubert bases. The special pieces carry a special \textbf{weight}, while the pieces in $\Hfrak$ can be seen as having weight 1. Then the \textbf{weight} of a puzzle is defined as the product of the weights of its pieces, and we find the structure constant by summing over the weights of all puzzles with the corresponding boundary conditions. 
We list the various puzzle piece sets below (sets 1-5 will come up throughout the paper, and we leave the others for future questions).

\begin{enumerate}
    \item $\Hfrak$ -- computes the structure constants for ordinary cohomology $H^*(\Gr(k;\C^n)$. 
    \item $\Hfrak \cup \left\{\raisebox{-.4\height}{
\begin{tikzpicture}[x=1pt,y=1pt,yscale=1,xscale=1,scale=.68]
\draw [fill={rgb, 255:red, 255; green, 246; blue, 74 }] (0.0,0.0) -- (15.0,25.98) -- (30.0,0.0) -- (15.0,-25.98) -- cycle;
\draw (7.5,12.99) node  [rotate=60]  {$\mathtt{0}$};
\draw (22.5,12.99) node  [rotate=-60]  {$\mathtt{1}$};
\draw (22.5,-12.99) node  [rotate=60]  {$\mathtt{0}$};
\draw (7.5,-12.99) node  [rotate=-60]  {$\mathtt{1}$};
\end{tikzpicture}
} \right\}$ -- computes the structure constants for $T$-equivariant cohomology $H_T^*(\Gr(k;\C^n))$ (see \cite{KTpuzzles}). We call the this special piece the \textit{equivariant piece}.
    \item $\Hfrak \cup \left\{\raisebox{-.4\height}{
\begin{tikzpicture}[x=1pt,y=1pt,yscale=1,xscale=1,scale=.68]
\draw [fill={rgb, 255:red, 255; green, 255; blue, 255 }] (0.0,0.0) -- (15.0,25.98) -- (30.0,0.0) -- cycle;
\draw (7.5,12.99) node  [rotate=60]  {$\mathtt{1\! 0}$};
\draw (22.5,12.99) node  [rotate=-60]  {$\mathtt{1\! 0}$};
\draw (15.0,0.0) node  [rotate=0]  {$\mathtt{1\! 0}$};
\end{tikzpicture}
} \right\}$ -- computes the structure constants for K-theory $K(\Gr(k;\C^n))$ in the Schubert structure sheaf basis $\{[\Ocal_\lambda]: \lambda\in\binom{[n]}{k}\}$ (see \cite[Theorem 4.6]{VakilKtheory}).
    \item $\Hfrak \cup \left\{\raisebox{-.4\height}{
\begin{tikzpicture}[x=1pt,y=1pt,yscale=1,xscale=1,scale=.68]
\draw [fill={rgb, 255:red, 255; green, 255; blue, 255 }] (30.0,25.98) -- (15.0,0.0) -- (0.0,25.98) -- cycle;
\draw (22.5,12.99) node  [rotate=60]  {$\mathtt{1\! 0}$};
\draw (7.5,12.99) node  [rotate=-60]  {$\mathtt{1\! 0}$};
\draw (15.0,25.98) node  [rotate=0]  {$\mathtt{1\! 0}$};
\end{tikzpicture}
} \right\}$ -- computes the structure constants for K-theory $K(\Gr(k;\C^n))$ in the Schubert ideal sheaf basis $\{[\Ical^\lambda]: \lambda\in\binom{[n]}{k}\}$ (see \cite[Theorem 1']{WheelerZJ}).
    \item $\Hfrak \cup \left\{\raisebox{-.4\height}{
\begin{tikzpicture}[x=1pt,y=1pt,yscale=1,xscale=1,scale=.68]
\draw [fill={rgb, 255:red, 255; green, 246; blue, 74 }] (0.0,0.0) -- (15.0,25.98) -- (30.0,0.0) -- (15.0,-25.98) -- cycle;
\draw (7.5,12.99) node  [rotate=60]  {$\mathtt{0}$};
\draw (22.5,12.99) node  [rotate=-60]  {$\mathtt{1}$};
\draw (22.5,-12.99) node  [rotate=60]  {$\mathtt{0}$};
\draw (7.5,-12.99) node  [rotate=-60]  {$\mathtt{1}$};
\end{tikzpicture}
}, \raisebox{-.4\height}{
\begin{tikzpicture}[x=1pt,y=1pt,yscale=1,xscale=1,scale=.68]
\draw [fill={rgb, 255:red, 255; green, 255; blue, 255 }] (0.0,0.0) -- (15.0,25.98) -- (30.0,0.0) -- cycle;
\draw (7.5,12.99) node  [rotate=60]  {$\mathtt{1\! 0}$};
\draw (22.5,12.99) node  [rotate=-60]  {$\mathtt{1\! 0}$};
\draw (15.0,0.0) node  [rotate=0]  {$\mathtt{1\! 0}$};
\end{tikzpicture}
} \right\}$ -- computes the structure constants for $T$-equivariant K-theory $K_T(\Gr(k;\C^n))$ in the Schubert structure sheaf basis $\{[\Ocal_\lambda]: \lambda\in\binom{[n]}{k}\}$ (see \cite[Theorem 2]{WheelerZJ}).
    \item $\Hfrak \cup \left\{\raisebox{-.4\height}{
\begin{tikzpicture}[x=1pt,y=1pt,yscale=1,xscale=1,scale=.68]
\draw [fill={rgb, 255:red, 255; green, 246; blue, 74 }] (0.0,0.0) -- (15.0,25.98) -- (30.0,0.0) -- (15.0,-25.98) -- cycle;
\draw (7.5,12.99) node  [rotate=60]  {$\mathtt{0}$};
\draw (22.5,12.99) node  [rotate=-60]  {$\mathtt{1}$};
\draw (22.5,-12.99) node  [rotate=60]  {$\mathtt{0}$};
\draw (7.5,-12.99) node  [rotate=-60]  {$\mathtt{1}$};
\end{tikzpicture}
}, \raisebox{-.4\height}{
\begin{tikzpicture}[x=1pt,y=1pt,yscale=1,xscale=1,scale=.68]
\draw [fill={rgb, 255:red, 255; green, 255; blue, 255 }] (30.0,25.98) -- (15.0,0.0) -- (0.0,25.98) -- cycle;
\draw (22.5,12.99) node  [rotate=60]  {$\mathtt{1\! 0}$};
\draw (7.5,12.99) node  [rotate=-60]  {$\mathtt{1\! 0}$};
\draw (15.0,25.98) node  [rotate=0]  {$\mathtt{1\! 0}$};
\end{tikzpicture}
} \right\}$ -- computes the structure constants for $T$-equivariant K-theory $K_T(\Gr(k;\C^n))$ in the Schubert ideal sheaf basis $\{[\Ical^\lambda]: \lambda\in\binom{[n]}{k}\}$ (see \cite[Theorem 2']{WheelerZJ}).
    \item $\Hfrak \cup \left\{\raisebox{-.4\height}{
\begin{tikzpicture}[x=1pt,y=1pt,yscale=1,xscale=1,scale=.68]
\draw [fill={rgb, 255:red, 255; green, 255; blue, 255 }] (0.0,0.0) -- (15.0,25.98) -- (30.0,0.0) -- cycle;
\draw (7.5,12.99) node  [rotate=60]  {$\mathtt{1\! 0}$};
\draw (22.5,12.99) node  [rotate=-60]  {$\mathtt{1\! 0}$};
\draw (15.0,0.0) node  [rotate=0]  {$\mathtt{1\! 0}$};
\end{tikzpicture}
},\raisebox{-.4\height}{
\begin{tikzpicture}[x=1pt,y=1pt,yscale=1,xscale=1,scale=.68]
\draw [fill={rgb, 255:red, 255; green, 255; blue, 255 }] (30.0,25.98) -- (15.0,0.0) -- (0.0,25.98) -- cycle;
\draw (22.5,12.99) node  [rotate=60]  {$\mathtt{1\! 0}$};
\draw (7.5,12.99) node  [rotate=-60]  {$\mathtt{1\! 0}$};
\draw (15.0,25.98) node  [rotate=0]  {$\mathtt{1\! 0}$};
\end{tikzpicture}
} \right\}$ -- computes the structure constants in the Segre--Schwartz--MacPherson (SSM) class basis of $\text{frac}\,H_{T \times C^\times}(T^*(G/P))$ (see \cite[\S 5]{knutsonzinnjustintwo}). 
\end{enumerate}

Note that in these extended sets, there are no additional rotations of the pieces included. For example, the yellow equivariant piece is only ever an upward-pointing rhombus and is never rotated. 

The pieces \raisebox{-.4\height}{
\begin{tikzpicture}[x=1pt,y=1pt,yscale=1,xscale=1,scale=.68]
\draw [fill={rgb, 255:red, 255; green, 255; blue, 255 }] (0.0,0.0) -- (15.0,25.98) -- (30.0,0.0) -- cycle;
\draw (7.5,12.99) node  [rotate=60]  {$\mathtt{1\! 0}$};
\draw (22.5,12.99) node  [rotate=-60]  {$\mathtt{1\! 0}$};
\draw (15.0,0.0) node  [rotate=0]  {$\mathtt{1\! 0}$};
\end{tikzpicture}
} and \raisebox{-.4\height}{
\begin{tikzpicture}[x=1pt,y=1pt,yscale=1,xscale=1,scale=.68]
\draw [fill={rgb, 255:red, 255; green, 255; blue, 255 }] (30.0,25.98) -- (15.0,0.0) -- (0.0,25.98) -- cycle;
\draw (22.5,12.99) node  [rotate=60]  {$\mathtt{1\! 0}$};
\draw (7.5,12.99) node  [rotate=-60]  {$\mathtt{1\! 0}$};
\draw (15.0,25.98) node  [rotate=0]  {$\mathtt{1\! 0}$};
\end{tikzpicture}
} each carry a weight of $-1$. The weight of the equivariant piece will be described next.

\justthesis{
When we use different puzzle piece sets to compute structure constants, the special pieces each carry a \textbf{weight}. Instead of just \textit{counting} the puzzles with the given boundary conditions, for each puzzle $P$ we take the product of the weights of all its pieces, which gives us the \textbf{weight} of the puzzle, $\wt(P)$. Then we sum over the weights of all the puzzles to find the structure constant. In fact, the pieces in the base set $\Hfrak$ can be thought of as having weight $1$, so when we count the puzzles we are still summing over their weights. The pieces \raisebox{-.4\height}{
\begin{tikzpicture}[x=1pt,y=1pt,yscale=1,xscale=1,scale=.68]
\draw [fill={rgb, 255:red, 255; green, 255; blue, 255 }] (0.0,0.0) -- (15.0,25.98) -- (30.0,0.0) -- cycle;
\draw (7.5,12.99) node  [rotate=60]  {$\mathtt{1\! 0}$};
\draw (22.5,12.99) node  [rotate=-60]  {$\mathtt{1\! 0}$};
\draw (15.0,0.0) node  [rotate=0]  {$\mathtt{1\! 0}$};
\end{tikzpicture}
} and \raisebox{-.4\height}{
\begin{tikzpicture}[x=1pt,y=1pt,yscale=1,xscale=1,scale=.68]
\draw [fill={rgb, 255:red, 255; green, 255; blue, 255 }] (30.0,25.98) -- (15.0,0.0) -- (0.0,25.98) -- cycle;
\draw (22.5,12.99) node  [rotate=60]  {$\mathtt{1\! 0}$};
\draw (7.5,12.99) node  [rotate=-60]  {$\mathtt{1\! 0}$};
\draw (15.0,25.98) node  [rotate=0]  {$\mathtt{1\! 0}$};
\end{tikzpicture}
} each carry a weight of $-1$. The weight of the equivariant piece will be described next.
}

\subsubsection{Equivariant puzzles}
\label{section:equivariantpuzzles}

Equivariant puzzles will be relevant in this paper and are the subject of Theorem \ref{thm:eqvtparallelogram}, so we will include more detailed discussion of them here.

The weight of an equivariant piece \raisebox{-.4\height}{
\begin{tikzpicture}[x=1pt,y=1pt,yscale=1,xscale=1,scale=.68]
\draw [fill={rgb, 255:red, 255; green, 246; blue, 74 }] (0.0,0.0) -- (15.0,25.98) -- (30.0,0.0) -- (15.0,-25.98) -- cycle;
\draw (7.5,12.99) node  [rotate=60]  {$\mathtt{0}$};
\draw (22.5,12.99) node  [rotate=-60]  {$\mathtt{1}$};
\draw (22.5,-12.99) node  [rotate=60]  {$\mathtt{0}$};
\draw (7.5,-12.99) node  [rotate=-60]  {$\mathtt{1}$};
\end{tikzpicture}
} appearing in a puzzle is a polynomial of the form $y_j-y_i$ where $i<j$, and $(i,j)$ correspond uniquely to the piece's position in the puzzle.
To find $i$, we draw a line downward from the SW edge of the piece, parallel to the NW side of the puzzle's boundary, until it hits the bottom of the boundary at the $i$th place. Likewise, to find $j$, we draw a line downward from the SE edge of the piece, parallel to the NE side of the puzzle's boundary, until it hits the bottom of the boundary at the $j$th place. An example is shown in Figure \ref{fig:equivariantpuzzleweight}.

Then the structure constants in $H^*_T(\Gr(k;\C^n))$ are given by
$$(c_T)_{\lambda,\mu}^\nu = \sum_{\inlinetri_{\lambda,\mu,\nu^\vee}\text{-puzzles } P} \wt(P) = \sum_{\inlinetri_{\lambda,\mu,\nu^\vee}\text{-puzzles } P} \left( \prod_{\substack{\text{equivariant} \\ \text{pieces }p \text{ in } P }} \wt(p) \right),$$ 
where $\wt$ denotes the weight of a puzzle or puzzle piece.

\begin{figure}[h]
    \centering
    \begin{subfigure}[t]{0.47\textwidth}
    \centering
     \input{examples/equivariantpuzzleweightexample}
    \caption{The weights of the equivariant pieces, from left to right, are $y_2-y_1$, $y_5-y_2$, and $y_6-y_4$. The weight of the puzzle is $(y_2-y_1)(y_5-y_2)(y_6-y_4)$.}
    \label{fig:equivariantpuzzleweightexample}   
    \end{subfigure}
    ~
    \begin{subfigure}[t]{0.47\textwidth}
    \centering
    \input{examples/equivariantpuzzleweightdiagram}
    \caption{To find the weight $y_j-y_i$ of an equivariant piece, we draw lines downward to the left and right from its position, and the two places where they hit the bottom give $(i,j)$ for that piece.}
    \label{fig:equivariantpuzzleweightdiagram}
    \end{subfigure}
    \caption{}
    \label{fig:equivariantpuzzleweight}
\end{figure}

\subsubsection{$d$-step puzzle rules}
There are puzzle rules, using larger sets of pieces, for Schubert calculus on 2-step, 3-step, and 4-step flag manifolds, and none yet for $d$-step where $d\geq 5$. (See \cite{BSKPtwostep}, \cite{Buchtwostep}, and \cite[Theorem 2]{knutsonzinnjustin} for 2-step puzzle rules for $H^*$, $H^*_T$, and $K_T$, respectively. See \cite[Theorem 3]{knutsonzinnjustin} for 3-step puzzle rules for $H^*$ and $K$. See \cite[\S 5.2]{knutsonzinnjustintwo} for a 4-step puzzle rule for $H^*$, which, unlike the rules for $d\leq 3$ is not manifestly positive.) 


Similarly to the case of Grassmannian (1-step) puzzles, the strings appearing along the boundary of a $d$-step puzzle must index Schubert varieties, so they consist only of the digits $\zerostr, \onestr,\twostr,\ldots,\dstr$ and cannot contain any \textit{composite} numbers, i.e. single numbers made up of multiple digits (such as the $\tenstr$ puzzle piece edge label for Grassmannian puzzles).

\explainalot{
References:
\begin{enumerate}
    \item $H^*$ 2-step -- \cite{BSKPtwostep}
    \item $H^*_T$ 2-step -- \cite{Buchtwostep}
    \item $K_T$ 2-step -- \cite[Theorem 2]{knutsonzinnjustin}
    \item $H$ and $K$ 3-step -- \cite[Theorem 3]{knutsonzinnjustin}
\end{enumerate}

\rbox{is this reference list complete?}
}



We will not go into much more detail; rather, we will include just enough discussion of 2-step and 3-step puzzles to be able to ultimately comment on how our results for Grassmannian puzzles extend to them.
Meanwhile, the situation for 4-step puzzles is more complicated 
and there is less explicitly written in the literature about their properties. For these reasons, we will abandon the topic of 4-step puzzles until Section \ref{section:furtherquestions}, where we discuss ``further questions.''





\subsection{Properties of puzzles}
\label{section:propertiesofpuzzles}

We now give a selection of puzzle properties that will be used in this paper.

\subsubsection{Basic symmetries of puzzles}
\label{section:basicsymmetriesofpuzzles}

\begin{enumerate}
    \item \textbf{Rotational symmetry.} Triangular puzzles using pieces in $\Hfrak \cup \left\{\raisebox{-.4\height}{
\begin{tikzpicture}[x=1pt,y=1pt,yscale=1,xscale=1,scale=.68]
\draw [fill={rgb, 255:red, 255; green, 255; blue, 255 }] (0.0,0.0) -- (15.0,25.98) -- (30.0,0.0) -- cycle;
\draw (7.5,12.99) node  [rotate=60]  {$\mathtt{1\! 0}$};
\draw (22.5,12.99) node  [rotate=-60]  {$\mathtt{1\! 0}$};
\draw (15.0,0.0) node  [rotate=0]  {$\mathtt{1\! 0}$};
\end{tikzpicture}
},\raisebox{-.4\height}{
\begin{tikzpicture}[x=1pt,y=1pt,yscale=1,xscale=1,scale=.68]
\draw [fill={rgb, 255:red, 255; green, 255; blue, 255 }] (30.0,25.98) -- (15.0,0.0) -- (0.0,25.98) -- cycle;
\draw (22.5,12.99) node  [rotate=60]  {$\mathtt{1\! 0}$};
\draw (7.5,12.99) node  [rotate=-60]  {$\mathtt{1\! 0}$};
\draw (15.0,25.98) node  [rotate=0]  {$\mathtt{1\! 0}$};
\end{tikzpicture}
} \right\}$ have $120^\circ$ rotational symmetry (also called 3-fold symmetry), where rotating a puzzle by a multiple of $120^\circ$ yields another puzzle with the same set of pieces. This operation gives bijections
    $$\{\inlinetri_{\lambda,\mu,\nu} \text{-puzzles}\} \leftrightarrow \{\inlinetri_{\mu,\nu,\lambda} \text{-puzzles}\} \leftrightarrow \{\inlinetri_{\nu,\lambda,\nu} \text{-puzzles}\}.$$ 

    For generalized polygonal puzzles using pieces in $\Hfrak \cup \left\{\raisebox{-.4\height}{
\begin{tikzpicture}[x=1pt,y=1pt,yscale=1,xscale=1,scale=.68]
\draw [fill={rgb, 255:red, 255; green, 255; blue, 255 }] (0.0,0.0) -- (15.0,25.98) -- (30.0,0.0) -- cycle;
\draw (7.5,12.99) node  [rotate=60]  {$\mathtt{1\! 0}$};
\draw (22.5,12.99) node  [rotate=-60]  {$\mathtt{1\! 0}$};
\draw (15.0,0.0) node  [rotate=0]  {$\mathtt{1\! 0}$};
\end{tikzpicture}
},\raisebox{-.4\height}{
\begin{tikzpicture}[x=1pt,y=1pt,yscale=1,xscale=1,scale=.68]
\draw [fill={rgb, 255:red, 255; green, 255; blue, 255 }] (30.0,25.98) -- (15.0,0.0) -- (0.0,25.98) -- cycle;
\draw (22.5,12.99) node  [rotate=60]  {$\mathtt{1\! 0}$};
\draw (7.5,12.99) node  [rotate=-60]  {$\mathtt{1\! 0}$};
\draw (15.0,25.98) node  [rotate=0]  {$\mathtt{1\! 0}$};
\end{tikzpicture}
} \right\}$, rotating in $60^\circ$ increments yields valid puzzles, though rotating by an odd multiple sends the piece \raisebox{-.4\height}{
\begin{tikzpicture}[x=1pt,y=1pt,yscale=1,xscale=1,scale=.68]
\draw [fill={rgb, 255:red, 255; green, 255; blue, 255 }] (0.0,0.0) -- (15.0,25.98) -- (30.0,0.0) -- cycle;
\draw (7.5,12.99) node  [rotate=60]  {$\mathtt{1\! 0}$};
\draw (22.5,12.99) node  [rotate=-60]  {$\mathtt{1\! 0}$};
\draw (15.0,0.0) node  [rotate=0]  {$\mathtt{1\! 0}$};
\end{tikzpicture}
} to \raisebox{-.4\height}{
\begin{tikzpicture}[x=1pt,y=1pt,yscale=1,xscale=1,scale=.68]
\draw [fill={rgb, 255:red, 255; green, 255; blue, 255 }] (30.0,25.98) -- (15.0,0.0) -- (0.0,25.98) -- cycle;
\draw (22.5,12.99) node  [rotate=60]  {$\mathtt{1\! 0}$};
\draw (7.5,12.99) node  [rotate=-60]  {$\mathtt{1\! 0}$};
\draw (15.0,25.98) node  [rotate=0]  {$\mathtt{1\! 0}$};
\end{tikzpicture}
}, and vice versa. When using the equivariant piece \raisebox{-.4\height}{
\begin{tikzpicture}[x=1pt,y=1pt,yscale=1,xscale=1,scale=.68]
\draw [fill={rgb, 255:red, 255; green, 246; blue, 74 }] (0.0,0.0) -- (15.0,25.98) -- (30.0,0.0) -- (15.0,-25.98) -- cycle;
\draw (7.5,12.99) node  [rotate=60]  {$\mathtt{0}$};
\draw (22.5,12.99) node  [rotate=-60]  {$\mathtt{1}$};
\draw (22.5,-12.99) node  [rotate=60]  {$\mathtt{0}$};
\draw (7.5,-12.99) node  [rotate=-60]  {$\mathtt{1}$};
\end{tikzpicture}
}, only $180^\circ$ rotation gives another valid puzzle.

   \item  \textbf{Duality symmetry.} Reflecting a puzzle across the $y$-axis while exchanging all $\zerostr$s and $\onestr$s is a reversible operation that results in another puzzle using the same set of pieces, for any puzzle piece set. Generalizing to the case of hexagons, and letting $^*$ denote the operation of reversing a binary string and exchanging $\zerostr$s and $\onestr$s, this yields a bijection 
    $$\{\inlinehex_{\alpha,\beta,\gamma,\delta,\epsilon,\zeta} \text{-puzzles}\} \leftrightarrow \{\inlinehex_{\epsilon^*,\delta^*,\gamma^*,\beta^*,\alpha^*,\zeta^*}\text{-puzzles}\}.$$

    \item \textbf{Commutativity.} For triangular puzzles using any set of pieces, the sum of the weights of all $\inlinetri_{\lambda,\mu,\nu}$-puzzles equals the sum of the weights of all $\inlinetri_{\mu,\lambda,\nu}$-puzzles. Furthermore, if the equivariant piece is not used, there is a bijection 
    $$\{\inlinetri_{\lambda,\mu,\nu} \text{-puzzles}\} \leftrightarrow \{\inlinetri_{\mu,\lambda,\nu} \text{-puzzles}\}.$$
    However, unlike in the previous two symmetries, the bijection for puzzle commutativity is more difficult (see \cite[\S 3.3]{Purbhoo} for a proof for puzzles using pieces in $\Hfrak$), and there is no visually apparent correspondence between the puzzles.

    The commutative properties of generalized polygonal puzzles will make up most of the main results of this paper, so will be discussed later.
\end{enumerate}

\justthesis{
\subsubsection{Basic symmetries of puzzles}
\label{section:basicsymmetriesofpuzzles}

\explainalot{Here we describe three fundamental symmetries of puzzles, as well as giving the corresponding statements in terms of structure constants and intersection theory.}

\begin{enumerate}
    \item \textbf{Rotational symmetry.} 

    \begin{itemize}
        \item \textbf{Classical puzzles.} Rotating a puzzle (without equivariant pieces) by $120^\circ$ is a reversible operation that results in another puzzle using the same set of pieces. 
        
        \textit{Example:
        }

    \begin{center}
    \raisebox{-.5\height}{\input{examples/01011_01101_01110}} {\large $\xleftrightarrow[]{\text{rotate}}$}
    \raisebox{-.5\height}{\input{examples/01011_01101_01110_rotate120}} {\large $\xleftrightarrow[]{\text{rotate}}$}
    \raisebox{-.5\height}{\input{examples/01011_01101_01110_rotate240}}
    \end{center}

    This operation gives bijections
    $$\{\inlinetri_{\lambda,\mu,\nu} \text{-puzzles}\} \leftrightarrow \{\inlinetri_{\mu,\nu,\lambda} \text{-puzzles}\} \leftrightarrow \{\inlinetri_{\nu,\lambda,\nu} \text{-puzzles}\}.$$ 
    We call this 3-fold rotational symmetry. This corresponds to the equalities of Littlewood-Richardson numbers
    $$c_{\lambda,\mu}^{\nu^\vee}= c_{\mu,\nu}^{\lambda^\vee}= c_{\nu,\lambda}^{\mu^\vee}$$
    and, geometrically,
    $$\underset{\Gr(k;\C^n)}{\int} [X_{\lambda}][X_{\mu}][X^{\nu^\vee}]= \underset{\Gr(k;\C^n)}{\int} [X_{\mu}][X_{\nu}][X^{\lambda^\vee}] = \underset{\Gr(k;\C^n)}{\int} [X_{\nu}][X_\lambda][X^{\mu^\vee}].$$
    This also works for all the extended puzzle piece sets in Section \ref{section:puzzlepiecesets} sans the equivariant piece, and we have the analogous statements in their associated cohomology theories.

    \item \textbf{Generalized polygonal puzzles.} Now, rotating a puzzle using pieces in $\Hfrak$ in $60^\circ$ increments results in more puzzles by our broadened definition. Generalizing to the case of hexagons, $60^\circ$ rotation gives bijections
    \begin{multline*}
        \{\inlinehex_{\lambda,\mu,\nu,\xi,o,\pi} \text{-puzzles}\} \leftrightarrow \{\inlinehex_{\mu,\nu,\xi,o,\pi,\lambda} \text{-puzzles}\} \leftrightarrow \{\inlinehex_{\nu,\xi,o,\pi,\lambda,\mu} \text{-puzzles}\}  \\ \leftrightarrow \{\inlinehex_{\xi,o,\pi,\lambda,\mu,\nu} \text{-puzzles}\} \leftrightarrow \{\inlinehex_{o,\pi,\lambda,\mu,\nu,\xi} \text{-puzzles}\} \leftrightarrow \{\inlinehex_{\pi,\lambda,\mu,\nu,\xi,0} \text{-puzzles}\}.
    \end{multline*}
    When using the equivariant piece, only $180^\circ$ rotation gives another valid puzzle. If using the set $\Hfrak \cup \left\{\raisebox{-.4\height}{
\begin{tikzpicture}[x=1pt,y=1pt,yscale=1,xscale=1,scale=.68]
\draw [fill={rgb, 255:red, 255; green, 255; blue, 255 }] (0.0,0.0) -- (15.0,25.98) -- (30.0,0.0) -- cycle;
\draw (7.5,12.99) node  [rotate=60]  {$\mathtt{1\! 0}$};
\draw (22.5,12.99) node  [rotate=-60]  {$\mathtt{1\! 0}$};
\draw (15.0,0.0) node  [rotate=0]  {$\mathtt{1\! 0}$};
\end{tikzpicture}
} \right\}$, rotating by odd multiples of $60^\circ$ produces puzzles using the set $\Hfrak \cup \left\{\raisebox{-.4\height}{
\begin{tikzpicture}[x=1pt,y=1pt,yscale=1,xscale=1,scale=.68]
\draw [fill={rgb, 255:red, 255; green, 255; blue, 255 }] (30.0,25.98) -- (15.0,0.0) -- (0.0,25.98) -- cycle;
\draw (22.5,12.99) node  [rotate=60]  {$\mathtt{1\! 0}$};
\draw (7.5,12.99) node  [rotate=-60]  {$\mathtt{1\! 0}$};
\draw (15.0,25.98) node  [rotate=0]  {$\mathtt{1\! 0}$};
\end{tikzpicture}
} \right\}$, and vice versa.
    \end{itemize}

    \item \textbf{Duality symmetry.}
    \begin{itemize}
        \item \textbf{Classical puzzles.} Let $^*$ be the operation of reversing a string and exchanging $\zerostr$s and $\onestr$s. Reflecting a puzzle across the $y$-axis and exchanging $\zerostr$s and $\onestr$s is a reversible operation that results in another puzzle using the same set of pieces.

        \textit{Example:}

    \begin{center}
    \raisebox{-.5\height}{\input{examples/01011_01101_01110}} {\large $\xleftrightarrow[\text{exchange } \zerostr\leftrightarrow \onestr]{\text{reflect across } y\text{-axis,}}$}
    \raisebox{-.5\height}{\input{examples/01011_01101_01110_dual}} 
    \end{center}

    This operation gives a bijection
    $$\{\inlinetri_{\lambda,\mu,\nu} \text{-puzzles}\} \leftrightarrow \{\inlinetri_{\mu^*,\lambda^*,\nu^*} \text{-puzzles}\}.$$
    This corresponds to the equality of Littlewood-Richardson numbers
    $$c_{\lambda,\mu}^{\nu^\vee}= c_{\mu^*,\lambda^*}^{(\nu^\vee)^*}$$ 
    and, geometrically,
    $$\underset{\Gr(k;\C^n)}{\int} [X_{\lambda}][X_{\mu}][X^{\nu^\vee}]=\underset{\Gr(n-k,(\C^n)^*)}{\int}[X_{\lambda^*}][X_{\mu^*}][X^{(\nu^\vee)^*}],$$
    which comes from the isomorphism $\Gr(k;\C^n) \rightarrow \Gr(n-k,(\C^n)^*)$ defined by $V\mapsto V^\perp$. 
    This symmetry works the same for all extended puzzle piece sets as well.

    \item \textbf{Generalized polygonal puzzles.} The operation of reflecting a puzzle across the $y$-axis and exchanging $\zerostr$s and $\onestr$s works the same as with classical puzzles. Generalizing to the case of hexagons, this operation gives a bijection 
    $$\{\inlinehex_{\alpha,\beta,\gamma,\delta,\epsilon,\zeta} \text{-puzzles}\} \leftrightarrow \{\inlinehex_{\epsilon^*,\delta^*,\gamma^*,\beta^*,\alpha^*,\zeta^*}\text{-puzzles}\}.$$
    \end{itemize}

    \item \textbf{Commutativity.}

    \begin{itemize}
        \item \textbf{Classical puzzles.} If we swap the NW and NE labels on a labeled boundary and look for all puzzles with that boundary labeling, we get the same structure constant as with the original boundary labeling. However, unlike in the previous two symmetries, the bijection for puzzle commutativity is more difficult (see \cite[\S 3.3]{Purbhoo} for a proof for puzzles using pieces in $\Hfrak$), and there is no visually apparent relationship between the puzzles.

    \textit{Example:}
    \begin{center}
    \raisebox{-.5\height}{\input{examples/01011_01101_01110}} {\large $\xleftrightarrow[]{\text{commute}}$}
    \raisebox{-.5\height}{\input{examples/01011_01101_01110_commute}} 
    \end{center}
    So there is a bijection
    $$\{\inlinetri_{\lambda,\mu,\nu} \text{-puzzles}\} \leftrightarrow \{\inlinetri_{\mu,\lambda,\nu} \text{-puzzles}\}.$$


    This corresponds to the equality of Littlewood-Richardson numbers
    $$c_{\lambda,\mu}^{\nu^\vee}= c_{\mu,\lambda}^{\nu^\vee}$$ 
    and, geometrically, 
    $$\underset{\Gr(k;\C^n)}{\int} [X_{\lambda}][X_{\mu}][X^{\nu^\vee}] = \underset{\Gr(k;\C^n)}{\int} [X_{\mu}][X_{\lambda}][X^{\nu^\vee}].$$
    Puzzles using the extended puzzle piece sets also have this property of commutativity.

        \item \textbf{Generalized polygonal puzzles.} Generalizations of commutativity for puzzles with different polygonal boundary shapes make up most of the main results of the paper, so will be discussed later.
    \end{itemize}

\end{enumerate}

} 

\subsubsection{Other relevant properties of puzzles}
\label{section:otherpuzzleproperties}

\begin{enumerate}
    \item \textbf{Unique identity puzzle.} 
    \justthesis{
    If we allow puzzle pieces in $\Hfrak \cup \left\{\raisebox{-.4\height}{
\begin{tikzpicture}[x=1pt,y=1pt,yscale=1,xscale=1,scale=.68]
\draw [fill={rgb, 255:red, 255; green, 246; blue, 74 }] (0.0,0.0) -- (15.0,25.98) -- (30.0,0.0) -- (15.0,-25.98) -- cycle;
\draw (7.5,12.99) node  [rotate=60]  {$\mathtt{0}$};
\draw (22.5,12.99) node  [rotate=-60]  {$\mathtt{1}$};
\draw (22.5,-12.99) node  [rotate=60]  {$\mathtt{0}$};
\draw (7.5,-12.99) node  [rotate=-60]  {$\mathtt{1}$};
\end{tikzpicture}
},\raisebox{-.4\height}{
\begin{tikzpicture}[x=1pt,y=1pt,yscale=1,xscale=1,scale=.68]
\draw [fill={rgb, 255:red, 255; green, 255; blue, 255 }] (0.0,0.0) -- (15.0,25.98) -- (30.0,0.0) -- cycle;
\draw (7.5,12.99) node  [rotate=60]  {$\mathtt{1\! 0}$};
\draw (22.5,12.99) node  [rotate=-60]  {$\mathtt{1\! 0}$};
\draw (15.0,0.0) node  [rotate=0]  {$\mathtt{1\! 0}$};
\end{tikzpicture}
} \right\}$ and we fix labels of $\lambda$ and $\sort(\lambda)$ on the NW and NE sides of a triangular boundary, then there exists a unique filling of the boundary with puzzle pieces. Furthermore, this filling uses only pieces in $\Hfrak$ and produces the label $\lambda^\vee$ on the South side. (See Figure \ref{fig:trivialtrianglefillingone}.) If we fix labels of $\lambda$ and $\sort(\lambda)$ on the South and NW sides respectively, then there exists a unique filling of the boundary with puzzle pieces. This filling uses only pieces in $\Hfrak$ and produces the label $\lambda^\vee$ on the NE side. (See Figure \ref{fig:trivialtrianglefillingtwo}.)
    }
This property is described in Figure \ref{fig:trivialtrianglefilling} and can be justified by the fact that the Schubert class corresponding to the string $\sort(\lambda)$ is the identity element in the cohomology theories that use any subset of those pieces.

On the other hand, this property does not hold when the piece \raisebox{-.4\height}{
\begin{tikzpicture}[x=1pt,y=1pt,yscale=1,xscale=1,scale=.68]
\draw [fill={rgb, 255:red, 255; green, 255; blue, 255 }] (30.0,25.98) -- (15.0,0.0) -- (0.0,25.98) -- cycle;
\draw (22.5,12.99) node  [rotate=60]  {$\mathtt{1\! 0}$};
\draw (7.5,12.99) node  [rotate=-60]  {$\mathtt{1\! 0}$};
\draw (15.0,25.98) node  [rotate=0]  {$\mathtt{1\! 0}$};
\end{tikzpicture}
} is included.

\begin{figure}[h!]
    \centering
    \begin{subfigure}[t]{0.4\textwidth}
    \centering
     \input{examples/trivialtrianglefillingone_no_opacity}
    \caption{}
    \label{fig:trivialtrianglefillingone}   
    \end{subfigure}
    ~
    \begin{subfigure}[t]{0.4\textwidth}
    \centering
    \input{examples/trivialtrianglefillingtwo_no_opacity}
    \caption{}
    \label{fig:trivialtrianglefillingtwo}
    \end{subfigure}
    \caption{In each case, if we fix the binary strings $\lambda$ and $\sort(\lambda)$ labeling the path highlighted in red, there exists a unique filling of the boundary using pieces in $\Hfrak \cup \left\{\protect\raisebox{-.4\height}{\input{puzzle_pieces/eqvtpiece}},\protect\raisebox{-.4\height}{\input{puzzle_pieces/cccdeltapiece}} \right\}$, which ends up only using pieces in $\Hfrak$ and produces $\lambda^\vee$ on the remaining side.
    }
    \label{fig:trivialtrianglefilling}
\end{figure}

\justthesis{
These facts will later be used for the operations of ``completing to a triangle,'' which will let us give a geometric interpretation to puzzles with polygonal boundary.
}


    \item \textbf{Discrete Green's theorem.} We assign a unit vector in $\R^2$ to each edge of a puzzle piece depending on its orientation and whether it is labeled $\zerostr$, $\onestr$, or $\tenstr$, as in \cite[Lemma 1]{knutsonzinnjustin}. Then the sum of the vectors around each puzzle piece equals $\vec{0}$, and by a discrete Green's theorem, this implies that the sum of the vectors around any boundary filled with puzzle pieces must also be $\vec{0}$.
    The following statements are not too difficult to derive from this, and they will be used in this paper:
    \begin{enumerate}
        \item 
        \justthesis{Suppose the NW and NE sides of a triangular boundary are labeled with binary strings $\lambda$ and $\mu$ respectively. If $\lambda$ has $m$ more $\onestr$s than $\mu$, then in any puzzle filling the boundary, the South side label must have exactly $m$ $\tenstr$ labels. On the other hand, if $\lambda$ has more $\zerostr$s than $\mu$, there can exist no puzzle filling the boundary. If $\lambda$ and $\mu$ have the same content, then for any puzzle filling the boundary, the South side label must have the same content as $\lambda$ and $\mu$.}
        
        Suppose two sides of a triangular boundary are labeled with binary strings $\lambda$ and $\mu$ (i.e. they have no $\tenstr$s).
        If $\lambda$ and $\mu$ do not have the same content, then there exists no filling of the boundary with puzzle pieces such that the South side is also labeled with a binary string (hence no puzzles). If $\lambda$ and $\mu$ have the same content, then for any filling of the boundary with puzzle pieces, the South side label must have the same content as $\lambda$ and $\mu$.

        \item For a trapezoidal boundary, if the two (necessarily) equal-length sides are labeled with binary strings $\lambda$ and $\mu$ with matching content, and the side opposite the base is labeled with a binary string $\rho$, then for any filling of the boundary with puzzle pieces, the label appearing on the base of the trapezoid must be a binary string with the content of $\lambda$ plus the content of $\rho$. 
        \item If all the sides of a trapezoidal boundary are labeled with binary strings, then the existence of a filling of the boundary with puzzle pieces implies that the strings on the pair of (necessarily) equal-length sides must have matching content.

    \end{enumerate}

    \explainalot{We would like to particularly stress that we always say ``binary string'' to mean a string of $\zerostr$s and $\onestr$s (with no $\tenstr$s).}

\end{enumerate}

\subsubsection{Properties of 2-step and 3-step puzzles}
\label{section:propertiesofpuzzles2step3step}
All of the properties of puzzles we have described also apply to 2-step and 3-step puzzles. They possess exactly the same basic symmetries given in Section \ref{section:basicsymmetriesofpuzzles}.
For Property 1 of Section \ref{section:otherpuzzleproperties}, we have the identical statements, with the operation ``$\sort$'' now more generally sorting all the entries of a string weakly from least to greatest. We also have discrete Green's theorem statements for 2-step and 3-step puzzles (see Lemma 2 and Lemma 3 of \cite{knutsonzinnjustin}), so that for Property 2 of Section \ref{section:otherpuzzleproperties}, we have the identical statements if we replace ``binary string'' with ``string of $\zerostr$s, $\onestr$s, and $\twostr$s'' for 2-step or ``string of $\zerostr$s, $\onestr$s, $\twostr$s, and $\threestr$s'' for 3-step. Another way to express this is by saying that the string has no composite entries. 


\subsection{Summary of results}

\begin{figure}[h]
\centering
\begin{subfigure}[b]{0.21\textwidth}
\centering
\input{boundaries_commute/triangleboundarycommutesplit}
\caption{$\protect\inlinetri_{\alpha\gamma,\nu,\delta\beta}$}
\label{fig:triangleboundarycommutesplit}
\end{subfigure}
\begin{subfigure}[b]{0.21\textwidth}
\centering
\input{boundaries_commute/trapezoidboundarycommute}
\caption{$\protect\inlinetrapside_{\beta,\gamma,\nu,\delta}$}
\label{fig:trapezoidboundarycommute}
\end{subfigure}
\begin{subfigure}[b]{0.17\textwidth}
\centering
\input{boundaries_commute/parallelogramboundarycommute}
\caption{$\protect\inlinepar_{\alpha,\gamma,\beta,\delta}$}
\label{fig:parallelogramboundarycommute}
\end{subfigure}
\begin{subfigure}[b]{0.18\textwidth}
\centering
\input{boundaries_commute/rhombusboundarycommute}
\caption{$\protect\inlinerhombus_{\alpha,\gamma,\beta,\delta}$}
\label{fig:rhombusboundarycommute}
\end{subfigure}
\begin{subfigure}[b]{0.20\textwidth}
\centering
\input{boundaries_commute/pentagonboundarycommute}
\caption{$\protect\inlinepentside_{\alpha,\beta,\gamma,\delta,\epsilon,\zeta}$}
\label{fig:pentagonboundarycommute}
\end{subfigure}

\vspace{\baselineskip}

\begin{subfigure}[b]{0.24\textwidth}
\centering
\input{boundaries_commute/hexagonboundarycommuteopposite}
\caption{$\protect\inlinehex_{\alpha,\beta,\gamma,\delta,\epsilon,\zeta}$}
\label{fig:hexagonboundarycommutetwoway}
\end{subfigure}
\begin{subfigure}[b]{0.24\textwidth}
\centering
\input{boundaries_commute/hexagonboundarycommutetwowayside}
\caption{$\protect\inlinehex_{\alpha,\beta,\gamma,\delta,\epsilon,\zeta}$}
\label{fig:hexagonboundarycommutethreeway}
\end{subfigure}
\begin{subfigure}[b]{0.24\textwidth}
\centering
\input{boundaries_commute/hexagonboundarycommutethreeway}
\caption{$\protect\inlinehex_{\alpha,\beta,\gamma,\delta,\epsilon,\zeta}$}
\label{fig:hexagonboundarycommuteopposite}
\end{subfigure}
\begin{subfigure}[b]{0.24\textwidth}
\centering
\input{boundaries_commute/hexagonboundarycommuteallway}
\caption{$\protect\inlinehex_{\alpha,\beta,\gamma,\delta,\epsilon,\zeta}$}
\label{fig:hexagonboundarycommuteopposite}
\end{subfigure}
\caption{Boundaries labeled with binary strings, where matching colored squiggly lines indicate which labels commute when they have matching content.}
\label{fig:boundaryshapescommute}
\end{figure}

Most of our main results are about puzzles using pieces in $\Hfrak$ and performing computations in $H^*$. These are summarized below.
\begin{itemize}
    \item We can commute the boundary labels of puzzles with the boundary types shown in Figure \ref{fig:boundaryshapescommute} when those labels have matching content, and the number of puzzles is preserved.
    \item For all of the commutative properties, we present purely geometric proofs. The geometric proofs will rely on operations that allow us to complete puzzles of these shapes to triangular puzzles, which have a familiar geometric Schubert calculus interpretation.
    \item For most of these boundary types we give a formula for the associated structure constant in terms of simpler Littlewood-Richardson numbers.
    \item We give a geometric interpretation and proof of $180^\circ$ rotational symmetry of hexagonal puzzles.
    \item We additionally give combinatorial, puzzle-based proofs for most statements.
\end{itemize}

Besides these, there is also a result about parallelogram-shaped equivariant puzzles, which describes how the associated structure constant in $H_T^*(\Gr(k;\C^n))$ is affected by the operation of commuting opposite boundary side labels. We also get that the number of puzzles is preserved (which is not the case when commuting the labels on the NW and NE sides of triangular equivariant puzzles).

We give a number of remarks about how our results translate to K-theory and puzzles using pieces in $\Hfrak \cup \left\{\raisebox{-.4\height}{
\begin{tikzpicture}[x=1pt,y=1pt,yscale=1,xscale=1,scale=.68]
\draw [fill={rgb, 255:red, 255; green, 255; blue, 255 }] (0.0,0.0) -- (15.0,25.98) -- (30.0,0.0) -- cycle;
\draw (7.5,12.99) node  [rotate=60]  {$\mathtt{1\! 0}$};
\draw (22.5,12.99) node  [rotate=-60]  {$\mathtt{1\! 0}$};
\draw (15.0,0.0) node  [rotate=0]  {$\mathtt{1\! 0}$};
\end{tikzpicture}
} \right\}$. In many cases the statements and arguments are nearly identical to those for ordinary cohomology.

Finally, in Section \ref{section:2step3step} we describe the extent to which our results are known to extend to 2-step and 3-step puzzles. The majority of our theorems go through with almost no change.

\section{Preliminaries}

\subsection{Converting between Young diagrams and binary strings}
\label{converting}

\begin{figure}[h]
    \centering
    \Large{$\lambda =\textcolor{red}{\zerostr\zerostr\zerostr}\textcolor{blue}{\onestr\onestr}\textcolor{red}{\zerostr}\textcolor{blue}{\onestr}\textcolor{red}{\zerostr}\textcolor{blue}{\onestr}$} \qquad\LARGE{$\leftrightarrow$} \qquad \raisebox{-.5\height}{\input{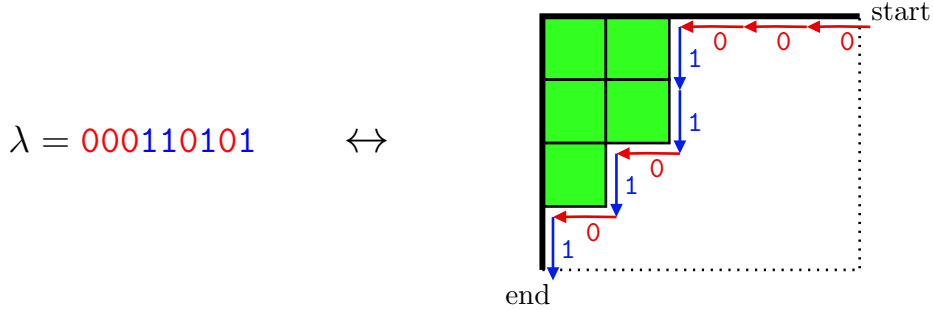}}
    \caption{An illustration of the bijection that identifies a binary string $\lambda\in\binom{[n]}{k}$ with a Young diagram in a $k\times(n-k)$ ambient rectangle.}
    \label{fig:stringyoungdiagrambijection}
\end{figure}

As mentioned in the introduction, the structure constants $c_{\lambda,\mu}^\nu$ in the Schubert basis in $H^*(\Gr(k;\C^n))$ are the Littlewood-Richardson numbers. Usually Littlewood-Richardson numbers are indexed by partitions, or equivalently, Young diagrams. But we can identify binary strings $\lambda\in\binom{[n]}{k}$ with Young diagrams fitting in a $k\times(n-k)$ ``ambient rectangle.'' This is done via a bijection that sends $\lambda$ to the Young diagram obtained by starting at the NW corner of the ambient rectangle and reading $\lambda$, stepping left one unit when encountering a $\zerostr$ and down one unit when encountering a $\onestr$. This is illustrated in Figure \ref{fig:stringyoungdiagrambijection}.

When Littlewood-Richardson numbers $c_{\lambda,\mu}^\nu$ are indexed by Young diagrams, there is no need to fix a choice of ``ambient rectangle'' to bound the three diagrams. However, when we use binary strings to index our structure constants for the purposes of Schubert calculus, it is important that the content of $\zerostr$s and $\onestr$ be consistent among $\lambda$, $\mu$, and $\nu$, and the content of the strings is what determines the size of the ambient rectangle that bounds their associated Young diagrams. Since the Young diagram for $\nu$ must always contain those for $\lambda$ and $\mu$ (lest $c_{\lambda,\mu}^\nu=0$ automatically), we can ``pad out'' the strings $\lambda$ and $\mu$ with additional $\zerostr$s in front and $\onestr$s in back in order to make their content match that of $\nu$,
and this preserves the associated Young diagrams. 

To be more explicit, suppose we have binary strings $\lambda\in\binom{[l_0+l_1]}{l_1}$, $\mu\in\binom{[m_0+m_1]}{m_1}$, and $\nu\in\binom{[n_0+n_1]}{n_1}$, where $n_0\geq l_0$, $n_1\geq l_1$, $n_0\geq m_0$, and $n_1\geq m_1$. Then we pad out $\lambda$ and $\mu$ with extra $\zerostr$s and $\onestr$s to get $\zerostr^{n_0-l_0}\lambda\onestr^{n_1-l_1}$ and $\zerostr^{n_0-m_0}\mu\onestr^{n_1-m_1}$ so that now we view them as living in the same ambient $n_1\times n_0$ rectangle as $\nu$. 

We now establish a notational convention that will be used going forward. When writing a structure constant, we will use $\pad$ to denote the operation of padding out the strings in the subscript so that their content matches the content of the string in the superscript. So we define
$$c_{\lambda\pad,\mu\pad}^\nu:= c_{\zerostr^{n_0-l_0}\lambda\onestr^{n_1-l_1},\zerostr^{n_0-m_0}\mu\onestr^{n_1-m_1}}^\nu.$$
Note that the notation $c_{\lambda\pad,\mu\pad}^\nu$ is independent of how $\lambda$, $\mu$, and $\nu$ were defined, as long as $\nu$ is assumed to have at least as much content as $\lambda$ and $\mu$.


\subsection{An equality of Littlewood-Richardson numbers: $c_{\sort(\lambda)\sort(\mu),\nu}^{\lambda\mu} = c_{\lambda\pad,\mu\pad}^\nu$}

In this section we will use Littlewood-Richardson rules that involve counting skew tableaux to establish an equality of Littlewood-Richardson numbers,
which will be useful in later proofs. We state the equality in the proposition below.

\begin{proposition}
    \label{prop:equalityofLR}
    For binary strings $\lambda\in\binom{[l_0+l_1]}{l_1}$, $\mu\in\binom{[m_0+m_1]}{m_1}$, and $\nu\in\binom{[l_0+l_1+m_0+m_1]}{l_1+m_1}$, we have an equality of Littlewood-Richardson numbers 
    $$c_{\sort(\lambda)\sort(\mu),\nu}^{\lambda\mu} = c_{\lambda\pad,\mu\pad}^\nu.$$
\end{proposition}

We will give the proof shortly, after discussing the relevant Littlewood-Richardson rules.

\subsubsection{Two Littlewood-Richardson rules}
Letting $\lambda$, $\mu$, and $\nu$ be Young diagrams, we have 
the following definitions from \cite{Fulton}:

\begin{definition}
For any semistandard Young tableau $U_\circ$ with shape $\mu$, define
$$\mathcal{S}(\nu/\lambda,U_\circ):=\{\text{skew tableaux } S \text{ on } \nu/\lambda : S \text{ rectifies to } U_\circ\}.$$
\end{definition}

\begin{definition}
For any semistandard Young tableau $V_\circ$ with shape $\nu$, define
$$\mathcal{T}(\lambda*\mu,V_\circ):=\{\text{skew tableaux } S \text{ on } \lambda*\mu : S \text{ rectifies to } V_\circ\},$$
where $\lambda*\mu$ is the skew diagram obtained by joining the northeast corner of $\lambda$ to the southwest corner of $\mu$.
\end{definition}

From \cite{Fulton} also, we then have two Littlewood-Richardson rules:
\[c_{\lambda,\mu}^\nu= |\mathcal{S}(\nu/\lambda,U_\circ)| \tag{Rule 1}\]
\[c_{\lambda,\mu}^\nu= |\mathcal{T}(\lambda*\mu,V_\circ)| \tag{Rule 2}\]

Above we are indexing with Young diagrams, so as written there is no consideration for ambient rectangles, but if we are going to commit strictly to talking about binary strings, as we will now proceed to do, we will write $c_{\lambda\pad,\mu\pad}^\nu$ instead of $c_{\lambda,\mu}^\nu$.

\subsubsection{Proof of Proposition \ref{prop:equalityofLR}}

\begin{figure}[]
    \centering
    \begin{subfigure}[t]{0.3\textwidth}
        \centering
        \input{young_diagrams/youngdiagram_mu}
        \caption{Young diagram for $\lambda=\onestr\zerostr\onestr\zerostr\zerostr$}
    \end{subfigure}
    \;
    \begin{subfigure}[t]{0.3\textwidth}
        \centering
        \input{young_diagrams/youngdiagram_lambda}
        \caption{Young diagram for $\mu=\zerostr\onestr\zerostr\zerostr\onestr\onestr\zerostr$}
    \end{subfigure}

     \vspace{2\baselineskip}
    \begin{subfigure}[t]{0.3\textwidth}
        \centering
        \input{young_diagrams/youngdiagram_muprime}
        \caption{Young diagram for $\lambda\pad=\zerostr\zerostr\zerostr\zerostr\onestr\zerostr\onestr\zerostr\zerostr\onestr\onestr\onestr$}
    \end{subfigure}
    ~
    \;
    \begin{subfigure}[t]{0.3\textwidth}
        \centering
        \input{young_diagrams/youngdiagram_lambdaprime}
        \caption{Young diagram for $\mu\pad=\zerostr\zerostr\zerostr\zerostr\onestr\zerostr\zerostr\onestr\onestr\zerostr\onestr\onestr$}
    \end{subfigure}
    ~
    \;
    \begin{subfigure}[t]{0.3\textwidth}
        \centering
        \input{young_diagrams/youngdiagram_nu}
        \caption{Young diagram for $\tilde{\mu}:=\nu=\zerostr\zerostr\zerostr\onestr\zerostr\onestr\zerostr\onestr\zerostr\onestr\zerostr\onestr$}
    \end{subfigure}

     \vspace{2\baselineskip}
    \begin{subfigure}[t]{0.31\textwidth}
        \centering
        \input{young_diagrams/youngdiagram_rect}
        \caption{Young diagram for 
        $\tilde{\lambda}:=\sort(\lambda)\sort(\mu)=\zerostr\zerostr\zerostr\zerostr\onestr\onestr\onestr\zerostr\zerostr\zerostr\onestr\onestr$}
    \end{subfigure}
    ~ 
    \begin{subfigure}[t]{0.31\textwidth}
        \centering
        \input{young_diagrams/youngdiagram_mulambda}
        \caption{Young diagram for $\tilde{\nu}:=\lambda\mu=\onestr\zerostr\onestr\zerostr\zerostr\zerostr\onestr\zerostr\zerostr\onestr\onestr\zerostr$}
    \end{subfigure}
    ~ 
    \begin{subfigure}[t]{0.33\textwidth}
        \centering
        \input{young_diagrams/youngdiagram_lambdastarmu}
        \caption{Skew diagram for $\tilde{\nu}/\tilde{\lambda}= \mu*\lambda=\onestr\zerostr\onestr\zerostr\zerostr\zerostr\onestr\zerostr\zerostr\onestr\onestr\zerostr/\zerostr\zerostr\zerostr\onestr\onestr\zerostr\zerostr\zerostr\zerostr\onestr\onestr\onestr$}
    \end{subfigure}
    \caption{An example of Young diagrams, shown in their ambient rectangles, corresponding to the binary strings defined in the proof in the case where $\lambda=\onestr\zerostr\onestr\zerostr\zerostr$, $\mu=\zerostr\onestr\zerostr\zerostr\onestr\onestr\zerostr$, and $\nu=\zerostr\zerostr\zerostr\onestr\zerostr\onestr\zerostr\onestr\zerostr\onestr\zerostr\onestr$.}
    \label{fig:youngdiagramsexample}
    
\end{figure}

\begin{proof}
We will show that applying Rule 2 to find $c_{\lambda\pad,\mu\pad}^\nu$ is equivalent to applying Rule 1 to find $c_{\tilde{\lambda},\tilde{\mu}}^{\tilde{\nu}}$, where $\tilde{\lambda}:= \sort(\lambda)\sort(\mu)= \zerostr^{l_0}\onestr^{l_1}\zerostr^{m_0}\onestr^{m_1}$ (its associated Young diagram is an $m_1\times l_0$ rectangle), $\tilde{\mu}:=\nu$, and $\tilde{\nu}:=\lambda\mu$. See Figure \ref{fig:youngdiagramsexample} for an example of the strings and associated Young diagrams at play, with ambient rectangles shown.

Let $V_\circ$ be a tableau with shape $\tilde{\mu}=\nu$. Then, noting that $\tilde{\nu}/\tilde{\lambda}= \mu*\lambda$, we have 
$$\mathcal{S}(\tilde{\nu}/\tilde{\lambda},V_\circ) = \mathcal{T}(\mu*\lambda,V_\circ).$$
This implies the equality of Littlewood-Richardson numbers
$$c_{\tilde{\lambda},\tilde{\mu}}^{\tilde{\nu}}= c_{\mu\pad,\lambda\pad}^\nu.$$ Since $c_{\mu\pad,\lambda\pad}^\nu=c_{\lambda\pad,\mu\pad}^\nu$ and by definition $c_{\tilde{\lambda},\tilde{\mu}}^{\tilde{\nu}} = c_{\sort(\lambda)\sort(\mu),\nu}^{\lambda\mu}$, we obtain the equality
$$c_{\sort(\lambda)\sort(\mu),\nu}^{\lambda\mu} = c_{\lambda\pad,\mu\pad}^\nu.$$

\end{proof}

 \begin{corollary}
\label{cor:pushforwardcoeffspuzzles}
    For $\lambda\in\binom{[l_0+l_1]}{l_1}$, $\mu\in\binom{[m_0+m_1]}{m_1}$, and $\nu\in\binom{[l_0+l_1+m_0+m_1]}{l_1+m_1}$, we have
    $$\#\{\inlinetri_{\sort(\lambda)\sort(\mu),\nu,(\lambda\mu)^\vee} \text{-puzzles}\} = \#\{\inlinetri_{\lambda\pad,\mu\pad,\nu^\vee} \text{-puzzles}\}.$$
    (See Figure \ref{fig:pushforwardpuzzles}).
\end{corollary}
\begin{proof}
    We simply translate Proposition \ref{prop:equalityofLR} into puzzle terms.
\end{proof}

\begin{figure*}[h]
    \centering
    \raisebox{-.5\height}{\input{proof_boundary_diagrams/pushforwardboundarycolor}} 
    \Huge $\quad \leftrightarrow \quad$ \normalsize \raisebox{-.5\height}{\input{proof_boundary_diagrams/lambdapadmupadnuboundarycolor}}
    \caption{Puzzles with the boundary $\protect\inlinetri_{\sort(\lambda)\sort(\mu),\nu,(\lambda\mu)^\vee}$ are in bijection with puzzles with the boundary $\protect\inlinetri_{\lambda\pad,\mu\pad,\nu^\vee}$. }
    \label{fig:pushforwardpuzzles}
\end{figure*}

\begin{remark}
In Figure \ref{fig:pushforwardpuzzles}, the thinner black line cutting across the left triangular boundary is drawn to indicate that there can be no $\tenstr$ puzzle piece edge labels appearing along that line in any puzzle with the given boundary labels, since $\lambda$ and $\sort(\lambda)$ have matching content. (This is by the discrete Green's theorem property given in Section \ref{section:otherpuzzleproperties}.) This also means we can divide the puzzle into a smaller triangular puzzle and a trapezoidal puzzle, glued together along that line, which is suggested by the division in color. We will continue with these conventions throughout the paper, drawing a thinner line cutting across a boundary shape to mean that no $\tenstr$s can appear along it and separating subpuzzles by color. 
\end{remark}

\begin{remark}
The equality of Littlewood-Richardson numbers in Proposition \ref{prop:equalityofLR} and the corresponding bijection of puzzles in Corollary \ref{cor:pushforwardcoeffspuzzles} are not understood by the author by any geometric argument of the flavor of those given in this paper, nor by any puzzle-based combinatorial argument. Finding such an understanding would be of interest.
\end{remark}





\subsection{Geometric lemmas}

The following two lemmas will form the foundation that all of our geometric arguments will build upon later. The bulk of the direct intersection-theoretic work underlying our results will be neatly confined to this section, and will just need to be referenced in later proofs.

For the entirety of this section, let $a_0,a_1,c_0,c_1$ be non-negative integers, and define $a:=a_0+a_1$ and $c:=c_0+c_1$. Let $\flag{F}$ and $\flag{\Ftilde}$ be the standard and anti-standard complete flags in $\C^{a+c}$, respectively.


\begin{lemma}
\label{lemma:directsum}
Let $T:=(\C^\times)^{a+c}$, and define the $T$-equivariant closed immersion
$$\Omega: \Gr(a_1;F_a)\times \Gr(c_1;\Ftilde_c) \hookrightarrow  \Gr(a_1+c_1;\C^{a+c}), \quad (V,W)\mapsto V\oplus W.$$ 
Then for $\alpha,\beta\in\binom{[a]}{a_1}$ and $\gamma,\delta\in\binom{[c]}{c_1}$, we have
\begin{equation*}
    X_{\alpha\gamma}\cap X^{(\delta\beta)^\vee} = \Omega((X_\alpha\cap X^{\beta^\vee}) \times (X_\gamma\cap X^{\delta^\vee})), \tag{1}
\end{equation*}
and 
\begin{equation*}
[X_{\alpha\gamma}][X^{(\delta\beta)^\vee}]= \Omega_*([X_\alpha][X^{\beta^\vee}]\otimes [X_\gamma][X^{\delta^\vee}]) \tag{2}
\end{equation*}
in $H_T^*(\Gr(a_1+c_1;\C^{a+c}))$ and in $H^*(\Gr(a_1+c_1;\C^{a+c}))$.



\end{lemma}

\begin{proof}

    We will show that $X_{\alpha\gamma}\cap X^{(\delta\beta)^\vee} = \Omega((X_\alpha\cap X^{\beta^\vee}) \times (X_\gamma\cap X^{\delta^\vee}))$ by proving an inclusion in both directions.
    
    First assume that $V\in X_{\alpha\gamma}\cap X^{(\delta\beta)^\vee}= X_{\alpha\gamma}(\flag{F})\cap X_{\delta\beta}(\flag{\Ftilde})$ (recall that by definition, $X_{\alpha\gamma}=X_{\alpha\gamma}(\flag{F})$ and $X^{(\delta\beta)^\vee}=X_{\delta\beta}(\flag{\Ftilde})$).  Consider the conditions encoded by the strings $\alpha\gamma$ and $\delta\beta$. Since $\alpha=\alpha_1\alpha_2\cdots\alpha_a\in\binom{[a]}{a_1}$, we have $\alpha_1+\alpha_2+\cdots+\alpha_a=a_1$, which encodes the condition that $\dim(V\cap F_a)\geq a_1$. And since $\delta=\delta_1\delta_2\cdots\delta_c \in\binom{[c]}{c_1}$), we have $\delta_1+\delta_2+\cdots+\delta_c=c_1$, which encodes the condition that $\dim(V\cap \Ftilde_c)\geq c_1$. We have that $F_a\cap \Ftilde_c=\{\mathbf{0}\}$, so $(V\cap F_a)\cap (V\cap \Ftilde_c)=\{\mathbf{0}\}$, and hence $\dim((V\cap F_a)\oplus (V\cap \Ftilde_c))= \dim(V\cap F_a) + \dim(V\cap \Ftilde_c) \geq a_1+c_1$. But since $(V\cap F_a)\oplus (V\cap \Ftilde_c)\subseteq V$ and $\dim(V)=a_1+c_1$, the above inequalities must in fact all be equalities, and we have that $(V\cap F_a) \oplus (V\cap \Ftilde_c) = V$.

    Let us simplify notation by writing $V_a:=V\cap F_a$ and $V_c:=V\cap \Ftilde_c$ going forward. We will now view $V_a$ and $V_c$ as elements of the smaller Grassmannians $\Gr(a_1;F_a)$ and $\Gr(c_1;\Ftilde_c)$ respectively, and we will speak of Schubert varieties in these Grassmannians.

    Let $\flag{F'}$ denote the flag $\{\mathbf{0}\}\subset F_1\subset F_2 \subset \cdots\subset F_a$, and let $\flag{\Ftilde'}$ denote the flag $\{\mathbf{0}\}\subset \langle\bfe_{a}\rangle \subset \langle\bfe_{a},\bfe_{a-1}\rangle \subset \cdots\subset \langle\bfe_{a},\bfe_{a-1},\ldots,\bfe_{1}\rangle = F_a$. Let $\flag{F''}$ denote the flag $\{\mathbf{0}\}\subset \langle\bfe_{a+1}\rangle \subset \langle\bfe_{a+1},\bfe_{a+2}\rangle \subset \cdots \subset \langle \mathbf{e}_{a+1},\ldots,\mathbf{e}_{a+c}\rangle= \Ftilde_c$, and let $\flag{\Ftilde''}$ denote the flag $\{\mathbf{0}\}\subset \Ftilde_1\subset \Ftilde_2 \subset \cdots\subset \Ftilde_c$. 


    Now, since $V\in X_{\alpha\gamma}(\flag{F})$, we know $\dim(V\cap F_i)\geq \alpha_i$ for all $1\leq i\leq a$. But then since $V\cap F_i=V\cap(F_a\cap F_i)= (V\cap F_a)\cap F_i= V_a\cap F_i$, we can say that $V_a$ also satisfies the intersection conditions encoded by $\alpha$ with respect to the flag $\flag{F'}$, so we have $V_a\in X_\alpha(\flag{F'})\subseteq \Gr(a_1;F_a)$. 

    For $1\leq j\leq c$, we have $\dim(V\cap F_{a+j})\geq a_1+\gamma_1+\cdots+\gamma_j$ due to the conditions encoded by $\alpha\gamma$. 
    We also have $V\cap F_{a+j}=(V_a\oplus V_c)\cap (F_{a}\oplus F''_j) = V_a\oplus (V_c\cap (F_{a}\oplus F''_j)) = V_a\oplus (V_c\cap F''_j)$, given that $V_a\subseteq (F_a\oplus F''_j)$ and $V_c\cap (F_a\oplus F''_j)=(V_c\cap \Ftilde_c)\cap (F_a\oplus F''_j) = V_c\cap (\Ftilde_c\cap (F_a\oplus F''_j)) = V_c\cap F''_j$, and hence $\dim(V\cap F_{a+j}) = \dim(V_a\oplus(V_c\cap F''_j))= \dim(V_a)+\dim (V_c\cap F''_j)=a_1+\dim (V_c\cap F''_j)$.
    Putting these together, we get $a_1+\dim (V_c\cap F''_j)\geq a_1+\gamma_1+\gamma_2+\cdots+\gamma_j$, and thus $\dim (V_c\cap F''_j)\geq \gamma_1+\gamma_2+\cdots+\gamma_j$.
    This tells us that $V_c$ satisfies the conditions encoded by $\gamma$ with respect to the flag $F''$. So we have $V_c\in X_\gamma(\flag{F''})\subseteq \Gr(c_1;\Ftilde_c)$. 

    A symmetric argument proves that we have $V_c\in X_{\delta}(\flag{\Ftilde''})\subseteq \Gr(c_1;\Ftilde_c)$ and $V_a\in X_\beta(\flag{\Ftilde'})\subseteq \Gr(a_1;F_a)$. To see this, simply let $\delta$ and $\Ftilde''_c$ take the place of $\alpha$ and $F_a$, and let $\beta$ and $\Ftilde'_j$ take the place of $\gamma$ and $F''_j$ in the previous two paragraphs. 
    
    Putting it all together, now we have $V_a\in X_\alpha(\flag{F'})\cap X_\beta(\flag{\Ftilde'})= X_\alpha\cap X^{\beta^\vee}\subseteq \Gr(a_1;F_a)$ and $V_c\in X_\gamma(\flag{F''})\cap X_{\delta}(\flag{\Ftilde''})= X_\gamma\cap X^{\delta^\vee}\subseteq \Gr(c_1;\Ftilde_c)$. So $(V_a,V_c)\in (X_\alpha\cap X^{\beta^\vee}) \times (X_\gamma\cap X^{\delta^\vee})$, and since $\Omega(V_a,V_c)=V$, this implies that $X_{\alpha\gamma}\cap X^{(\delta\beta)^\vee} \subseteq \Omega((X_\alpha\cap X^{\beta^\vee}) \times (X_\gamma\cap X^{\delta^\vee}))$.





    Now we will prove the opposite inclusion. We begin with the assumption that $(U,W)\in (X_\alpha\cap X^{\beta^\vee}) \times (X_\gamma\cap X^{\delta^\vee})= (X_\alpha\times X_\gamma)\cap (X^{\beta^\vee}\times X^{\delta^\vee})= (X_\alpha(\flag{F'})\times X_\gamma(\flag{F''}))\cap (X_{\beta}(\flag{\Ftilde'})\times X_{\delta}(\flag{\Ftilde''}))$. Let $V=\Omega(U,W)=U\oplus W$. We have $\dim(V\cap F_i)\geq \dim(U\cap F_i)= \dim(U\cap F'_i) \geq \alpha_1+\alpha_2+\cdots+\alpha_i$ for $1\leq i \leq a$. We also have $\dim(V\cap F_{a+j})\geq \dim(U\cap F_a\oplus W\cap F''_{j}) = \dim(U)+\dim(W\cap F''_j)\geq a_1+\gamma_1+\cdots+\gamma_j= \alpha_1+\alpha_2+\cdots+\alpha_a+ \gamma_1+\cdots+\gamma_j$ for $1\leq j\leq c$. So $V$ satisfies the conditions encoded by $\alpha\gamma$ with respect to the flag $\flag{F}$, hence $V\in X_{\alpha\gamma}$. 
    
    We have $\dim(V\cap \Ftilde_i)\geq \dim(W\cap \Ftilde_i)= \dim(W\cap \Ftilde''_i) \geq \delta_1+\delta_2+\cdots+\delta_i$ for $1\leq i \leq c$. We also have $\dim(V\cap \Ftilde_{c+j})\geq \dim(W\cap \Ftilde_c \oplus U\cap \Ftilde'_{j}) = \dim(W)+ \dim(U\cap \Ftilde'_j) \geq c_1+\beta_1+\beta_2+\cdots+\beta_j= \delta_1+\delta_2+\cdots+\delta_c + \beta_1+\beta_2+\cdots+\beta_j$ for $1\leq j\leq a$. So $V$ satisfies the conditions encoded by $\delta\beta$ with respect to the flag $\flag{\Ftilde}$, hence $V\in X^{(\delta\beta)^\vee}$. 

    So finally we have $\Omega(U,W)=V\in X_{\alpha\gamma}\cap X^{(\delta\beta)^\vee}$, and thus we conclude $\Omega((X_\alpha\cap X^{\beta^\vee}) \times (X_\gamma\cap X^{\delta^\vee}))\subseteq  X_{\alpha\gamma}\cap X^{(\delta\beta)^\vee}$.

    With both directions of inclusion, we have proved that $$X_{\alpha\gamma}\cap X^{(\delta\beta)^\vee} = \Omega((X_\alpha\cap X^{\beta^\vee}) \times (X_\gamma\cap X^{\delta^\vee})).$$

    Recall from Section \ref{section:intersectiontheory} that $[X_\lambda][X^\mu]=[X_\lambda\cap X^\mu]$ for any choice of $\lambda,\mu$. Then in $H^*_T(\Gr(a_1+c_1;\C^{a+c})$ or in $H^*(\Gr(a_1+c_1;\C^{a+c})$, we have 
    \begin{multline*}
    [X_{\alpha\gamma}][X^{(\delta\beta)^\vee}] =[X_{\alpha\gamma}\cap X^{(\delta\beta)^\vee}]= [\Omega((X_\alpha\cap X^{\beta^\vee}) \times (X_\gamma\cap X^{\delta^\vee}))] =\Omega_*([(X_\alpha\cap X^{\beta^\vee}) \times (X_\gamma\cap X^{\delta^\vee})]) \\
    =  \Omega_*([X_\alpha\cap X^{\beta^\vee}] \otimes [X_\gamma\cap X^{\delta^\vee}]) =  \Omega_*([X_\alpha][X^{\beta^\vee}] \otimes [X_\gamma][X^{\delta^\vee}]),
    \end{multline*} 
    as claimed.
\end{proof}

\begin{lemma}
\label{lemma:commutesplit}

Define block matrices
$$\Phi_a:=\left[
\begin{array}{c|c}
J_{a} & \mathbf{0} \\ \hline
\mathbf{0} & I_{c}
\end{array}\right]
\quad \text{and} \quad
\Phi_c:=\left[
\begin{array}{c|c}
I_{a} & \mathbf{0} \\ \hline
\mathbf{0} & J_{c}
\end{array}\right]
$$
in $\GL(\C^{a+c})$, where $I_a$ and $J_a$ (resp. $I_c$ and $J_c$) denote the $a\times a$ (resp. $c\times c$) identity and anti-diagonal permutation matrices, respectively.

Let $\alpha,\beta\in\binom{[a]}{a_1}$ and $\gamma,\delta\in\binom{[c]}{c_1}$. Then,
\begin{enumerate}[(i)]
    \item in $H_T^*(\Gr(a_1+c_1;\C^{a+c}))$, we have 
$$[X_{\alpha\gamma}][X^{(\delta\beta)^\vee}]= \Phi_a\cdot ([X_{\beta\gamma}][X^{(\delta\alpha)^\vee}]) = \Phi_c\cdot ([X_{\alpha\delta}][X^{(\gamma\beta)^\vee}]) = \Phi_c\cdot\Phi_a\cdot ([X_{\beta\delta}][X^{(\gamma\alpha)^\vee}]).$$
    \item in $H^*(\Gr(a_1+c_1;\C^{a+c}))$, we have 
$$[X_{\alpha\gamma}][X^{(\delta\beta)^\vee}]=  [X_{\beta\gamma}][X^{(\delta\alpha)^\vee}] =  [X_{\alpha\delta}][X^{(\gamma\beta)^\vee}] = [X_{\beta\delta}][X^{(\gamma\alpha)^\vee}].$$

\end{enumerate}

\end{lemma}

\begin{proof}
   $\GL(\C^{a+c})$ acts on $\Gr(a_1+c_1;\C^{a+c})$, and $\GL(F_a)\times \GL(\Ftilde_c)$ acts on $\Gr(a_1;F_a)\times \Gr(c_1;\Ftilde_c)$.

Let $\Omega$ be as defined in Lemma \ref{lemma:directsum}. Then note for $(g_1,g_2)\in \GL(F_a)\times \GL(\Ftilde_c)$, we have $\Omega((g_1,g_2)\cdot X)=(g_1\oplus g_2)\cdot \Omega(X)$ for a subvariety $X\subseteq \GL(F_a)\times \GL(\Ftilde_c)$.
And $\Omega_*[((g_1,g_2)\cdot X)]= [\Omega((g_1,g_2)\cdot X)] = [(g_1\oplus g_2)\cdot \Omega(X)]= (g_1\oplus g_2)\cdot [\Omega(X)]= (g_1\oplus g_2)\cdot \Omega_*([X])$.

From Lemma \ref{lemma:directsum}, we have $[X_{\alpha\gamma}][X^{\delta\delta)^\vee}]= \Omega_*([X_\alpha][X^{\beta^\vee}]\otimes [X_\gamma][X^{\delta^\vee}])$.


We also have $(J_a,I_c)\cdot ([X_\beta][X^{\alpha^\vee}]\otimes [X_\gamma][X^{\delta^\vee}])= (J_a,I_c)\cdot ([X_\beta\cap X^{\alpha^\vee}]\otimes [X_\gamma\cap X^{\delta^\vee}]) = (J_a,I_c)\cdot ([X_\beta\cap X^{\alpha^\vee}\times X_\gamma\cap X^{\delta^\vee}]) = [(J_a,I_c)\cdot (X_\beta\cap X^{\alpha^\vee}\times X_\gamma\cap X^{\delta^\vee})] =[J_a\cdot (X_\beta\cap X^{\alpha^\vee})\times I_c\cdot (X_\gamma\cap X^{\delta^\vee})] = [J_a\cdot (X_\beta\cap X^{\alpha^\vee})]\otimes [I_c\cdot (X_\gamma\cap X^{\delta^\vee})] = [X_\alpha\cap X^{\beta^\vee}]\otimes [X_\gamma\cap X^{\delta^\vee}] = [X_\alpha][X^{\beta^\vee}]\otimes [X_\gamma][X^{\delta^\vee}]$.


Then, putting all these together, we get $[X_{\alpha\gamma}][X^{(\delta\beta)^\vee}] = \Omega_*([X_\alpha][X^{\beta^\vee}]\otimes [X_\gamma][X^{\delta^\vee}])= \Omega_*((J_a,I_c)\cdot ([X_\beta][X^{\alpha^\vee}]\otimes [X_\gamma][X^{\delta^\vee}])) = \Phi_a\cdot \Omega_*([X_\beta][X^{\alpha^\vee}]\otimes [X_\gamma][X^{\delta^\vee}]) = \Phi_a\cdot ([X_{\beta\gamma}][X^{(\delta\alpha)^\vee}])$. 


With a symmetric argument, we can prove $[X_{\alpha\gamma}][X^{(\delta\beta)^\vee}]= \Phi_c\cdot ([X_{\alpha\delta}][X^{(\gamma\beta)^\vee}])$. 
In this case we have that $(I_a,J_c)\cdot ([X_\alpha][X^{\beta^\vee}]\otimes [X_\delta][X^{\gamma^\vee}])= [X_\alpha][X^{\beta^\vee}]\otimes [X_\gamma][X^{\delta^\vee}]$. 
Then $[X_{\alpha\gamma}][X^{(\delta\beta)^\vee}] = \Omega_*([X_\alpha][X^{\beta^\vee}]\otimes [X_\gamma][X^{\delta^\vee}])= \Omega_*((I_a,J_c)\cdot ([X_\alpha][X^{\beta^\vee}]\otimes [X_\delta][X^{\gamma^\vee}])) = \Phi_c\cdot \Omega_*([X_\alpha][X^{\beta^\vee}]\otimes [X_\delta][X^{\gamma^\vee}]) = \Phi_c\cdot ([X_{\alpha\delta}][X^{(\gamma\beta)^\vee}])$.

Finally we have $(J_a,J_c)\cdot ([X_\beta][X^{\alpha^\vee}]\otimes [X_\delta][X^{\gamma^\vee}])= (I_a,J_c) \cdot(J_a,I_c)\cdot ([X_\beta][X^{\alpha^\vee}]\otimes [X_\delta][X^{\gamma^\vee}])= [X_\alpha][X^{\beta^\vee}]\otimes [X_\gamma][X^{\delta^\vee}]$, and $[X_{\alpha\gamma}][X^{(\delta\beta)^\vee}] = \Omega_*([X_\alpha][X^{\beta^\vee}]\otimes [X_\gamma][X^{\delta^\vee}])= \Omega_*((J_a,J_c)\cdot ([X_\beta][X^{\alpha^\vee}]\otimes [X_\delta][X^{\gamma^\vee}])) = \Phi_c\cdot\Phi_a\cdot \Omega_*([X_\beta][X^{\alpha^\vee}]\otimes [X_\delta][X^{\gamma^\vee}]) = \Phi_c\cdot\Phi_a\cdot ([X_{\beta\delta}][X^{(\gamma\alpha)^\vee}])$.

In ordinary cohomology $H^*(\Gr(a_1+c_1;\C^{a+c})$, $\Phi_a$ and $\Phi_c$ act trivially on the classes. We have $\Phi_a\cdot ([X_{\beta\gamma}][X^{(\delta\alpha)^\vee}]) = \Phi_a\cdot ([X_{\beta\gamma}\cap X^{(\delta\alpha)^\vee}])= [\Phi_a\cdot (X_{\beta\gamma}\cap X^{(\delta\alpha)^\vee})] =[(\Phi_a\cdot X_{\beta\gamma})\cap (\Phi_a\cdot X^{(\delta\alpha)^\vee})]= [\Phi_a\cdot X_{\beta\gamma}][\Phi_a\cdot X^{(\delta\alpha)^\vee}]$, but $[\Phi_a\cdot X_{\beta\gamma}]=[ X_{\beta\gamma}]$, and $[\Phi_a\cdot X^{(\delta\alpha)^\vee}]=[ X^{(\delta\alpha)^\vee}]$ as discussed in Section \ref{section:schubertvarieties}. The same argument applies with $\Phi_c$. This proves the version of the statement for $H^*(\Gr(a_1+c_1;\C^{a+c})$.
\end{proof}

\subsubsection{$d$-step statements}
\label{Section:geometriclemmasdstep}


Here we will generalize Lemmas \ref{lemma:directsum} and \ref{lemma:commutesplit} to the case of $d$-step flag manifolds. For this section, let $a_1\leq a_2\leq\cdots\leq a_d\leq a$ and $c_1\leq c_2\leq\cdots\leq c_d\leq c$ be non-negative integers. Let $\alpha,\beta$ be strings with the content of $\zerostr^{a-a_d}\onestr^{a_d-a_{d-1}}\twostr^{a_{d-1}-a_{d-2}}\cdots\dstr^{a_1}$, and let $\gamma,\delta$ be strings with the content of $\zerostr^{c-c_d}\onestr^{c_d-c_{d-1}}\twostr^{c_{d-1}-c_{d-2}}\cdots\dstr^{c_1}$.  Let $\flag{F}$ and $\flag{\Ftilde}$ be the standard and anti-standard complete flags in $\C^{a+c}$, respectively.

\begin{customlemma}{1'}
\label{lemma:directsumdstep}
Let $T:=(\C^\times)^{a+c}$, and define the $T$-equivariant closed immersion
    \begin{gather*}
    \Omega: \Fl(a_1,a_2,\ldots,a_d;F_a)\times \Fl(c_1,c_2,\ldots,c_d;\Ftilde_c) \hookrightarrow  \Fl(a_1+c_1,a_2+c_2,\ldots,a_d+c_d; \C^{a+c}), \\
    ((V_1\subset V_2\subset \cdots\subset V_d), (W_1\subset W_2\subset \cdots\subset W_d)) \mapsto ((V_1\oplus W_1)\subset (V_2\oplus W_2)\subset \cdots\subset (V_d\oplus W_d)).
    \end{gather*}

    Then we have
    \begin{equation*}
    X_{\alpha\gamma}\cap X^{(\delta\beta)^\vee} = \Omega((X_\alpha\cap X^{\beta^\vee}) \times (X_\gamma\cap X^{\delta^\vee})), \tag{1}
\end{equation*}
and 
\begin{equation*}
[X_{\alpha\gamma}][X^{(\delta\beta)^\vee}]= \Omega_*([X_\alpha][X^{\beta^\vee}]\otimes [X_\gamma][X^{\delta^\vee}]) \tag{2}
\end{equation*}
in $H_T^*(\Fl(a_1+c_1,a_2+c_2,\ldots,a_d+c_d; \C^{a+c}))$ and in $H^*(\Fl(a_1+c_1,a_2+c_2,\ldots,a_d+c_d; \C^{a+c}))$.
\end{customlemma}

\begin{proof}
    For each $1\leq i\leq d$, we have a commutative diagram
 
    \begin{tikzcd}
{\Fl(a_1,a_2,\ldots,a_d;F_a)\times \Fl(c_1,c_2,\ldots,c_d;\Ftilde_c)} \arrow[r, "\Omega", hook] \arrow[d, "p_a^i\times p_c^i", two heads] & {\Fl(a_1+c_1,a_2+c_2,\ldots,a_d+c_d; \C^{a+c})} \arrow[d, "p_{a+c}^i", two heads] \\
\Gr(a_i;F_a)\times \Gr(c_i;\Ftilde_c) \arrow[r, "\Omega^i", hook]                                                                         & \Gr(a_i+c_i;\C^{a+c}),                                                            
\end{tikzcd}

where, recall, $p^i$ is used to denote a map that sends each $d$-step flag to its $i$th subspace, and $\Omega^i$ is the direct sum map.
We can directly apply statement (1) of Lemma \ref{lemma:directsum} (where here $\Omega^i$ takes the place of $\Omega$), along with the commutative diagram, to write
\begin{align*}
p_{a+c}^i\left(X_{\alpha\gamma}\cap X^{(\delta^\beta)^\vee}\right) &=\Omega^i\left((p_a^i\times p_c^i)((X_{\alpha}\cap X^{\beta^\vee}) \times (X_{\gamma}\cap X^{\delta^\vee}))\right) \\
&= p_{a+c}^i\left(\Omega((X_{\alpha}\cap X^{\beta^\vee}) \times (X_{\gamma}\cap X^{\delta^\vee}))\right).
\end{align*}
Given the general fact that $\flag{V}\in X_\lambda$ if and only if $p^i(\flag{V})\in p^i(X_\lambda)$ for all $i$, it is not difficult from this point to stitch these $d$ components together to deduce that overall
$$X_{\alpha\gamma}\cap X^{(\delta^\beta)^\vee} = \Omega((X_{\alpha}\cap X^{\beta^\vee}) \times (X_{\gamma}\cap X^{\delta^\vee})),$$
and the cohomological statements easily follow.
\end{proof}

\begin{customlemma}{2'}
\label{lemma:commutesplitdstep}
Define $\Phi_a$ and $\Phi_c$ as in Lemma \ref{lemma:commutesplit}. Then,
\begin{enumerate}[(i)]
    \item in $H_T^*(\Fl(a_1+c_1,a_2+c_2,\ldots,a_d+c_d; \C^{a+c}))$, we have 
$$[X_{\alpha\gamma}][X^{(\delta\beta)^\vee}]= \Phi_a\cdot ([X_{\beta\gamma}][X^{(\delta\alpha)^\vee}]) = \Phi_c\cdot ([X_{\alpha\delta}][X^{(\gamma\beta)^\vee}]) = \Phi_c\cdot\Phi_a\cdot ([X_{\beta\delta}][X^{(\gamma\alpha)^\vee}]).$$
    \item in $H^*(\Fl(a_1+c_1,a_2+c_2,\ldots,a_d+c_d; \C^{a+c}))$, we have 
$$[X_{\alpha\gamma}][X^{(\delta\beta)^\vee}]=  [X_{\beta\gamma}][X^{(\delta\alpha)^\vee}] =  [X_{\alpha\delta}][X^{(\gamma\beta)^\vee}] = [X_{\beta\delta}][X^{(\gamma\alpha)^\vee}].$$
\end{enumerate}
\end{customlemma}
\begin{proof}
    The proof is identical to that of Lemma \ref{lemma:commutesplit}.
\end{proof}

\subsubsection{K-theory statements}
\begin{remark}
    Lemmas \ref{lemma:directsum} and \ref{lemma:directsumdstep} have analogous statements in K-theory, namely that in $K(\Gr(a_1+c_1;\C^{a+c}))$ or $K(\Fl(a_1+c_1,a_2+c_2,\ldots,a_d+c_d; \C^{a+c}))$, we have
    $$[\Ocal_{\alpha\gamma}][J_{a+c}\cdot \Ocal_{\delta\beta}]= \Omega_*([\Ocal_\alpha][J_a\cdot \Ocal_{\beta}]\otimes [\Ocal_\gamma][J_c\cdot \Ocal_{\delta}]),$$ 
    where $J_{a+c}$, $J_a$, and $J_c$ are the $(a+c)\times (a+c)$, $a\times a$, and $c\times c$ anti-diagonal permutation matrices, respectively. Equivalently we have
    $$[\Ocal_{\alpha\gamma}][\Ocal_{\delta\beta}]= \Omega_*([\Ocal_\alpha][\Ocal_{\beta}]\otimes [\Ocal_\gamma][\Ocal_{\delta}]).$$ 
    The proof is essentially identical to those of the lemmas, as the arguments are geometric.
\end{remark}

\begin{remark}
\label{remark:lemmatwoktheory}
    Lemmas \ref{lemma:commutesplit} and \ref{lemma:commutesplitdstep} have analogous statements in K-theory, namely that in $K(\Gr(a_1+c_1;\C^{a+c}))$ or $K(\Fl(a_1+c_1,a_2+c_2,\ldots,a_d+c_d; \C^{a+c}))$, we have
     $$[\Ocal_{\alpha\gamma}][J_{a+c}\cdot \Ocal_{\delta\beta}]=   [\Ocal_{\beta\gamma}][J_{a+c}\cdot \Ocal_{\delta\alpha}]=  [\Ocal_{\alpha\delta}][J_{a+c}\cdot \Ocal_{\gamma\beta}] =  [\Ocal_{\beta\delta}][J_{a+c}\cdot \Ocal_{\gamma\alpha}],$$
     and equivalently,
     $$[\Ocal_{\alpha\gamma}][\Ocal_{\delta\beta}]=   [\Ocal_{\beta\gamma}][\Ocal_{\delta\alpha}]=  [\Ocal_{\alpha\delta}][\Ocal_{\gamma\beta}] =  [\Ocal_{\beta\delta}][\Ocal_{\gamma\alpha}].$$
    The proof is essentially identical again.
\end{remark}

\subsection{Pushforward and pullback coefficients for the Grassmannian direct sum map}

The content of this section will be used to derive the formulas we will give later for the structure constants associated to different boundary shapes in terms of other Littlewood-Richardson numbers.

\explainalot{The following proposition also is the same for $d$-step}

\begin{proposition}
    \label{prop:pushforwardcoeff}
    Let $\lambda\in\binom{[l_0+l_1]}{l_1}$, $\mu\in\binom{[m_0+m_1]}{m_1}$, and $\nu\in\binom{[l_0+l_1+m_0+m_1]}{l_1+m_1}$, and define the map
    $$\Omega: \Gr(a_1;F_a)\times \Gr(c_1;\Ftilde_c) \hookrightarrow  \Gr(a_1+c_1;\C^{a+c}), \quad (V,W)\mapsto V\oplus W.$$
    Then for the Schubert basis expansions
$$\Omega_*([X^\lambda]\otimes[X^\mu])= \sum_{\nu} d_{\lambda,\mu}^\nu [X^\nu]$$
in $H^*(\Gr(a_1+c_1;\C^{a+c}))$ and
$$\Omega^*([X_\nu])=\sum_{\lambda,\mu} d_{\lambda,\mu}^\nu ([X_\lambda]\otimes[X_\mu])$$
in $H^*(\Gr(a_1;F_a)\times \Gr(c_1;\Ftilde_c))$, the coefficients are 
$$d_{\lambda,\mu}^\nu = c_{\sort(\lambda)\sort(\mu),\nu}^{\lambda\mu}.$$

\end{proposition}

\begin{proof}
    By Lemma \ref{lemma:directsum}, we have $[X_{\sort(\lambda)\sort(\mu)}][X^{\lambda\mu}]= \Omega_*([X_{\sort(\lambda)}][X^\lambda]\otimes [X_{\sort(\mu)}][X^\mu])$. We have that $X_{\sort(\lambda)}=\Gr(a_1;F_a)$ since $\sort(\lambda)$ encodes no essential conditions, so $[X_{\sort(\lambda)}]=[\Gr(a_1;F_a)]= 1\in H^*(\Gr(a_1;F_a))$. Similarly, $[X_{\sort(\mu)}]=[\Gr(c_1;\Ftilde_c)]=1\in H^*(\Gr(c_1;\Ftilde_c))$. So we can continue on to rewrite the equality as $[X_{\sort(\lambda)\sort(\mu)}][X^{\lambda\mu}]= \Omega_*(1\cdot [X^\lambda]\otimes 1\cdot [X^\mu]) = \Omega_*([X^\lambda]\otimes [X^\mu])$. 

    Then we have
    \begin{multline*}
    d_{\lambda,\mu}^\nu = \underset{\Gr(a_1+c_1;\C^{a+c})}{\int} \left(\sum_{\nu} d_{\lambda,\mu}^\nu [X^\nu]\right)[X_\nu] = \underset{\Gr(a_1+c_1;\C^{a+c})}{\int} \Omega_*([X^\lambda]\otimes[X^\mu])[X_\nu]\\
    = \underset{\Gr(a_1+c_1;\C^{a+c})}{\int} [X_{\sort(\lambda)\sort(\mu)}][X^{\lambda\mu}][X_\nu] = c_{\sort(\lambda)\sort(\mu),\nu}^{\lambda\mu}.
    \end{multline*}

    Now we add an explanation for why the coefficients in the Schubert basis expansion of $\Omega^*([X_\nu])$ should be the same.
    We use the projection formula (i.e. $f_*(\varphi)\smile \omega=f_*(\varphi\smile f^*(\omega))$ for any map $f$) to write $\Omega_*([X^\lambda]\otimes[X^\mu])[X_\nu]= \Omega_*\left(([X^\lambda]\otimes[X^\mu])\Omega^*([X_\nu])\right)$, and then we have 
    $$d_{\lambda,\mu}^\nu =\underset{\Gr(a_1+c_1;\C^{a+c})}{\int} \Omega_*([X^\lambda]\otimes[X^\mu])[X_\nu] = \underset{\Gr(a_1+c_1;\C^{a+c})}{\int} \Omega_*\left(([X^\lambda]\otimes[X^\mu])\Omega^*([X_\nu])\right).$$
    Since homology is a covariant functor and integration is pushforward to a point, and the composition $\Gr(a_1;F_a)\times \Gr(c_1;\Ftilde_c) \xrightarrow{\Omega} \Gr(a_1+c_1;\C^{a+c}) \rightarrow \pt$ is equal to the map to a point $\Gr(a_1;F_a)\times \Gr(c_1;\Ftilde_c)\rightarrow \pt$, we have
    $$\underset{\Gr(a_1+c_1;\C^{a+c})}{\int} \Omega_*(\varphi) = \underset{\Gr(a_1;F_a)\times \Gr(c_1;\Ftilde_c)}{\int} \varphi$$
    for any class $\varphi$. It then follows that 
    $$d_{\lambda,\mu}^\nu= \underset{\Gr(a_1+c_1;\C^{a+c})}{\int} \Omega_*\left(([X^\lambda]\otimes[X^\mu])\Omega^*([X_\nu])\right) = \underset{\Gr(a_1;F_a)\times \Gr(c_1;\Ftilde_c)}{\int} ([X^\lambda]\otimes[X^\mu])\Omega^*([X_\nu]),$$
    which implies that the coefficient on $[X_\lambda]\otimes[X_\mu]$ in the Schubert basis expansion of $\Omega^*([X_\nu])$ is equal to $d_{\lambda,\mu}^\nu$. 
\end{proof}

\begin{corollary}
\label{cor:pushforwardcoeffs}
The coefficients of Proposition \ref{prop:pushforwardcoeff} are the Littlewood-Richardson numbers:
    $$d_{\lambda,\mu}^\nu = c_{\lambda\pad,\mu\pad}^\nu.$$
\end{corollary}
\begin{proof}
    Combine the statements of Proposition \ref{prop:equalityofLR} and Proposition \ref{prop:pushforwardcoeff}.
\end{proof}

\begin{remark}
    Corollary \ref{cor:pushforwardcoeffs} is not a novel statement (e.g. see \cite[\S 1.1]{KnutsonLederer}), but we have not encountered a simple proof like ours in the literature. 
\end{remark}

\section{Split symmetry}

\begin{figure}[h]
    \centering
    \input{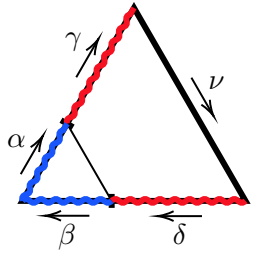}
    \caption{$\protect\inlinetri_{\alpha\gamma,\nu,\delta\beta}$, displaying ``split symmetry.''}
    \label{fig:triangleboundarycommutesplitsolo}
\end{figure}

For triangular puzzles, what we will refer to as ``split symmetry'' in this paper is a situation where the boundary labels are such that we can split the puzzle into two parts: a smaller triangular puzzle, and a trapezoidal puzzle. This is shown in Figure \ref{fig:triangleboundarycommutesplitsolo}, and we have also already seen an example of this in Figure \ref{fig:pushforwardpuzzles}.

This symmetry occurs when the puzzle has boundary $\inlinetri_{\alpha\gamma,\nu,\delta\beta}$, where $\alpha$ and $\beta$ have matching content, and $\gamma$ and $\delta$ have matching content. We will find that we can independently swap $\alpha$ with $\beta$ and $\gamma$ with $\delta$.

This is perhaps the most fundamental symmetry we will discuss, since it is actually the main building block for our proofs of the commutative properties of puzzles with all the various boundary shapes shown in Figure \ref{fig:boundaryshapescommute} in the introduction.

For the entirety of this section, let $a_0,a_1,c_0,c_1$ be non-negative integers, and define $a:=a_0+a_1$ and $c:=c_0+c_1$. Let us also shorten the notation for the spaces we will work with by defining $\Gr_{a+c}:=\Gr(a_1+c_1;\C^{a+c})$, $\Gr_a:=\Gr(a_1;F_a)$, and $\Gr_c:=\Gr(c_1;\Ftilde_c)$, where $\flag{F}$ and $\flag{\Ftilde}$ are the standard and anti-standard complete flags.

\subsection{Geometric interpretation of puzzles with split symmetry}

We will give triangular puzzles with split symmetry a geometric interpretation by looking at the corresponding structure constants and intersection theory. Using puzzle pieces in $\Hfrak$, we have

$$\#\{\inlinetri_{\alpha\gamma,\nu,\delta\beta} \text{-puzzles}\} = c_{\alpha\gamma,\nu}^{(\delta\beta)^\vee}.$$
With the next two propositions, we present a commutativity statement for these Littlewood-Richardson numbers and a formula for them in terms of simpler Littlewood-Richardson numbers, which will be used in our proofs about more complex boundary shapes later. 

\begin{proposition}
\label{prop:splitLRnumbers}
    For $\alpha,\beta\in\binom{[a]}{a_1}$, $\gamma,\delta\in\binom{[c]}{c_1}$, and $\nu\in\binom{[c+a]}{a_1+c_1}$, we have the following equalities of Littlewood-Richardson numbers:
    $$c_{\alpha\gamma,\nu}^{(\delta\beta)^\vee}= c_{\beta\gamma,\nu}^{(\delta\alpha)^\vee} = c_{\alpha\delta,\nu}^{(\gamma\beta)^\vee} = c_{\beta\delta,\nu}^{(\gamma\alpha)^\vee}.$$
\end{proposition}

\begin{proof}
    
    By using Lemma \ref{lemma:commutesplit} as it applies to $H^*$, we have $$\underset{\Gr_{a+c}}{\int} [X_{\alpha\gamma}][X^{(\delta\beta)^\vee}][X_\nu] = \underset{\Gr_{a+c}}{\int} [X_{\beta\gamma}][X^{(\delta\alpha)^\vee}][X_\nu] = \underset{\Gr_{a+c}}{\int} [X_{\alpha\delta}][X^{(\gamma\beta)^\vee}][X_\nu] = \underset{\Gr_{a+c}}{\int} [X_{\beta\delta}][X^{(\gamma\alpha)^\vee}][X_\nu].$$
    The above integrals, from left to right, equal the Littlewood-Richardson numbers $c_{\alpha\gamma,\nu}^{(\beta\delta)^\vee}$, $c_{\beta\gamma,\nu}^{(\alpha\delta)^\vee}$, $c_{\alpha\delta,\nu}^{(\beta\gamma)^\vee}$, $c_{\beta\delta,\nu}^{(\alpha\gamma)^\vee}$. Thus the claim immediately follows.
    
\end{proof}

\subsubsection{A formula in terms of Littlewood-Richardson numbers}

\begin{proposition}
\label{prop:splitsymmetryformula}
For $\alpha,\beta\in\binom{[a]}{a_1}$, $\gamma,\delta\in\binom{[c]}{c_1}$, and $\nu\in\binom{[c+a]}{a_1+c_1}$, and using puzzle pieces in $\Hfrak$, we have
    $$\#\{\inlinetri_{\alpha\gamma,\nu,\delta\beta} \text{-puzzles}\} =   \sum_{\lambda,\mu}c_{\lambda\pad,\mu\pad}^\nu \cdot c_{\alpha,\lambda}^{\beta^\vee}\cdot c_{\gamma,\mu}^{\delta^\vee}.$$
\end{proposition}

\begin{proof}
Using Lemma \ref{lemma:directsum} and Corollary \ref{cor:pushforwardcoeffs} along with the projection formula and basic steps, we have
\begin{align*}
\#\{\inlinetri_{\alpha\gamma,\nu,\delta\beta} \text{-puzzles}\} 
&= c_{\alpha\gamma,\nu}^{(\delta\beta)^\vee} \\
&= \underset{\Gr_{a+c}}{\int} [X_{\alpha\gamma}][X^{(\delta\beta)^\vee}][X_\nu]\\
&= \underset{\Gr_{a+c}}{\int} \Omega_*([X_\alpha][X^{\beta^\vee}]\otimes [X_\gamma][X^{\delta^\vee}])[X_\nu]\\
&=\underset{\Gr_{a}\times\Gr_{c}}{\int} ([X_\alpha][X^{\beta^\vee}]\otimes [X_\gamma][X^{\delta^\vee}])\Omega^*([X_\nu])\\
&= \underset{\Gr_{a}\times\Gr_{c}}{\int} ([X_\alpha][X^{\beta^\vee}]\otimes [X_\gamma][X^{\delta^\vee}])\left(\sum_{\lambda,\mu} d_{\lambda,\mu}^\nu ([X_\lambda]\otimes[X_\mu])\right)\\
&= \underset{\Gr_{a}\times\Gr_{c}}{\int} ([X_\alpha][X^{\beta^\vee}]\otimes [X_\gamma][X^{\delta^\vee}])\left(\sum_{\lambda,\mu} c_{\lambda\pad,\mu\pad}^\nu ([X_\lambda]\otimes[X_\mu])\right)\\
&= \sum_{\lambda,\mu}c_{\lambda\pad,\mu\pad}^\nu \underset{\Gr_{a}\times\Gr_{c}}{\int} [X_\alpha][X^{\beta^\vee}][X_\lambda] \otimes [X_\gamma][X^{\delta^\vee}][X_\mu] \\
&= \sum_{\lambda,\mu}c_{\lambda\pad,\mu\pad}^\nu \left(\underset{\Gr_{a}}{\int} [X_\alpha][X^{\beta^\vee}][X_\lambda]\right) \left(\underset{\Gr_{c}}{\int} [X_\gamma][X^{\delta^\vee}][X_\mu]\right) \\
&= \sum_{\lambda,\mu}c_{\lambda\pad,\mu\pad}^\nu \cdot c_{\alpha,\lambda}^{\beta^\vee}\cdot c_{\gamma,\mu}^{\delta^\vee}.
\end{align*}
\end{proof}

\subsection{Commutative property of puzzles with split symmetry}

Now we return to puzzles with split symmetry, and we restate our commutativity result in the language of puzzles.

\begin{theorem}[Commutative Property of Puzzles with Split Symmetry]
\label{thm:puzzlecommutesplit}
    For $\alpha,\beta\in\binom{[a]}{a_1}$, $\gamma,\delta\in\binom{[c]}{c_1}$, and $\nu\in\binom{[c+a]}{a_1+c_1}$, and puzzle pieces in $\Hfrak$, we have
    $$\#\{\inlinetri_{\alpha\gamma,\nu,\delta\beta} \text{-puzzles}\}= \#\{\inlinetri_{\beta\gamma,\nu,\delta\alpha} \text{-puzzles}\} = \#\{\inlinetri_{\alpha\delta,\nu,\gamma\beta} \text{-puzzles}\} = \#\{\inlinetri_{\beta\delta,\nu,\gamma\alpha} \text{-puzzles}\}$$
In other words, the pairs $\alpha,\beta$ and $\gamma,\delta$ commute independently. 
\end{theorem}

\begin{proof}
We give three proof methods. One is geometric, while the other two are combinatorial. The last argument uses a property we will prove in the next section (we still list it here because it is no less valid a method, despite the chronological inconvenience). 
\begin{enumerate}
    \item \textbf{Geometric proof.} The statement immediately follows from Proposition \ref{prop:splitLRnumbers}.

    \item \textbf{Proof by the commutative property of triangular puzzles.} It can also be proved by using the basic commutative property of triangular puzzles. 
    First we note that in any $\inlinetri_{\alpha\gamma,\nu,\delta\beta}$-puzzle, there can be no $\tenstr$s appearing along the red line shown in Figure \ref{fig:splitsymmetrypuzzle}, and if we call the string appearing along that line $\lambda$, then $\lambda$ must have the same content as $\alpha$ and $\beta$. (This comes from statement (a) under the discrete Green's theorem property given in Section \ref{section:otherpuzzleproperties}.) So we can divide the puzzle into a $\inlinetri_{\alpha,\lambda,\beta}$-puzzle and a trapezoidal $\inlinetrapside_{\lambda^\vee,\gamma,\nu,\delta}$-puzzle.
    Then we can count the number of $\inlinetri_{\nu,\sort(\beta)\gamma,\delta\beta}$-puzzles as 
\begin{align*}
\#\{\inlinetri_{\nu,\sort(\beta)\gamma,\delta\beta} \text{-puzzles}\} = \sum_\lambda \left(\#\{\inlinetri_{\alpha,\lambda,\beta} \text{-puzzles}\}\right) \left(\#\{\inlinetrapside_{\lambda^\vee,\gamma,\nu,\delta} \text{-puzzles}\}\right)  \\
= \sum_\lambda c_{\alpha,\lambda}^{\beta^\vee}\cdot \left(\#\{\inlinetrapside_{\lambda^\vee,\gamma,\nu,\delta} \text{-puzzles}\}\right) .
\end{align*}
Since $c_{\alpha,\lambda}^{\beta^\vee} = c_{\beta,\lambda}^{\alpha^\vee}$ (by commutativity plus rotational symmetry), this implies that we can commute $\alpha$ and $\beta$ as in the theorem statement.

To prove that $\gamma$ and $\delta$ commute, we first apply the commutative property of triangular puzzles to the NW and NE sides of the outer boundary to get
$$\#\{\inlinetri_{\alpha\gamma,\nu,\delta\beta} \text{-puzzles}\}= \#\{\inlinetri_{\nu,\alpha\gamma,\delta\beta} \text{-puzzles}\}.$$
Then we can similarly divide the puzzle into a triangular subpuzzle and a trapezoidal subpuzzle as shown in Figure \ref{fig:splitsymmetrypuzzlecommuted}, and we obtain 
\begin{align*}
\#\{\inlinetri_{\nu,\alpha\gamma,\delta\beta} \text{-puzzles}\} = \sum_\mu \left(\#\{\inlinetrapsiderotatetwo_{\nu,\alpha,\mu^\vee,\beta} \text{-puzzles}\}\right) \left(\#\{\inlinetri_{\mu,\gamma,\delta} \text{-puzzles}\}\right) \\
= \sum_\mu \left(\#\{\inlinetrapsiderotatetwo_{\nu,\alpha,\mu^\vee,\beta} \text{-puzzles}\}\right) \cdot c_{\mu,\gamma}^{\delta^\vee}.
\end{align*}
Since $c_{\mu,\gamma}^{\delta^\vee} = c_{\mu,\delta}^{\gamma^\vee}$ (by commutativity plus rotational symmetry), this implies that we can commute $\gamma$ and $\delta$ as in the theorem statement.

To make both swaps, we apply these arguments in succession.

\begin{figure*}[]
    \centering
    \begin{subfigure}[t]{0.47\textwidth}
    \centering
     \input{proof_boundary_diagrams/splitsymmetrypuzzle}
    \caption{$\alpha$ and $\beta$ commute via the basic commutative property of triangular puzzles (plus rotational symmetry) applied to the green triangular subpuzzle.}
    \label{fig:splitsymmetrypuzzle}   
    \end{subfigure}
    ~
    \begin{subfigure}[t]{0.47\textwidth}
    \centering
    \input{proof_boundary_diagrams/splitsymmetrypuzzlecommuted}
    \caption{$\gamma$ and $\delta$ commute via the basic commutative property of triangular puzzles (plus rotational symmetry) applied to the green triangular subpuzzle.}
    \label{fig:splitsymmetrypuzzlecommuted}
    \end{subfigure}
    \caption{ 
    }
    \label{fig:splitsymmetrypuzzles}
\end{figure*}

\item \textbf{Proof by the commutative property of triangular puzzles and the commutative property of trapezoidal puzzles.} Once the latter property is proved in the next section, we will be able to use it to commute $\gamma$ and $\delta$ directly in Figure \ref{fig:splitsymmetrypuzzle}. We use the technique from the previous argument to commute $\alpha$ and $\beta$. 

\end{enumerate} 
\end{proof}

\subsubsection{K-theory version}
\begin{remark}
\label{remark:splitsymmetryKtheory}
We have an analogous statement of Proposition \ref{prop:splitLRnumbers} for K-theory structure constants in the Schubert structure sheaf basis, and thus the analogous statement of Theorem \ref{thm:puzzlecommutesplit} for puzzles using pieces in $\Hfrak\cup\left\{\raisebox{-.4\height}{
\begin{tikzpicture}[x=1pt,y=1pt,yscale=1,xscale=1,scale=.68]
\draw [fill={rgb, 255:red, 255; green, 255; blue, 255 }] (0.0,0.0) -- (15.0,25.98) -- (30.0,0.0) -- cycle;
\draw (7.5,12.99) node  [rotate=60]  {$\mathtt{1\! 0}$};
\draw (22.5,12.99) node  [rotate=-60]  {$\mathtt{1\! 0}$};
\draw (15.0,0.0) node  [rotate=0]  {$\mathtt{1\! 0}$};
\end{tikzpicture}
}\right\}$ (where now the equalities are of the sums of the puzzle weights rather than just the number of puzzles). This comes from using the analogous statement of Lemma \ref{lemma:commutesplit} for K-theory (see Remark \ref{remark:lemmatwoktheory}) to write $$\underset{\Gr_{a+c}}{\Kint} [\Ocal_{\alpha\gamma}][\Ocal_{\delta\beta}][\Ical^{\nu^\vee}] = \underset{\Gr_{a+c}}{\Kint} [\Ocal_{\beta\gamma}][X_{\delta\alpha}][\Ical^{\nu^\vee}] = \underset{\Gr_{a+c}}{\Kint} [\Ocal_{\alpha\delta}][\Ocal_{\gamma\beta}][\Ical^{\nu^\vee}] = \underset{\Gr_{a+c}}{\Kint} [\Ocal_{\beta\delta}][\Ocal_{\gamma\alpha}][\Ical^{\nu^\vee}].$$
Note that it was necessary to commute $\nu$ and $\delta\beta$ first, in contrast to the statements and proofs in ordinary cohomology, since the K-theory Schubert bases are not self-dual.
\end{remark}


\section{Trapezoids}

We will see that a trapezoidal puzzle can be thought of as a special case of a triangular puzzle with split symmetry. Also, for most of the boundary shapes in Figure \ref{fig:boundaryshapescommute}, the puzzle can be split up into a trapezoidal subpuzzle and a trianglar subpuzzle or two trapezoidal subpuzzles, and our trapezoid commutativity property will provide an easy proof of the commutativity in those cases where a trapezoid is a building block. 



For the entirety of this section, let $a_0,a_1,c_0,c_1$ be non-negative integers, and define $a:=a_0+a_1$ and $c:=c_0+c_1$. Let us also shorten the notation for the spaces we will work with by defining $\Gr_{a+c}:=\Gr(a_1+c_1;\C^{a+c})$, $\Gr_a:=\Gr(a_1;F_a)$, and $\Gr_c:=\Gr(c_1;\Ftilde_c)$, where $\flag{F}$ and $\flag{\Ftilde}$ are the standard and anti-standard complete flags.

\subsection{Geometric interpretation of trapezoidal puzzles}

\subsubsection{Completion of a trapezoid to a triangle}

Now we will see our first instance of the operation of ``completing to a triangle,'' where we take a puzzle with a more complicated boundary shape and trivially glue on unique identity puzzles (see property 1 of Section \ref{section:otherpuzzleproperties}) to turn it into a triangular puzzle, which has a classical geometric Schubert calculus interpretation. 

\begin{proposition}
    \label{prop:trapezoidcomplete}
    For $\beta\in\binom{[a]}{a_1}$, $\gamma,\delta\in\binom{[c]}{c_1}$, and $\nu\in\binom{[c+a]}{a_1+c_1}$, and puzzle pieces in $\Hfrak \cup \left\{\raisebox{-.4\height}{
\begin{tikzpicture}[x=1pt,y=1pt,yscale=1,xscale=1,scale=.68]
\draw [fill={rgb, 255:red, 255; green, 246; blue, 74 }] (0.0,0.0) -- (15.0,25.98) -- (30.0,0.0) -- (15.0,-25.98) -- cycle;
\draw (7.5,12.99) node  [rotate=60]  {$\mathtt{0}$};
\draw (22.5,12.99) node  [rotate=-60]  {$\mathtt{1}$};
\draw (22.5,-12.99) node  [rotate=60]  {$\mathtt{0}$};
\draw (7.5,-12.99) node  [rotate=-60]  {$\mathtt{1}$};
\end{tikzpicture}
},\raisebox{-.4\height}{
\begin{tikzpicture}[x=1pt,y=1pt,yscale=1,xscale=1,scale=.68]
\draw [fill={rgb, 255:red, 255; green, 255; blue, 255 }] (0.0,0.0) -- (15.0,25.98) -- (30.0,0.0) -- cycle;
\draw (7.5,12.99) node  [rotate=60]  {$\mathtt{1\! 0}$};
\draw (22.5,12.99) node  [rotate=-60]  {$\mathtt{1\! 0}$};
\draw (15.0,0.0) node  [rotate=0]  {$\mathtt{1\! 0}$};
\end{tikzpicture}
} \right\}$, we have a bijection
    $$\{\inlinetrapside_{\beta,\gamma,\nu,\delta}\text{-puzzles}\} \leftrightarrow \{\inlinetri_{\sort(\beta)\gamma,\nu,\delta\beta} \text{-puzzles}\}.$$
\end{proposition}

\begin{figure*}[h!]
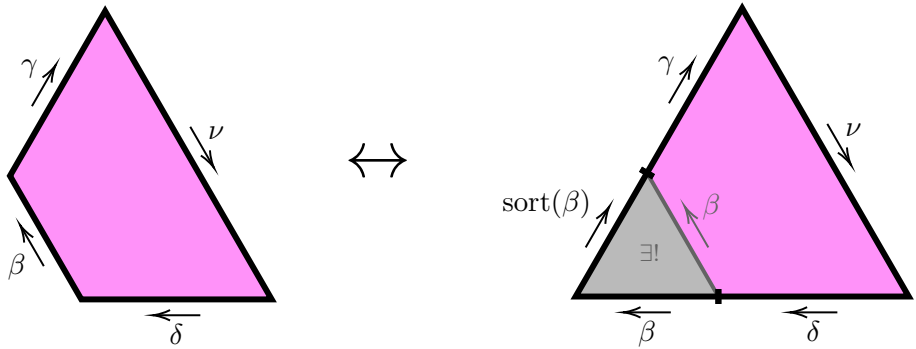

    \centering
    \raisebox{-.5\height}{\input{proof_boundary_diagrams/complete_to_triangle/trapezoidfilling}} 
    \Huge $\quad \leftrightarrow \quad$ \normalsize \raisebox{-.5\height}{\input{proof_boundary_diagrams/complete_to_triangle/trapezoidcomplete}}
    \caption{Puzzles with the boundary $\protect\inlinetrapside_{\beta,\gamma,\nu,\delta}$ are in bijection with puzzles with the boundary $\protect\inlinetri_{\sort(\beta)\gamma,\nu,\delta\beta}$. On the righthand side, given the outer boundary labels, there exists a unique filling of the grey triangular region, which forces the labels around the inner pink region to replicate those for the trapezoid itself (see Property 1 of Section \ref{section:otherpuzzleproperties}). The bijection is to glue on the uniquely filled grey region, and the inverse is to cut it off.}
    \label{fig:trapezoidcomplete}
\end{figure*}

\begin{proof}
    See Figure \ref{fig:trapezoidcomplete}.
\end{proof}

\begin{corollary}
\label{cor:trapezoidpushforwardcoeffs}
    A special case of trapezoidal puzzles compute the coefficients of Proposition \ref{prop:pushforwardcoeff}. Namely, using puzzle pieces in $\Hfrak$, we have
    $$\#\{\inlinetrapside_{\lambda^\vee,\sort(\mu),\nu,\mu^\vee} \text{-puzzles}\} = d_{\lambda,\mu}^\nu.$$
\end{corollary}
\begin{proof}
     Combining Proposition \ref{prop:trapezoidcomplete}, Proposition \ref{prop:equalityofLR}, and Corollary \ref{cor:pushforwardcoeffs}, we have
    $$\#\{\inlinetrapside_{\lambda^\vee,\sort(\mu),\nu,\mu^\vee} \text{-puzzles}\} = \#\{\inlinetri_{\sort(\lambda)\sort(\mu),\nu,(\lambda\mu)^\vee} \text{-puzzles}\} = c_{\sort(\lambda)\sort(\mu),\nu}^{\lambda\mu} = c_{\lambda\pad,\mu\pad}^{\nu}=d_{\lambda,\mu}^\nu.$$
\end{proof}

\subsubsection{A formula in terms of Littlewood-Richardson numbers}

\begin{proposition}
    \label{prop:trapezoidgeometric}
    For $\beta\in\binom{[a]}{a_1}$, $\gamma,\delta\in\binom{[c]}{c_1}$, and $\nu\in\binom{[c+a]}{a_1+c_1}$, and using puzzle pieces in $\Hfrak$, we have
    $$\#\{\inlinetrapside_{\beta,\gamma,\nu,\delta}\text{-puzzles}\} =  \sum_{\mu}c_{(\beta^\vee)\pad,\mu\pad}^\nu \cdot c_{\gamma,\mu}^{\delta^\vee}.$$
\end{proposition}

\begin{proof}
    After completing a trapezoidal puzzle to a triangular puzzle as in Proposition \ref{prop:trapezoidcomplete}, we obtain a special case of a puzzle with split symmetry. We apply Propositions \ref{prop:trapezoidcomplete} and \ref{prop:splitsymmetryformula}, setting $\alpha=\sort(\beta)$, to write

     $$\#\{\inlinetrapside_{\beta,\gamma,\nu,\delta}\text{-puzzles}\} = \#\{\inlinetri_{\sort(\beta)\gamma,\nu,\delta\beta} \text{-puzzles}\} = c_{\sort(\beta)\gamma,\nu}^{(\delta\beta)^\vee} =   \sum_{\lambda,\mu}c_{\lambda\pad,\mu\pad}^\nu \cdot c_{\sort(\beta),\lambda}^{\beta^\vee}\cdot c_{\gamma,\mu}^{\delta^\vee}.$$

Now noting that $[X_{\sort(\beta)}]=[\Gr_a]=1\in H^*(\Gr_a)$ since the string $\sort(\beta)$ encodes no essential conditions, we write
     $$c_{\sort(\beta),\lambda}^{\beta^\vee} = \underset{\Gr_{a}}{\int} [X_{\sort(\beta)}][X^{\beta^\vee}][X_\lambda]= \underset{\Gr_{a}}{\int}[X^{\beta^\vee}][X_\lambda] = \langle \lambda,\beta^\vee\rangle = \begin{cases}1, &\text{ if } \lambda=\beta^\vee \\0, &\text{ if }\lambda\neq\beta^\vee \end{cases}.$$
     Therefore the summand above equals $0$ unless $\lambda=\beta^\vee$, so we reduce to the cases where it is nonzero and rewrite our expression as 
    $$c_{\sort(\beta)\gamma,\nu}^{(\delta\beta)^\vee} =   \sum_{\lambda,\mu}c_{\lambda\pad,\mu\pad}^\nu \cdot c_{\sort(\beta),\lambda}^{\beta^\vee}\cdot c_{\gamma,\mu}^{\delta^\vee} = \sum_{\mu}c_{(\beta^\vee)\pad,\mu\pad}^\nu \cdot c_{\gamma,\mu}^{\delta^\vee}.$$
\end{proof}

\subsection{Commutative property of trapezoidal puzzles}

\begin{figure}[h]
    \centering
    \input{boundaries_commute/trapezoidboundarycommute}
    \caption{$\protect\inlinetrapside_{\beta,\gamma,\nu,\delta}$}
    \label{fig:trapezoidboundarycommutesolo}
\end{figure}

\begin{theorem}[Commutative Property of Trapezoidal Puzzles]
\label{thm:trapezoidcommute}
     For $\beta\in\binom{[a]}{a_1}$, $\gamma,\delta\in\binom{[c]}{c_1}$, and $\nu\in\binom{[c+a]}{a_1+c_1}$, and using puzzle pieces in $\Hfrak$, we have
    $$\#\{\inlinetrapside_{\beta,\gamma,\nu,\delta} \text{-puzzles}\} = \#\{\inlinetrapside_{\beta,\delta,\nu,\gamma} \text{-puzzles}\}.$$
    In other words we can commute the labels on the equal-length sides of the trapezoid and preserve the number of puzzles.
\end{theorem}

\begin{proof}
We give several different proof methods. 
\begin{enumerate}
    \item \textbf{Proof by the commutative property of triangular puzzles.} We use the same method as in Argument 2 in the proof of Theorem \ref{thm:puzzlecommutesplit} to commute $\gamma$ and $\delta$.

    \item \textbf{Proof by the commutative property of puzzles with split symmetry.} By completing to a triangle, we may think of our trapezoidal puzzle as a $\inlinetri_{\sort(\beta)\gamma,\nu,\delta\beta}$-puzzle, which has split symmetry. We apply Theorem \ref{thm:puzzlecommutesplit}, and we have 
    $$\#\{\inlinetri_{\sort(\beta)\gamma,\nu,\delta\beta} \text{-puzzles}\}= \#\{\inlinetri_{\sort(\beta)\delta,\nu,\gamma\beta} \text{-puzzles}\}.$$
    The claim follows immediately.

    \item \textbf{Geometric proof.} We could prove this theorem directly, geometrically, in exactly the same way as we did Proposition \ref{prop:splitLRnumbers} and Theorem \ref{thm:puzzlecommutesplit}, just setting $\nu=\sort(\beta)$.

    \item \textbf{Geometric proof.} This theorem is also immediate from Proposition \ref{prop:trapezoidgeometric}, since $c_{\gamma,\mu}^{\delta^\vee}=c_{\delta,\mu}^{\gamma^\vee}$.

\explainalot{
    \item \textbf{Proof by the commutative property of triangular puzzles.} It can also be proved by using the basic commutative property of triangular puzzles. We complete the trapezoid to a triangle as in Figure \ref{fig:trapezoidcomplete}, then we apply the commutative property to the NW and NE sides of the triangular boundary to get
    $$\#\{\inlinetri_{\sort(\beta)\gamma,\nu,\delta\beta} \text{-puzzles}\}= \#\{\inlinetri_{\nu,\sort(\beta)\gamma,\delta\beta} \text{-puzzles}\}.$$
    
    Assume we have a $\inlinetri_{\nu,\sort(\beta)\gamma,\delta\beta}$-puzzle. Then because $\gamma$ and $\delta$ have matching content, there can be no $\tenstr$s appearing along the red line shown in Figure \ref{fig:trapezoidcommuted}, and if we call the string appearing along that line $\mu$, then $\mu$ must have the same content as $\gamma$ and $\delta$. (This comes from statement (a) under the discrete Green's theorem property given in Section \ref{section:otherpuzzleproperties}.) So we can divide the puzzle into a $\inlinetri_{\mu,\gamma,\delta}$-puzzle and a trapezoidal $\inlinetrapsiderotatetwo_{\nu,\sort(\beta),\mu^\vee,\beta}$-puzzle.
    
\begin{figure*}[h]
       \centering 
        \input{proof_boundary_diagrams/trapezoidcommuted}
        \caption{In any $\protect\inlinetri_{\nu,\sort(\beta)\gamma,\delta\beta}$-puzzle, there must be a triangular subpuzzle filling the green region, as well as a trapezoidal subpuzzle filling the yellow region.}
       \label{fig:trapezoidcommuted}
    \end{figure*}
    
Then we can count the number of $\inlinetri_{\nu,\sort(\beta)\gamma,\delta\beta}$-puzzles as 
\begin{align*}
\#\{\inlinetri_{\nu,\sort(\beta)\gamma,\delta\beta} \text{-puzzles}\} = \sum_\mu \left(\#\{\inlinetrapsiderotatetwo_{\nu,\sort(\beta),\mu^\vee,\beta} \text{-puzzles}\}\right) \left(\#\{\inlinetri_{\mu,\gamma,\delta} \text{-puzzles}\}\right) \\
= \sum_\mu \left(\#\{\inlinetrapsiderotatetwo_{\nu,\sort(\beta),\mu^\vee,\beta} \text{-puzzles}\}\right) \cdot c_{\mu,\gamma}^{\delta^\vee}.
\end{align*}
Since $c_{\mu,\gamma}^{\delta^\vee} = c_{\mu,\delta}^{\gamma^\vee}$, this implies that we can commute $\gamma$ and $\delta$ as in the theorem statement.
}

\end{enumerate}
\end{proof}

\begin{remark}
    If the pair $\gamma,\delta$ were not specified to have matching content in Theorem \ref{thm:trapezoidcommute}, the statement would still hold true in the case of mismatched content, since there would exist no puzzles in that case. This is a consequence of the discrete Green's theorem, statement (c), in Section \ref{section:otherpuzzleproperties}.
\end{remark}


\subsubsection{K-theory version}
\begin{remark}
\label{remark:Ktheorytrapezoidpuzzles}
    Theorem \ref{thm:trapezoidcommute} has an analogous statement for K-theory puzzles using pieces in $\Hfrak\cup\left\{\raisebox{-.4\height}{
\begin{tikzpicture}[x=1pt,y=1pt,yscale=1,xscale=1,scale=.68]
\draw [fill={rgb, 255:red, 255; green, 255; blue, 255 }] (0.0,0.0) -- (15.0,25.98) -- (30.0,0.0) -- cycle;
\draw (7.5,12.99) node  [rotate=60]  {$\mathtt{1\! 0}$};
\draw (22.5,12.99) node  [rotate=-60]  {$\mathtt{1\! 0}$};
\draw (15.0,0.0) node  [rotate=0]  {$\mathtt{1\! 0}$};
\end{tikzpicture}
}\right\}$, where the sum of the weights of the trapezoidal puzzles (rather than just the number) on either side of the commuting operation is equal. Arguments 1, 2, and 3 (given Remark \ref{remark:splitsymmetryKtheory}) work in a nearly identical way for K-theory as well.
\end{remark}

\section{Parallelograms}

We will see that, once given a geometric interpretation, parallelogram-shaped puzzles can be viewed as a special case of trapezoidal puzzles. The commutative property of parallelogram-shaped puzzles will be implied by that of trapezoidal puzzles, but the parallelogram as a special case allows us to take our results further than with the trapezoid. In Theorem \ref{thm:parallelogramLRnumbers}, we get a more detailed formula, in terms of Littlewood-Richardson numbers, for the number of parallelogram-shaped puzzles. We include both a purely geometric proof of that theorem, as well as a puzzle-based proof.

With Theorem \ref{thm:eqvtparallelogram}, we are also able to prove a commutativity result for parallelogram-shaped equivariant puzzles, where we state how commuting a pair of opposite boundary labels acts on the associated structure constant in the $T$-equivariant cohomology. 

For the entirety of this section, let $a_0,a_1,c_0,c_1$ be non-negative integers, and define $a:=a_0+a_1$ and $c:=c_0+c_1$. Let $\flag{F}$ and $\flag{\Ftilde}$ be the standard and anti-standard flags in $\C^{a+c}$, respectively. Let us also shorten the notation for the spaces we will work with by defining $\Gr_{a+c}:=\Gr(a_1+c_1;\C^{a+c})$, $\Gr_a:=\Gr(a_1;F_a)$, and $\Gr_c:=\Gr(c_1;\Ftilde_c)$.

\subsection{Geometric interpretation of parallelogram-shaped puzzles}

\subsubsection{Completion of a parallelogram to a triangle}

\begin{proposition}
    \label{prop:parallelogramcomplete}
    For $\alpha,\beta\in\binom{[a]}{a_1}$ and $\gamma,\delta\in\binom{[c]}{c_1}$, and using puzzle pieces in $\Hfrak \cup \left\{\raisebox{-.4\height}{
\begin{tikzpicture}[x=1pt,y=1pt,yscale=1,xscale=1,scale=.68]
\draw [fill={rgb, 255:red, 255; green, 246; blue, 74 }] (0.0,0.0) -- (15.0,25.98) -- (30.0,0.0) -- (15.0,-25.98) -- cycle;
\draw (7.5,12.99) node  [rotate=60]  {$\mathtt{0}$};
\draw (22.5,12.99) node  [rotate=-60]  {$\mathtt{1}$};
\draw (22.5,-12.99) node  [rotate=60]  {$\mathtt{0}$};
\draw (7.5,-12.99) node  [rotate=-60]  {$\mathtt{1}$};
\end{tikzpicture}
},\raisebox{-.4\height}{
\begin{tikzpicture}[x=1pt,y=1pt,yscale=1,xscale=1,scale=.68]
\draw [fill={rgb, 255:red, 255; green, 255; blue, 255 }] (0.0,0.0) -- (15.0,25.98) -- (30.0,0.0) -- cycle;
\draw (7.5,12.99) node  [rotate=60]  {$\mathtt{1\! 0}$};
\draw (22.5,12.99) node  [rotate=-60]  {$\mathtt{1\! 0}$};
\draw (15.0,0.0) node  [rotate=0]  {$\mathtt{1\! 0}$};
\end{tikzpicture}
} \right\}$, we have a bijection
    $$\{\inlinepar_{\alpha,\gamma,\beta,\delta} \text{-puzzles}\} \leftrightarrow \{\inlinetri_{\sort(\alpha)\gamma,\beta\sort(\delta),\delta\alpha} \text{-puzzles}\}.$$
\end{proposition}

\begin{figure*}[h]
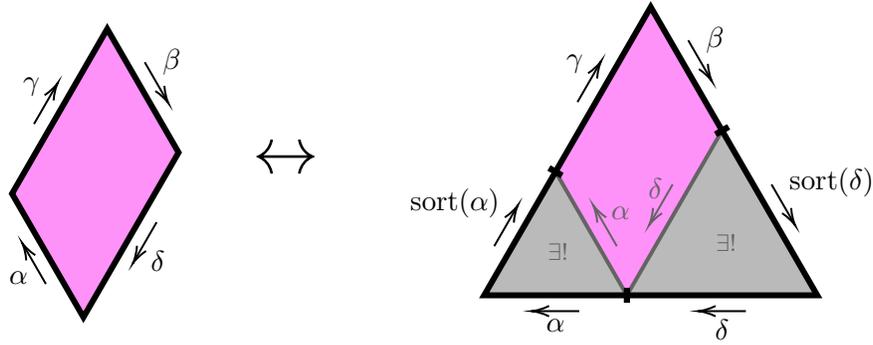

    \centering
    \raisebox{-.5\height}{\input{proof_boundary_diagrams/complete_to_triangle/parallelogramfilling}} 
    \Huge $\quad \leftrightarrow \quad$ \normalsize \raisebox{-.5\height}{\input{proof_boundary_diagrams/complete_to_triangle/parallelogramcomplete}}
    \caption{Puzzles with the boundary $\protect\inlinepar_{\alpha,\gamma,\beta,\delta}$ are in bijection with puzzles with the boundary $\protect\inlinetri_{\sort(\alpha)\gamma,\beta\sort(\delta),\delta\alpha}$. 
    On the righthand side, given the outer boundary labels, there exists a unique filling of each grey triangular region, and these force the labels around the inner pink region to replicate those for the parallelogram itself (see Property 1 of Section \ref{section:otherpuzzleproperties}). The bijection is to glue on the uniquely filled grey regions, and the inverse is to cut them off.}
    \label{fig:parallelogramcomplete}
\end{figure*}

\begin{proof}
    See Figure \ref{fig:parallelogramcomplete}.
\end{proof}

\subsubsection{A formula in terms of Littlewood-Richardson numbers}

\begin{theorem}
    \label{thm:parallelogramLRnumbers}

    For $\alpha,\beta\in\binom{[a_0+a_1]}{a_1}$ and $\gamma,\delta\in\binom{[c_0+c_1]}{c_1}$, and using puzzle pieces in $\Hfrak$, we have 
    $$\#\{\inlinepar_{\alpha,\gamma,\beta,\delta} \text{-puzzles}\} =  \begin{cases}
        \sum\limits_{\mu}  c_{\alpha,\beta}^{\mu\onestr^{a_1-c_1}\zerostr^{a_0-c_0}}\cdot c_{\gamma,\delta}^{\mu^\vee} & \text{if } a_0\geq c_0 \text{ and } a_1>c_1 \\
        \sum\limits_{\lambda}  c_{\alpha,\beta}^{\lambda^\vee}\cdot c_{\gamma,\delta}^{\onestr^{c_1-a_1}\zerostr^{c_0-a_0}\lambda} & \text{if } c_1\geq a_1 \text{ and } c_0>a_0 \\
        \sum\limits_{\eta}  c_{\alpha,\beta}^{\eta^\vee\zerostr^{a_0-c_0}}\cdot c_{\gamma,\delta}^{\onestr^{c_1-a_1}\eta} & \text{if } a_0\geq c_0 \text{ and } c_1\geq a_1 \\
        0 & \text{if } c_0>a_0 \text{ and } a_1>c_1
    \end{cases}.$$
\end{theorem}

\begin{proof}

In this section we give a geometric proof. We will give an alternate puzzle-based proof in Section \ref{section:uniquefillingproof}.

\textbf{Geometric proof.}
By completing to a triangle as in Proposition \ref{prop:parallelogramcomplete}, we start with
$$\#\{\inlinepar_{\alpha,\gamma,\beta,\delta} \text{-puzzles}\} =  \#\{\inlinetri_{\sort(\alpha)\gamma,\beta\sort(\delta),\delta\alpha} \text{-puzzles}\}  = \underset{\Gr_{a+c}}{\int} [X_{\sort(\alpha)\gamma}][X_{\beta\sort(\delta)}][X^{(\delta\alpha)^\vee}].$$

Next we will use Lemma \ref{lemma:commutesplit} to write $[X_{\sort(\alpha)\gamma}][X^{(\delta\alpha)^\vee}]= [X_{\alpha\gamma}][X^{(\delta\sort(\alpha))^\vee}]$ and $[X_{\beta\sort(\delta)}][X^{(\delta\sort(\alpha))^\vee}]= [X_{\sort(\alpha)\sort(\delta)}][X^{(\delta\beta)^\vee}]$. The purpose of this is to commute labels so that the problem looks like the situation of split symmetry rather than a case of a trapezoid, which happens to help with this proof. We have
\begin{align*}
    \#\{\inlinepar_{\alpha,\gamma,\beta,\delta} \text{-puzzles}\} &= \underset{\Gr_{a+c}}{\int} [X_{\sort(\alpha)\gamma}][X_{\beta\sort(\delta)}][X^{(\delta\alpha)^\vee}]= \underset{\Gr_{a+c}}{\int}  [X_{\alpha\gamma}][X_{\beta\sort(\delta)}][X^{(\delta\sort(\alpha))^\vee}]\\ 
    &= \underset{\Gr_{a+c}}{\int}  [X_{\alpha\gamma}][X_{\sort(\alpha)\sort(\delta)}][X^{(\delta\beta)^\vee}] =c_{\alpha\gamma,\sort(\alpha)\sort(\delta)}^{(\delta\beta)^\vee}.
\end{align*}

Now this Littlewood-Richardson number corresponds to puzzles with split symmetry, so we can apply Proposition \ref{prop:splitsymmetryformula}, setting $\nu=\sort(\alpha)\sort(\delta)$, and write
$$\#\{\inlinepar_{\alpha,\gamma,\beta,\delta} \text{-puzzles}\} =c_{\alpha\gamma,\sort(\alpha)\sort(\delta)}^{(\delta\beta)^\vee} = \sum_{\lambda,\mu}c_{\lambda\pad,\mu\pad}^{\sort(\alpha)\sort(\delta)} \cdot c_{\alpha,\lambda}^{\beta^\vee}\cdot c_{\gamma,\mu}^{\delta^\vee} = \sum_{\lambda,\mu}c_{\lambda\pad,\mu\pad}^{\sort(\alpha)\sort(\delta)} \cdot c_{\alpha,\beta}^{\lambda^\vee}\cdot c_{\gamma,\delta}^{\mu^\vee}.$$
 
Let us look more closely at $c_{\lambda\pad,\mu\pad}^{\sort(\alpha)\sort(\delta)}$. Recall that here $\lambda\pad=\zerostr^{c_0}\lambda \onestr^{c_1}$, $\mu\pad=\zerostr^{a_0}\mu \onestr^{a_1}$, and $\sort(\alpha)\sort(\delta)= \zerostr^{a_0}\onestr^{a_1}\zerostr^{c_0} \onestr^{c_1}$. 
We have
\begin{align*}
    c_{\lambda\pad,\mu\pad}^{\sort(\alpha)\sort(\delta)} &= \underset{\Gr_{a+c}}{\int} [X_{\lambda\pad}][X_{\mu\pad}][X^{\sort(\alpha)\sort(\delta)}] \\
    &=\underset{\Gr_{a+c}}{\int}[X_{\zerostr^{c_0}\lambda \onestr^{c_1}}][X_{\zerostr^{a_0}\mu \onestr^{a_1}}][X^{\zerostr^{a_0}\onestr^{a_1}\zerostr^{c_0} \onestr^{c_1}}]\\
    &= \underset{\Gr_{a+c}}{\int} [X_{\zerostr^{c_0}\lambda \onestr^{c_1}}(\flag{F})][X_{\zerostr^{a_0}\mu \onestr^{a_1}}(\flag{F})][X_{\onestr^{c_1}\zerostr^{c_0}\onestr^{a_1}\zerostr^{a_0}}(\flag{\Ftilde})].
\end{align*}


We will reduce to cases where the integral is nonzero. First of all, if $\zerostr^{c_0}\lambda \onestr^{c_1}$ does not begin with $\zerostr^{a_0}$, then we have $[X_{\zerostr^{c_0}\lambda \onestr^{c_1}}(\flag{F})][X_{\onestr^{c_1}\zerostr^{c_0}\onestr^{a_1}\zerostr^{a_0}}(\flag{\Ftilde})]= [X_{\zerostr^{c_0}\lambda \onestr^{c_1}}(\flag{F})\cap X_{\onestr^{c_1}\zerostr^{c_0}\onestr^{a_1}\zerostr^{a_0}}(\flag{\Ftilde})]=[\emptyset]=0$. This is because the conditions $\onestr^{c_1}\zerostr^{c_0}\onestr^{a_1}\zerostr^{a_0}$ say that for  $V\in X_{\onestr^{c_1}\zerostr^{c_0}\onestr^{a_1}\zerostr^{a_0}}(\flag{\Ftilde})$, we have $V\subseteq \Ftilde_{c_1+c_0+a_1}$, so then $V\cap F_{a_0}=\{\mathbf{0}\}$. Then if $\zerostr^{c_0}\lambda \onestr^{c_1}$ has a $\onestr$ anywhere in the first $a_0$ spots, having $V\in X_{\zerostr^{c_0}\lambda \onestr^{c_1}}(\flag{F})$ would require $\dim(V\cap F_{a_0})\geq 1$, which is impossible if $V\in X_{\onestr^{c_1}\zerostr^{c_0}\onestr^{a_1}\zerostr^{a_0}}(\flag{\Ftilde})$. 
Therefore, let us assume that $\zerostr^{c_0}\lambda \onestr^{c_1}$ does begin with $\zerostr^{a_0}$, so we can rewrite $\zerostr^{c_0}\lambda \onestr^{c_1}=\zerostr^{a_0}\eta\onestr^{c_1}$, where $\eta=\zerostr^{c_0-a_0}\lambda$ if $c_0\geq a_0$ and $\zerostr^{a_0-c_0}\eta=\lambda$ if $c_0< a_0$.

In addition, if $\zerostr^{a_0}\mu \onestr^{a_1}$ does not end with $\onestr^{c_1}$, then $[X_{\zerostr^{a_0}\mu \onestr^{a_1}}(\flag{F})][X_{\onestr^{c_1}\zerostr^{c_0}\onestr^{a_1}\zerostr^{a_0}}(\flag{\Ftilde})]= [X_{\zerostr^{a_0}\mu \onestr^{a_1}}(\flag{F})\cap X_{\onestr^{c_1}\zerostr^{c_0}\onestr^{a_1}\zerostr^{a_0}}(\flag{\Ftilde})]=[\emptyset]=0$. This is because the conditions $\onestr^{c_1}\zerostr^{c_0}\onestr^{a_1}\zerostr^{a_0}$ say that for  $V\in X_{\onestr^{c_1}\zerostr^{c_0}\onestr^{a_1}\zerostr^{a_0}}(\flag{\Ftilde})$, we have $\Ftilde_{c_1}\subseteq V$, so $\dim(V\cap \Ftilde_{c_1})=c_1$. But then if $\zerostr^{a_0}\mu \onestr^{a_1}$ has a $\zerostr$ anywhere in the last $c_1$ spots, having $V\in X_{\zerostr^{a_0}\mu \onestr^{a_1}}(\flag{F})$ would mean $\dim(V\cap F_{a_0+a_1+c_0})>a_1$. Since $F_{a_0+a_1+c_0}\cap \Ftilde_{c_1}=\{\mathbf{0}\}$ and $\dim(V)=a_1+c_1$, we cannot have $\dim(V\cap F_{a_0+a_1+c_0})>a_1$ and $\dim(V\cap \Ftilde_{c_1})=c_1$ simultaneously. Therefore, let us assume that $\zerostr^{a_0}\mu \onestr^{a_1}$ does end with $\onestr^{c_1}$, so we can rewrite $\zerostr^{a_0}\mu \onestr^{a_1}=\zerostr^{a_0}\theta\onestr^{c_1}$, where $\theta\onestr^{c_1-a_1}=\mu$ if $c_1\geq a_1$ and $\theta=\mu\onestr^{a_1-c_1}$ if $c_1< a_1$.



Now define the maps 
\begin{gather*}
    \Psi: \Gr(0,F_{a_0})\times \Gr(a_1,\langle \bfe_{a_0+1},\ldots,\bfe_{a_0+a_1+c_0}\rangle) \times \Gr(c_1;\Ftilde_{c_1}) \hookrightarrow \Gr(a_1,F_{a_0+a_1+c_0}) \times \Gr(c_1;\Ftilde_{c_1}), \\ 
(U,V,W)=(\{\mathbf{0}\},V,\Ftilde_{c_1})\mapsto (U\oplus V,W) =(\{\mathbf{0}\} \oplus V, \Ftilde_{c_1})
\end{gather*}
and 
\begin{gather*}
\Omega:\Gr(a_1,F_{a_0+a_1+c_0}) \times \Gr(c_1;\Ftilde_{c_1}) \hookrightarrow \Gr(a_1+c_1;\C^{a+c}),\\
(V,W)= (V,\Ftilde_{c_1}) \mapsto V\oplus W= V\oplus \Ftilde_{c_1}.
\end{gather*}
Going forward, let us shorten the notation for these spaces to $\Gr_{a_0}:=\Gr(0,F_{a_0})$, $\Gr_{a_1+c_0}:=\Gr(a_1,\langle \bfe_{a_0+1},\ldots,\bfe_{a_0+a_1+c_0}\rangle)$, and $\Gr_{c_1}:=\Gr(c_1;\Ftilde_{c_1})$.



Using Lemma \ref{lemma:directsum}, we have $[X_{\zerostr^{a_0}\theta\onestr^{c_1}}][X^{\zerostr^{a_0}\onestr^{a_1}\zerostr^{c_0} \onestr^{c_1}}]=\Omega_*([X_{\zerostr^{a_0}\theta}][X^{\zerostr^{a_0}\onestr^{a_1}\zerostr^{c_0}}]\otimes [X_{\onestr^{c_1}}][X^{\onestr^{c_1}}])$, and $[X_{\zerostr^{a_0}\theta}][X^{\zerostr^{a_0}\onestr^{a_1}\zerostr^{c_0}}]= \Omega'_*([X_{\zerostr^{a_0}}][X^{\zerostr^{a_0}}]\otimes [X_{\theta}][X^{\onestr^{a_1}\zerostr^{c_0}}])$. 
Then $[X_{\zerostr^{a_0}\theta\onestr^{c_1}}][X^{\zerostr^{a_0}\onestr^{a_1}\zerostr^{c_0} \onestr^{c_1}}] = (\Omega_*\circ \Psi_*)([X_{\zerostr^{a_0}}][X^{\zerostr^{a_0}}]\otimes [X_{\theta}][X^{\onestr^{a_1}\zerostr^{c_0}}]\otimes [X_{\onestr^{c_1}}][X^{\onestr^{c_1}}])= (\Omega\circ \Psi)_*([X_{\zerostr^{a_0}}]\otimes [X_{\theta}]\otimes [X_{\onestr^{c_1}}])$.

We also want to find $(\Omega \circ \Psi)^*([X_{\zerostr^{a_0}\eta \onestr^{c_1}}])$, which will come up later. We can write it as a linear combination of basis elements $[X_{\zerostr^{a_0}}]\otimes [X_{\rho}]\otimes [X^{\onestr_{c_1}}]$ of $H^*(\Gr_{a_0}\times \Gr_{a_1+c_0}\times \Gr_{c_1})$. The coefficient on $[X_{\zerostr^{a_0}}]\otimes [X_{\rho}]\otimes [X_{\onestr^{c_1}}]$ is given by 
\begin{multline*}
    \underset{\Gr_{a_0}\times \Gr_{a_1+c_0}\times \Gr_{c_1}}{\int}  (\Omega \circ \Psi)^*([X_{\zerostr^{a_0}\eta\onestr^{c_1}}])([X^{\zerostr^{a_0}}]\otimes [X^{\rho}]\otimes [X^{\onestr^{c_1}}]) \\
    = \underset{\Gr_{a+c}}{\int} [X_{\zerostr^{a_0}\eta\onestr^{c_1}}](\Omega \circ \Psi)_*([X^{\zerostr^{a_0}}]\otimes [X^{\rho}]\otimes [X^{\onestr^{c_1}}]) = \underset{\Gr_{a+c}}{\int} [X_{\zerostr^{a_0}\eta\onestr^{c_1}}][X^{\zerostr^{a_0}\rho\onestr^{c_1}}]=\langle \eta,\rho\rangle.
\end{multline*}

This implies that $(\Omega \circ \Psi)^*([X_{\zerostr^{a_0}\eta\onestr^{c_1}}])= [X_{\zerostr^{a_0}}]\otimes [X_{\eta}]\otimes [X_{\onestr^{c_1}}]$. 

Now we return to our goal of analyzing $c_{\lambda\pad,\mu\pad}^{\sort(\alpha)\sort(\delta)}$.  Putting everything we established above together, we write
{\allowdisplaybreaks 
\begin{align*}
c_{\lambda\pad,\mu\pad}^{\sort(\alpha)\sort(\delta)} &= \underset{\Gr_{a+c}}{\int} [X_{\zerostr^{c_0}\lambda \onestr^{c_1}}][X_{\zerostr^{a_0}\mu \onestr^{a_1}}][X^{\zerostr^{a_0}\onestr^{a_1}\zerostr^{c_0} \onestr^{c_1}}]\\
&=\underset{\Gr_{a+c}}{\int} [X_{\zerostr^{a_0}\eta \onestr^{c_1}}][X_{\zerostr^{a_0}\theta \onestr^{c_1}}][X^{\zerostr^{a_0}\onestr^{a_1}\zerostr^{c_0} \onestr^{c_1}}]\\
&=\underset{\Gr_{a+c}}{\int} [X_{\zerostr^{a_0}\eta \onestr^{c_1}}](\Omega\circ \Psi)_*([X_{\zerostr^{a_0}}]\otimes [X_{\theta}]\otimes [X_{\onestr^{c_1}}])\\
&=\underset{\Gr_{a_0}\times \Gr_{a_1+c_0}\times \Gr_{c_1}}{\int} (\Omega \circ \Psi)^*([X_{\zerostr^{a_0}\eta \onestr^{c_1}}]) ([X_{\zerostr^{a_0}}]\otimes [X_{\theta}]\otimes [X_{\onestr^{c_1}}])\\
&=\underset{\Gr_{a_0}\times \Gr_{a_1+c_0}\times \Gr_{c_1}}{\int} ([X_{\zerostr^{a_0}}]\otimes [X_{\eta}]\otimes [X_{\onestr^{c_1}}]) ([X_{\zerostr^{a_0}}]\otimes [X_{\theta}]\otimes [X_{\onestr^{c_1}}])\\
&=\underset{\Gr_{a_0}\times \Gr_{a_1+c_0}\times \Gr_{c_1}}{\int} ([X_{\zerostr^{a_0}}][X_{\zerostr^{a_0}}]\otimes [X_{\eta}][X_{\theta}]\otimes [X_{\onestr^{c_1}}][X_{\onestr^{c_1}}])\\
&=\underset{\Gr_{a_0}\times \Gr_{a_1+c_0}\times \Gr_{c_1}}{\int} ([X_{\zerostr^{a_0}}][X^{\zerostr^{a_0}}]\otimes [X_{\eta}][X^{\theta^\vee}]\otimes [X_{\onestr^{c_1}}][X^{\onestr^{c_1}}]) \\
&=\underset{\Gr_{a_0}}{\int} [X_{\zerostr^{a_0}}][X^{\zerostr^{a_0}}] \underset{\Gr_{a_1+c_0}}{\int} [X_{\eta}][X^{\theta^\vee}] \underset{ \Gr_{c_1}}{\int} [X_{\onestr^{c_1}}][X^{\onestr^{c_1}}] \\
&=\underset{\Gr_{a_1+c_0}}{\int} [X_{\eta}][X^{\theta^\vee}]=\langle \eta,\theta^\vee\rangle.
\end{align*}}

This says that $c_{\lambda\pad,\mu\pad}^{\sort(\alpha)\sort(\delta)}= 1$ if $\eta=\theta^\vee$ and 0 otherwise. We will see what this condition means in terms of $\lambda$ and $\mu$ and finish writing the formula. This will depend on the values of $a_0,a_1,c_0,c_1$ relative to each other. The four possible cases are as follows.

\begin{enumerate}
    \item Case: $a_0\geq c_0$ and $a_1>c_1$. We have $\zerostr^{a_0-c_0}\eta=\lambda$ and $\theta=\mu\onestr^{a_1-c_1}$. Then $\eta=\theta^\vee$ if and only if $\lambda=\zerostr^{a_0-c_0}\eta= \zerostr^{a_0-c_0}\theta^\vee = \zerostr^{a_0-c_0}\onestr^{a_1-c_1}\mu^\vee$. We sum over $\mu$.
    $$\#\{\inlinepar_{\alpha,\gamma,\beta,\delta} \text{-puzzles}\} =\sum_{\lambda,\mu} c_{\lambda\pad,\mu\pad}^{\sort(\alpha)\sort(\delta)} \cdot c_{\alpha,\beta}^{\lambda^\vee}\cdot c_{\gamma,\delta}^{\mu^\vee} = \sum_{\mu}  c_{\alpha,\beta}^{\mu\onestr^{a_1-c_1}\zerostr^{a_0-c_0}}\cdot c_{\gamma,\delta}^{\mu^\vee}$$

    \item Case: $c_1\geq a_1$ and $c_0> a_0$. We have $\eta=\zerostr^{c_0-a_0}\lambda$ and $\theta\onestr^{c_1-a_1}=\mu$. Then $\eta=\theta^\vee$ if and only if $\mu=\theta\onestr^{c_1-a_1}= \eta^\vee\onestr^{c_1-a_1} = \lambda^\vee\zerostr^{c_0-a_0}\onestr^{c_1-a_1}$. We sum over $\lambda$.
    $$\#\{\inlinepar_{\alpha,\gamma,\beta,\delta} \text{-puzzles}\} =\sum_{\lambda,\mu} c_{\lambda\pad,\mu\pad}^{\sort(\alpha)\sort(\delta)} \cdot c_{\alpha,\beta}^{\lambda^\vee}\cdot c_{\gamma,\delta}^{\mu^\vee} = \sum_{\lambda}  c_{\alpha,\beta}^{\lambda^\vee}\cdot c_{\gamma,\delta}^{\onestr^{c_1-a_1}\zerostr^{c_0-a_0}\lambda}$$

    \item Case: $a_0\geq c_0$ and $c_1\geq a_1$. We have $\zerostr^{a_0-c_0}\eta=\lambda$ and $\theta\onestr^{c_1-a_1}=\mu$. Then $\eta=\theta^\vee$ if and only if $\mu=\theta\onestr^{c_1-a_1}= \eta^\vee\onestr^{c_1-a_1}$. We write the sum in terms of $\eta$.

    $$\#\{\inlinepar_{\alpha,\gamma,\beta,\delta} \text{-puzzles}\} =\sum_{\lambda,\mu} c_{\lambda\pad,\mu\pad}^{\sort(\alpha)\sort(\delta)}\cdot  c_{\alpha,\beta}^{\lambda^\vee}\cdot c_{\gamma,\delta}^{\mu^\vee} = \sum_{\eta}  c_{\alpha,\beta}^{\eta^\vee\zerostr^{a_0-c_0}}\cdot c_{\gamma,\delta}^{\onestr^{c_1-a_1}\eta}$$

    \item Case: $c_0>a_0$ and $a_1>c_1$. We have $\eta=\zerostr^{c_0-a_0}\lambda$ and $\theta=\mu\onestr^{a_1-c_1}$. Then $\eta=\theta^\vee$ if and only if $\zerostr^{c_0-a_0}\lambda = \onestr^{a_1-c_1}\mu^\vee$, but this is impossible. So in this case we always have $c_{\lambda\pad,\mu\pad}^{\sort(\alpha)\sort(\delta)}=0$, so 
    $$\#\{\inlinepar_{\alpha,\gamma,\beta,\delta} \text{-puzzles}\} =\sum_{\lambda,\mu} c_{\lambda\pad,\mu\pad}^{\sort(\alpha)\sort(\delta)} \cdot c_{\alpha,\beta}^{\lambda^\vee}\cdot c_{\gamma,\delta}^{\mu^\vee} =0.$$  
\end{enumerate}
\end{proof}

\subsection{Commutative property of parallelogram-shaped puzzles}

\begin{figure}[h]
    \centering
    \input{boundaries_commute/parallelogramboundarycommute}   \caption{$\protect\inlinepar_{\alpha,\gamma,\beta,\delta}$}
    \label{fig:parallelogramboundarycommutesolo}
\end{figure}

\begin{theorem}[Commutative Property of Parallelogram-Shaped Puzzles]
    \label{thm:Hparallelogramcommute}
    For $\alpha,\beta\in\binom{[a]}{a_1}$ and $\gamma,\delta\in\binom{[c]}{c_1}$, and using puzzle pieces in $\Hfrak$, we have
    $$ \#\{\inlinepar_{\alpha,\gamma,\beta,\delta} \text{-puzzles}\} =  \#\{\inlinepar_{\beta,\gamma,\alpha,\delta} \text{-puzzles}\} =  \#\{\inlinepar_{\alpha,\delta,\beta,\gamma} \text{-puzzles}\} =  \#\{\inlinepar_{\beta,\delta,\alpha,\gamma} \text{-puzzles}\}.$$
    In other words, we can commute the pairs $\alpha,\beta$ and $\gamma,\delta$ independently and preserve the number of puzzles.
\end{theorem}

\begin{proof}
We give several different arguments.
\begin{enumerate}
    \item \textbf{Geometric proof.} We complete to a triangle as in Proposition \ref{prop:parallelogramcomplete} and have 
    $$\#\{\inlinepar_{\alpha,\gamma,\beta,\delta} \text{-puzzles}\} =  \#\{\inlinetri_{\sort(\alpha)\gamma,\beta\sort(\delta),\delta\alpha} \text{-puzzles}\} = \underset{\Gr_{a+c}}{\int} [X_{\sort(\alpha)\gamma}][X_{\beta\sort(\delta)}][X^{(\delta\alpha)^\vee}].$$
    We use Lemma \ref{lemma:commutesplit} to write 
    $$[X_{\sort(\alpha)\gamma}][X_{\beta\sort(\delta)}][X^{(\delta\alpha)^\vee}] = [X_{\sort(\alpha)\gamma}][X_{\alpha\sort(\delta)}][X^{(\delta\beta)^\vee}],$$
    $$[X_{\sort(\alpha)\gamma}][X_{\beta\sort(\delta)}][X^{(\delta\alpha)^\vee}] = [X_{\sort(\alpha)\delta}][X_{\beta\sort(\delta)}][X^{(\gamma\alpha)^\vee}],$$
    and
    $$[X_{\sort(\alpha)\gamma}][X_{\beta\sort(\delta)}][X^{(\delta\alpha)^\vee}] = [X_{\sort(\alpha)\gamma}][X_{\alpha\sort(\delta)}][X^{(\delta\beta)^\vee}] = [X_{\sort(\alpha)\delta}][X_{\alpha\sort(\delta)}][X^{(\gamma\beta)^\vee}].$$
    Then we have
    \begin{multline*}
        \underset{\Gr_{a+c}}{\int} [X_{\sort(\alpha)\gamma}][X_{\beta\sort(\delta)}][X^{(\delta\alpha)^\vee}] = \underset{\Gr_{a+c}}{\int} [X_{\sort(\alpha)\gamma}][X_{\alpha\sort(\delta)}][X^{(\delta\beta)^\vee}] \\
        = \underset{\Gr_{a+c}}{\int} [X_{\sort(\alpha)\delta}][X_{\beta\sort(\delta)}][X^{(\gamma\alpha)^\vee}] = \underset{\Gr_{a+c}}{\int} [X_{\sort(\alpha)\delta}][X_{\alpha\sort(\delta)}][X^{(\gamma\beta)^\vee}].
    \end{multline*}
    The theorem statement immediately follows.

    \item \textbf{Geometric proof.} It also follows from Theorem \ref{thm:parallelogramLRnumbers}, since the pairs $\alpha,\beta$ and $\gamma,\delta$ commute in each Littlewood-Richardson number. 
    \item \textbf{Proof by commutative property of trapezoidal puzzles.} Commutativity of $\gamma$ and $\delta$ follows directly from Theorem \ref{thm:trapezoidcommute}, since a parallelogram can be seen as a special case of a trapezoid after completing to a triangle as in Figure \ref{fig:parallelogramcomplete}. Then to commute $\alpha$ and $\beta$, we just have to conjugate the trapezoid commutativity with a $120^\circ$ rotation.
    \item \textbf{Proof by commutative property of triangular puzzles.} It can also be proved by using the basic commutative property of triangular puzzles. We complete the parallelogram to a triangle as in Figure \ref{fig:parallelogramcomplete}, then we apply the commutative property to the NW and NE sides of the triangular boundary to get
    $$\#\{\inlinetri_{\sort(\alpha)\gamma,\beta\sort(\delta),\delta\alpha} \text{-puzzles}\}= \#\{\inlinetri_{\beta\sort(\delta),\sort(\alpha)\gamma,\delta\alpha} \text{-puzzles}\}.$$
    
    Now assume we have a $\inlinetri_{\beta\sort(\delta),\sort(\alpha)\gamma,\delta\alpha}$-puzzle. Then because $\alpha$ and $\beta$ have matching content and $\gamma$ and $\delta$ have matching content, there can be no $\tenstr$s appearing along the red lines shown in Figure \ref{fig:parallelogramcommuted}, and if we call the strings appearing along those lines $\lambda$ and $\mu$, then $\lambda$ must have the same content as $\alpha$ and $\beta$ and $\mu$ must have the same content as $\gamma$ and $\delta$. (This comes from statement (a) under the discrete Green's theorem property given in Section \ref{section:otherpuzzleproperties}.) So there must be an individual $\inlinetri_{\beta,\lambda,\alpha}$-puzzle and an individual $\inlinetri_{\mu,\gamma,\delta}$-puzzle in the left and right green triangular regions, respectively.
    
    \begin{figure*}[h!]
       \centering 
        \input{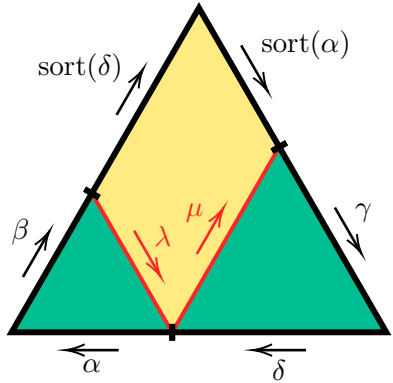}
        \caption{In any $\protect\inlinetri_{\beta\sort(\delta),\sort(\alpha)\gamma,\delta\alpha}$-puzzle, there must be a triangular subpuzzle filling each of the green regions, as well as a paralleogram-shaped subpuzzle filling the yellow region.}
       \label{fig:parallelogramcommuted}
    \end{figure*}

    
    For each fixed choice of $\lambda$ and $\mu$, we can count all $\inlinetri_{\beta\sort(\delta),\sort(\alpha)\gamma,\delta\alpha}$-puzzles that have those labels along those internal lines. To do this, we take the product of the number of puzzles in the two green triangular regions and the number of puzzles in the yellow parallelogram-shaped region, namely $\left(\#\{\inlinepar_{\sort(\delta),\sort(\alpha),\mu^\vee,\lambda^\vee} \text{-puzzles}\}\right) \left(\#\{\inlinetri_{\beta,\lambda,\alpha} \text{-puzzles}\}\right) \left(\#\{\inlinetri_{\mu,\gamma,\delta} \text{-puzzles}\}\right)$. Then to find the number of all $\inlinetri_{\beta\sort(\delta),\sort(\alpha)\gamma,\delta\alpha}$-puzzles, we sum over each choice of fixed pair $\lambda,\mu$. Explicitly, we have the formula
    \begin{multline*}
    \#\{\inlinetri_{\beta\sort(\delta),\sort(\alpha)\gamma,\delta\alpha} \text{-puzzles}\} \\ 
    = \sum_{\lambda,\mu} \left(\#\{\inlinepar_{\sort(\delta),\,\sort(\alpha),\,\mu^\vee,\,\lambda^\vee} \text{-puzzles}\}\right) \left(\#\{\inlinetri_{\beta,\lambda,\alpha} \text{-puzzles}\}\right) \left(\#\{\inlinetri_{\mu,\gamma,\delta} \text{-puzzles}\}\right)\\
    =\sum_{\lambda,\mu} \left(\#\{\inlinepar_{\sort(\delta),\,\sort(\alpha),\,\mu^\vee,\,\lambda^\vee} \text{-puzzles}\}\right) \cdot c_{\beta,\lambda}^{\alpha^\vee} \cdot c_{\mu,\gamma}^{\delta^\vee}.
    \end{multline*}
    Now the fact that $c_{\beta,\lambda}^{\alpha^\vee}= c_{\alpha,\lambda}^{\beta^\vee}$ and $c_{\mu,\gamma}^{\delta^\vee}= c_{\mu,\delta}^{\gamma^\vee}$ implies that the pairs $\alpha,\beta$ and $\gamma,\delta$ can be commuted as in the theorem statement.
\end{enumerate}
\end{proof}

\begin{remark}
    In fact, if the pair $\alpha,\beta$ (or likewise $\gamma,\delta$) do not have matching content, then there exists no filling of $\inlinepar_{\alpha,\gamma,\delta,\beta}$. This is a consequence of the discrete Green's theorem given in Section \ref{section:otherpuzzleproperties}. So in Theorem \ref{thm:Hparallelogramcommute} we did not have to specify matching content for the statement to hold true, vacuously in the case of mismatched content, but that is not a case we wished to emphasize.
\end{remark}

\begin{remark}
    Theorem \ref{thm:Hparallelogramcommute} also implies the commutative property of trapezoidal puzzles, i.e. Theorem \ref{thm:trapezoidcommute}. This is because, due to the unique identity triangle property given in Section \ref{section:otherpuzzleproperties}, there is a bijection
    $$\{\inlinetrapside_{\beta,\gamma,\nu,\delta} \text{-puzzles}\} \leftrightarrow \{\inlinepar_{\sort(\delta)\beta,\gamma,\nu,\delta} \text{-puzzles}\}.$$
    Then we use this and Theorem \ref{thm:Hparallelogramcommute} to write
    $$\#\{\inlinetrapside_{\beta,\gamma,\nu,\delta} \text{-puzzles}\} = \#\{\inlinepar_{\sort(\delta)\beta,\gamma,\nu,\delta} \text{-puzzles}\} = \#\{\inlinepar_{\sort(\delta)\beta,\delta,\nu,\gamma} \text{-puzzles}\}= \#\{\inlinetrapside_{\beta,\delta,\nu,\gamma} \text{-puzzles}\}.$$
\end{remark}

\subsubsection{K-theory version}

\begin{remark}
    \label{remark:Ktheoryparallelogrampuzzles}
    Theorem \ref{thm:Hparallelogramcommute} has an analogous statement for K-theory puzzles using pieces in $\Hfrak\cup\left\{\raisebox{-.4\height}{
\begin{tikzpicture}[x=1pt,y=1pt,yscale=1,xscale=1,scale=.68]
\draw [fill={rgb, 255:red, 255; green, 255; blue, 255 }] (0.0,0.0) -- (15.0,25.98) -- (30.0,0.0) -- cycle;
\draw (7.5,12.99) node  [rotate=60]  {$\mathtt{1\! 0}$};
\draw (22.5,12.99) node  [rotate=-60]  {$\mathtt{1\! 0}$};
\draw (15.0,0.0) node  [rotate=0]  {$\mathtt{1\! 0}$};
\end{tikzpicture}
}\right\}$. Arguments 3 and 4 in the proof of Theorem \ref{thm:Hparallelogramcommute} work the same for these types of puzzles. Argument 1 works similarly but with a little added annoyance of needing to apply rotational symmetry and commutativity of triangular puzzles since the K-theory Schubert bases are not self-dual. 
\end{remark}

\subsection{Puzzle-based proof of Theorem \ref{thm:parallelogramLRnumbers}}
\label{section:uniquefillingproof}

Here we give an alternate proof of Theorem \ref{thm:parallelogramLRnumbers} which uses basic properties of puzzles and combinatorics rather than any geometric arguments. First we prove a lemma and then give the main proof of the theorem.

\subsubsection{Lemma on the unique filling of a partially labeled pentagonal boundary}

\begin{lemma}
    \label{lemma:uniquefilling}


    For integers $a_0,a_1,c_0,c_1$, let $\Pcal(a_0,a_1,c_0,c_1)$ denote a pentagonal boundary with side lengths shown in Figure \ref{fig:pentagonalboundary}. Traveling clockwise, label its NW side with $\zerostr^{c_0}\onestr^{c_1}$ and its NE side with $\zerostr^{a_0}\onestr^{a_1}$. Then, 
    \begin{itemize}
        \item if $a_0\geq c_0$ or $c_1\geq a_1$, there exists a unique filling of this partially labeled pentagonal boundary using puzzle pieces in $\Hfrak\cup\left\{\raisebox{-.4\height}{
\begin{tikzpicture}[x=1pt,y=1pt,yscale=1,xscale=1,scale=.68]
\draw [fill={rgb, 255:red, 255; green, 255; blue, 255 }] (0.0,0.0) -- (15.0,25.98) -- (30.0,0.0) -- cycle;
\draw (7.5,12.99) node  [rotate=60]  {$\mathtt{1\! 0}$};
\draw (22.5,12.99) node  [rotate=-60]  {$\mathtt{1\! 0}$};
\draw (15.0,0.0) node  [rotate=0]  {$\mathtt{1\! 0}$};
\end{tikzpicture}
}\right\}$. Moreover, this filling uses only pieces in $\Hfrak$ and produces SE, South, and SW boundary labels of $\onestr^{\max\{0,c_1-a_1\}}\zerostr^{\max\{0,c_0-a_0\}}$, $\onestr^{\min\{a_1,c_1\}}\zerostr^{\min\{a_0,c_0\}}$, and $\onestr^{\max\{0,a_1-c_1\}}\zerostr^{\max\{0,a_0-c_0\}}$
        , respectively.
        This is illustrated in Figure \ref{fig:pentagonuniquefilling}.
        
        \item if $c_0>a_0$ and $a_1>c_1$, there exists no filling of this partially labeled boundary using puzzle pieces in $\Hfrak\cup\left\{\raisebox{-.4\height}{
\begin{tikzpicture}[x=1pt,y=1pt,yscale=1,xscale=1,scale=.68]
\draw [fill={rgb, 255:red, 255; green, 255; blue, 255 }] (0.0,0.0) -- (15.0,25.98) -- (30.0,0.0) -- cycle;
\draw (7.5,12.99) node  [rotate=60]  {$\mathtt{1\! 0}$};
\draw (22.5,12.99) node  [rotate=-60]  {$\mathtt{1\! 0}$};
\draw (15.0,0.0) node  [rotate=0]  {$\mathtt{1\! 0}$};
\end{tikzpicture}
}\right\}$.
    \end{itemize}
\end{lemma}

\begin{figure*}[h]
    \centering
    \begin{subfigure}[t]{0.475\textwidth}
        \centering
        \input{unique_filling_figures/pentagonalboundary}
        \caption{The pentagonal boundary $\Pcal(a_0,a_1,c_0,c_1)$, with side lengths shown. }
        \label{fig:pentagonalboundary}
    \end{subfigure}%
    ~\, 
    \begin{subfigure}[t]{0.475\textwidth}
        \centering
        \input{unique_filling_figures/pentagonuniquefilling}
        \caption{The unique filling of $\Pcal(a_0,a_1,c_0,c_1)$ with fixed labels on the NW and NE sides, for $a_0\geq c_0$ or $c_1\geq a_1$. The implied boundary labels are shown in lime.}
        \label{fig:pentagonuniquefilling}
    \end{subfigure}
    \caption{}
\end{figure*}

\begin{proof}

The proof will proceed by induction on $a_0,a_1,c_0,c_1$. We will call any case where $a_0\geq c_0$ or $c_1\geq a_1$ a ``good case'' and the case where $c_0>a_0$ and $a_1>c_1$ the ``bad case.''

Also, we use the duality symmetry property of puzzles (see Section \ref{section:basicsymmetriesofpuzzles}), to avoid basically proving everything twice needlessly. Alongside each case we will list the dual case, which can be proved simply by conjugating the provided argument with the operation of reflecting puzzles across the $y$-axis and exchanging $\zerostr$s and $\onestr$s. 
For example, if we take the entire setup in the case where $c_1\geq a_1$ and $c_0>a_0$ and follow the duality symmetry bijection, we end up in the case where $a_0\geq c_0$ and $a_1\geq c_1$. We can prove the theorem in that case and then follow the bijection backwards to verify that the theorem must hold for $c_1\geq a_1$ and $c_0>a_0$ as well.



\begin{figure*}[]
    \centering
    \begin{subfigure}[t]{0.34\textwidth}
        \centering
        \input{unique_filling_figures/operationAtwo_no_opacity}
        \caption{Case: $a_0\geq c_0$ and $a_1>c_1$}
        \label{fig:operationAtwo}
    \end{subfigure}
    ~ 
    \begin{subfigure}[t]{0.38\textwidth}
        \centering
        \input{unique_filling_figures/operationAone_no_opacity}
        \caption{Case: $a_0\geq c_0$, $c_1\geq a_1$, and $a_1>0$}
        \label{fig:operationAone}
    \end{subfigure}%
    ~ 
    \begin{subfigure}[t]{0.24\textwidth}
        \centering
        \input{unique_filling_figures/operationB_no_opacity}
        \caption{Case: $a_1=c_0=0$ and $a_0>0$}
        \label{fig:operationB}
    \end{subfigure}
    \caption{Inductive steps for the proof. The original fixed boundaries and labels are shown in black. A red-highlighted path indicates where fixing the labels implies a unique filling of the orange region, which produces the labels shown in orange. The particular unique fillings of the orange regions are shown in Figure \ref{fig:uniquefillingrows}.}
\end{figure*}

\textbf{Good case:} $a_0\geq c_0$ or $c_1\geq a_1$.

\begin{enumerate}[(a)]
    \item Subcase: $a_0\geq c_0$ and $a_1>c_1$ (dual $c_1\geq a_1$ and $c_0>a_0$).

     We remark that in this case the pentagonal boundary $\Pcal(a_0,a_1,c_0,c_1)$ is degenerate so as to actually be a trapezoid. We draw a smaller pentagonal boundary $\Pcal(a_0,a_1-1,c_0,c_1)$ (also a trapezoid in this case) whose NW and NE sides have fixed labels $\zerostr^{c_0}\onestr^{c_1}$ and $\zerostr^{a_0}\onestr^{a_1-1}$ inside our original boundary $\Pcal(a_0,a_1,c_0,c_1)$. This region has a unique filling with the properties stated in the theorem, according to the inductive hypothesis. Then there is a unique way to extend this to a filling of our original boundary, and this filling satisfies the claim. This is illustrated in Figures \ref{fig:operationAtwo} and \ref{fig:operationAtworow}.

     \item Subcase: $a_0\geq c_0$, $c_1\geq a_1$, and $a_1>0$ (dual $a_0\geq c_0$, $c_1\geq a_1$, and $c_0>0$).

     We draw a smaller pentagonal boundary $\Pcal(a_0,a_1-1,c_0,c_1)$ whose NW and NE sides have fixed labels $\zerostr^{c_0}\onestr^{c_1}$ and $\zerostr^{a_0}\onestr^{a_1-1}$ inside our original boundary. This region has a unique filling with the properties stated in the theorem, according to the inductive hypothesis. Then there is a unique way to extend this to a filling of our original boundary, and this filling satisfies the claim. This is illustrated in Figures \ref{fig:operationAone} and \ref{fig:operationAonerow}.

     \item Subcase: $a_1=c_0=0$ and $a_0>0$ (dual $a_1=c_0=0$ and $c_1>0$).

    We remark that in this case the pentagonal boundary $\Pcal(a_0,a_1,c_0,c_1)$ is degenerate so as to actually be a parallelogram. There exists a unique filling of the top NW row inside our original partially labeled boundary, which is given in Figure \ref{fig:operationBrow}, and we fix that in place. Then we draw a smaller pentagonal boundary $\Pcal(a_0-1,0,0,c_1)$ (also a parallelogram in this case) whose NW side is the SW side of that top row. Carrying over the labels we have fixed so far, the smaller pentagonal boundary's NW and NE sides are labeled with $\onestr^{c_1}$ and $\zerostr^{a_0-1}$, respectively. This region has a unique filling with the properties stated in the theorem, according to the inductive hypothesis, and taking the union of this filling and the filling of our fixed top row, we have that our original partially labeled boundary also has a unique filling that satisfies the claim. This is illustrated in Figure \ref{fig:operationB}. The difference here is that in the first two cases, the purple region helps imply the orange region, but in this case, it is the other way around.

    \item Base case: $a_0=a_1=c_0=c_1=0$.

    There is a unique empty filling satisfying the claim in this case.

\end{enumerate}

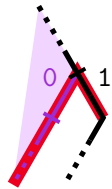
\begin{figure*}[]
    \centering
    \begin{subfigure}[t]{0.38\textwidth}
        \centering
        \input{unique_filling_figures/operationAtworow_no_opacity}
        \caption{Unique filling of the orange region in Figure \ref{fig:operationAtwo}. The string of labels on the North side is replicated on the South side. }
        \label{fig:operationAtworow}
    \end{subfigure}
 ~\,
    \begin{subfigure}[t]{0.26\textwidth}
        \centering
        \input{unique_filling_figures/operationAonerow_no_opacity}
        \caption{Unique filling of the orange region in Figure \ref{fig:operationAone}. The string of labels on the SE side has one fewer $\onestr$ than on the NW side.}
        \label{fig:operationAonerow}
    \end{subfigure}%
~\,
    \begin{subfigure}[t]{0.24\textwidth}
        \centering
        \input{unique_filling_figures/operationBrow_no_opacity}
        \caption{Unique filling of the orange region in Figure \ref{fig:operationB}. The string of labels on the NW side is replicated on the SE side.}
        \label{fig:operationBrow}
    \end{subfigure}
    \caption{Each figure shows the unique filling, using puzzle pieces in $\Hfrak\cup\left\{\protect\raisebox{-.4\height}{\input{puzzle_pieces/cccdeltapiece}}\right\}$, of the given row, assuming we've fixed the labels on the path highlighted in red. Only pieces in $\Hfrak$ end up being used. 
    }
    \label{fig:uniquefillingrows}
\end{figure*}

\textbf{Bad case:} $c_0>a_0$ and $a_1>c_1$.


In this case we prove that there is no way to fill the partially labeled pentagonal boundary.

\begin{enumerate}[(a)]
    \item Subcase: $a_1-c_1>1$ (dual $c_0-a_0>1$).

    We can draw a smaller pentagonal boundary $\Pcal(a_0,a_1-1,c_0,c_1)$ whose NW and NE sides have fixed labels $\zerostr^{c_0}\onestr^{c_1}$ and $\zerostr^{a_0}\onestr^{a_1-1}$ inside of our original boundary. However, this itself is still in the bad case since $a_1-1>c_1$, so by the inductive hypothesis, there exists no filling. Then since there is no possible filling of this smaller partially labeled boundary, there can be no filling of the larger one which contains it. 

    \item Subcase: $a_1-c_1=1$ (dual $c_0-a_0=1$).
    
    We can draw a smaller partially labeled pentagonal boundary $\Pcal(a_0,a_1-1,c_0,c_1)$ whose NW and NE sides have fixed labels $\zerostr^{c_0}\onestr^{c_1}$ and $\zerostr^{a_0}\onestr^{a_1-1}$ inside of our original boundary. Since $c_1=a_1-1$ and $c_0>a_0$, this smaller partially labeled boundary falls into the dual version of subcase (a) of the ``good case,'' so it has a unique filling as described in the theorem. Considering our parameters and referring to Figure \ref{fig:pentagonuniquefilling}, we can see there will be no $\onestr$s on the SE side of this filling, only $\zerostr$s. We end up with a situation where the rightmost portion of the original partially labeled boundary, with the unique filling of the smaller boundary shown in purple, looks like:
    \begin{center}\input{unique_filling_figures/badcase_no_opacity}\end{center} 
    
    We would need to find a way to fill the remaining empty region here, but given the fixed labels on the path highlighted in red, this is impossible because there is no allowed puzzle piece with a $\zerostr$ on its NW side and a $\onestr$ on its NE side. Therefore, there is no way to extend the unique filling of the smaller region to a filling of our original partially labeled boundary.

    \item Base case: $c_0=a_0=1$ and $a_1=c_1=0$.
    
    The boundary $\Pcal(0,1,1,0)$ with the partial labeling described in the theorem is shown below:

    \begin{center}\input{unique_filling_figures/badcasebasecase}\end{center} 
    
    There is no allowed puzzle piece with  $\zerostr$ on its NW side and a $\onestr$ on its NE side, so it is impossible to fill this partially labeled boundary.  
\end{enumerate}
\end{proof}

\begin{remark}
    Lemma \ref{lemma:uniquefilling} does not hold if we are using the piece \raisebox{-.4\height}{
\begin{tikzpicture}[x=1pt,y=1pt,yscale=1,xscale=1,scale=.68]
\draw [fill={rgb, 255:red, 255; green, 255; blue, 255 }] (30.0,25.98) -- (15.0,0.0) -- (0.0,25.98) -- cycle;
\draw (22.5,12.99) node  [rotate=60]  {$\mathtt{1\! 0}$};
\draw (7.5,12.99) node  [rotate=-60]  {$\mathtt{1\! 0}$};
\draw (15.0,25.98) node  [rotate=0]  {$\mathtt{1\! 0}$};
\end{tikzpicture}
}. In Figure \ref{fig:operationBrow}, it can be seen that if we allow this piece, there is another way to fill the row, where we replace the very bottom piece with a \raisebox{-.4\height}{
\begin{tikzpicture}[x=1pt,y=1pt,yscale=1,xscale=1,scale=.68]
\draw [fill={rgb, 255:red, 255; green, 255; blue, 255 }] (30.0,25.98) -- (15.0,0.0) -- (0.0,25.98) -- cycle;
\draw (22.5,12.99) node  [rotate=60]  {$\mathtt{1\! 0}$};
\draw (7.5,12.99) node  [rotate=-60]  {$\mathtt{1\! 0}$};
\draw (15.0,25.98) node  [rotate=0]  {$\mathtt{1\! 0}$};
\end{tikzpicture}
}. It also does not hold if we use the equivariant piece \raisebox{-.4\height}{
\begin{tikzpicture}[x=1pt,y=1pt,yscale=1,xscale=1,scale=.68]
\draw [fill={rgb, 255:red, 255; green, 246; blue, 74 }] (0.0,0.0) -- (15.0,25.98) -- (30.0,0.0) -- (15.0,-25.98) -- cycle;
\draw (7.5,12.99) node  [rotate=60]  {$\mathtt{0}$};
\draw (22.5,12.99) node  [rotate=-60]  {$\mathtt{1}$};
\draw (22.5,-12.99) node  [rotate=60]  {$\mathtt{0}$};
\draw (7.5,-12.99) node  [rotate=-60]  {$\mathtt{1}$};
\end{tikzpicture}
}, since then there \textit{is} a way to fill the partially labeled boundary in the base case of the ``bad case'' above.
\end{remark}

\subsubsection{Proof of the theorem}

Now we return to the main part of the proof. We will actually begin by picking up where we left off at the end of Argument 4 in the proof of Theorem \ref{thm:Hparallelogramcommute}, which this proof will take further. 

Let us now analyze the set $\{\inlinepar_{\sort(\delta),\,\sort(\alpha),\,\mu^\vee,\,\lambda^\vee} \text{-puzzles}\}$ of puzzles that can fill the yellow parallelogram-shaped region in Figure \ref{fig:parallelogramcommuted}. This region can be divided into the pentagonal region bounded by $\Pcal(a_0,a_1,c_0,c_1)$ as defined in Lemma \ref{lemma:uniquefilling} and an inverted triangular region with side length $\min\{a_0,c_0\}+\min\{a_1,c_1\}$, as shown in Figure \ref{fig:yellowparallelogramboundary}.

We know from Lemma \ref{lemma:uniquefilling} that there is a unique filling of the pentagonal region except in the ``bad case'' where $c_0>a_0$ and $a_1>c_1$. So in the bad case, there can be no filling of the yellow parallelogram-shaped region, and thus no $\inlinetri_{\beta\sort(\delta),\sort(\alpha)\gamma,\delta\alpha}$-puzzles. 

Now let us assume we are in a ``good case'' where $a_0\geq c_0$ or $c_1\geq a_0$, so we have the unique filling of the pentagonal region described in Lemma \ref{lemma:uniquefilling}. Fixing this in place, now we consider what fillings may appear in the inverted triangular region on the bottom.
The label $\zerostr^{\min\{a_0,c_0\}}\onestr^{\min\{a_1,c_1\}}$ on the North side of the inverted triangular region (inherited from the South side of the unique filling of the pentagonal region) implies that for any choice of label $\rho$ on its SE side, with same content, there is a unique filling of the inverted triangular region, and this implies the label $\rho^\vee$ on the SW side. (This is due to the unique identity puzzle property described in Section \ref{section:otherpuzzleproperties}.) This is all illustrated in Figure \ref{fig:yellowparallelogramfilling}. (Note that if $\rho$ does not have the same content as $\zerostr^{\min\{a_0,c_0\}}\onestr^{\min\{a_1,c_1\}}$, there would be $\tenstr$s on the SW side for any filling of the inverted triangular region, which will not be allowed given that the labels $\lambda$ and $\mu$ in Figure \ref{fig:parallelogramcommuted} must not have $\tenstr$s. So we have excluded that possibility.) Taking the union of the unique filling of the pentagonal region and this filling of the inverted triangular region gives us a filling of the yellow parallelogram-shaped region. Varying over $\rho$ gives us all possible fillings of this region.


\begin{figure}[h]
    \centering
    \begin{subfigure}[t]{0.47\textwidth}
    \centering
     \input{unique_filling_figures/yellowparallelogramboundary}
    \caption{The yellow parallelogram-shaped region in Figure \ref{fig:parallelogramcommuted} can be divided into the pentagonal region of Lemma \ref{lemma:uniquefilling} and an inverted triangular region, with side lengths shown.}
    \label{fig:yellowparallelogramboundary}   
    \end{subfigure}
    ~\;
    \begin{subfigure}[t]{0.47\textwidth}
    \centering
    \input{unique_filling_figures/yellowparallelogramfilling_no_opacity}
    \caption{Given the unique filling of the pentagonal region (as in Lemma \ref{lemma:uniquefilling}) and a fixed choice of $\rho$, the labels on the path highlighted in red imply there is a unique filling of the orange inverted triangular region, which produces the label $\rho^\vee$ as shown.}
    \label{fig:yellowparallelogramfilling}
    \end{subfigure}
    \caption{}
\end{figure}

We have just shown that, assuming $a_0\geq c_0$ or $c_1\geq a_1$, there exists a $\inlinepar_{\sort(\delta),\,\sort(\alpha),\,\mu^\vee,\,\lambda^\vee}$-puzzle if and only if 
$$\mu^\vee=\onestr^{\max\{0,c_1-a_1\}}\zerostr^{\max\{0,c_0-a_0\}}\rho$$
and 
$$\lambda^\vee=\rho^\vee\onestr^{\max\{0,a_1-c_1\}}\zerostr^{\max\{0,a_0-c_0\}}$$ 
for some $\rho$ with the content of $\zerostr^{\min\{a_0,c_0\}}\onestr^{\min\{a_1,c_1\}}$, and in this case the puzzle is unique.

Then returning to the formula we had at the end of argument 3 in the proof of Theorem \ref{thm:Hparallelogramcommute}, we have
\begin{multline*}
    \#\{\inlinetri_{\beta\sort(\delta),\sort(\alpha)\gamma,\delta\alpha} \text{-puzzles}\} 
    =\sum_{\lambda,\mu} \left(\#\{\inlinepar_{\sort(\delta),\,\sort(\alpha),\,\mu^\vee,\,\lambda^\vee} \text{-puzzles}\}\right) \cdot c_{\beta,\lambda}^{\alpha^\vee} \cdot c_{\mu,\gamma}^{\delta^\vee} \\
    =\sum_{\lambda,\mu} \left(\#\{\inlinepar_{\sort(\delta),\,\sort(\alpha),\,\mu^\vee,\,\lambda^\vee} \text{-puzzles}\}\right) \cdot c_{\alpha,\beta}^{\lambda^\vee} \cdot c_{\gamma,\delta}^{\mu^\vee} \\
    =\begin{cases}
    \sum\limits_{\rho} c_{\alpha,\beta}^{\rho^\vee\onestr^{\max\{0,a_1-c_1\}}\zerostr^{\max\{0,a_0-c_0\}}} \cdot c_{\gamma,\delta}^{\onestr^{\max\{0,c_1-a_1\}}\zerostr^{\max\{0,c_0-a_0\}}\rho} & \text{if } a_0\geq c_0 \text{ or } c_1\geq a_1 \\ 0 & \text{ if } c_0>a_0 \text{ and } a_1>c_1
\end{cases}.
\end{multline*}

Now if we split the case where $a_0\geq c_0$ or $c_1\geq a_1$ into subcases, we immediately get the statement of the theorem we are proving, Theorem \ref{thm:parallelogramLRnumbers}.

\subsubsection{K-theory version}
\begin{remark}
    We can use the same general proof method to get a formula in the case of K-theory puzzles using pieces in $\Hfrak\cup\left\{\raisebox{-.4\height}{
\begin{tikzpicture}[x=1pt,y=1pt,yscale=1,xscale=1,scale=.68]
\draw [fill={rgb, 255:red, 255; green, 255; blue, 255 }] (0.0,0.0) -- (15.0,25.98) -- (30.0,0.0) -- cycle;
\draw (7.5,12.99) node  [rotate=60]  {$\mathtt{1\! 0}$};
\draw (22.5,12.99) node  [rotate=-60]  {$\mathtt{1\! 0}$};
\draw (15.0,0.0) node  [rotate=0]  {$\mathtt{1\! 0}$};
\end{tikzpicture}
}\right\}$. Lemma \ref{lemma:uniquefilling} allows for use of these pieces, not just pieces in $\Hfrak$. The main proof can proceed identically until we get to the discussion of fillings of the inverted triangular region. If we rotate the situation $180^\circ$, the inverted triangular region allowing the \raisebox{-.4\height}{
\begin{tikzpicture}[x=1pt,y=1pt,yscale=1,xscale=1,scale=.68]
\draw [fill={rgb, 255:red, 255; green, 255; blue, 255 }] (0.0,0.0) -- (15.0,25.98) -- (30.0,0.0) -- cycle;
\draw (7.5,12.99) node  [rotate=60]  {$\mathtt{1\! 0}$};
\draw (22.5,12.99) node  [rotate=-60]  {$\mathtt{1\! 0}$};
\draw (15.0,0.0) node  [rotate=0]  {$\mathtt{1\! 0}$};
\end{tikzpicture}
} piece becomes an upright triangular region allowing the \raisebox{-.4\height}{
\begin{tikzpicture}[x=1pt,y=1pt,yscale=1,xscale=1,scale=.68]
\draw [fill={rgb, 255:red, 255; green, 255; blue, 255 }] (30.0,25.98) -- (15.0,0.0) -- (0.0,25.98) -- cycle;
\draw (22.5,12.99) node  [rotate=60]  {$\mathtt{1\! 0}$};
\draw (7.5,12.99) node  [rotate=-60]  {$\mathtt{1\! 0}$};
\draw (15.0,25.98) node  [rotate=0]  {$\mathtt{1\! 0}$};
\end{tikzpicture}
} piece. The unique identity puzzle property does not hold when the latter piece is involved. In this case, for a given choice of $\rho$, there is not a unique filling of the triangular region, and the other label need not be $\rho^\vee$. So in this case, we just need to account for and sum over each possible \textit{pair} of labels $\rho,\sigma$ on the SE and SW sides of the inverted triangular region in our formula.

    \explainalot{ 
    The contribution of weights of the puzzles in that inverted triangular region will then be
    $$\sum_{\rho,\sigma} \underset{\Gr({\min\{a_1,c_1\},\C^{\min\{a_0,c_0\}+\min\{a_1,c_1\}}})}{\Kint} [\Ical^{\rho^\vee}][\Ical^{\sigma^\vee}][\Ocal_{\zerostr^{\min\{a_0,c_0\}}\onestr^{\min\{a_1,c_1\}}}].$$
    And then $\lambda$ and $\mu$ can vary independently of each other and will contain whatever $\rho$ and $\sigma$ are. }
\end{remark}

\subsection{Commutative property of parallelogram-shaped equivariant puzzles}

\begin{figure}[h]
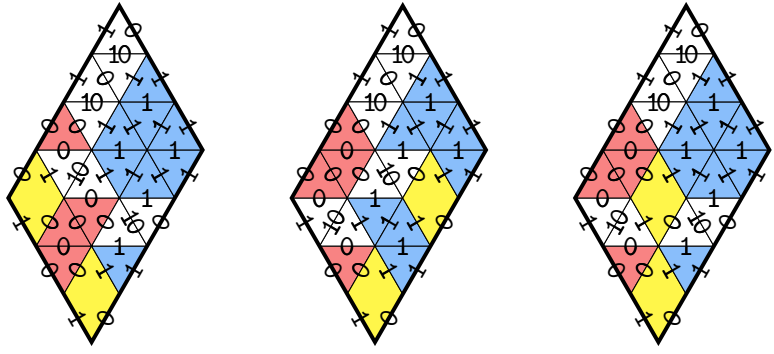
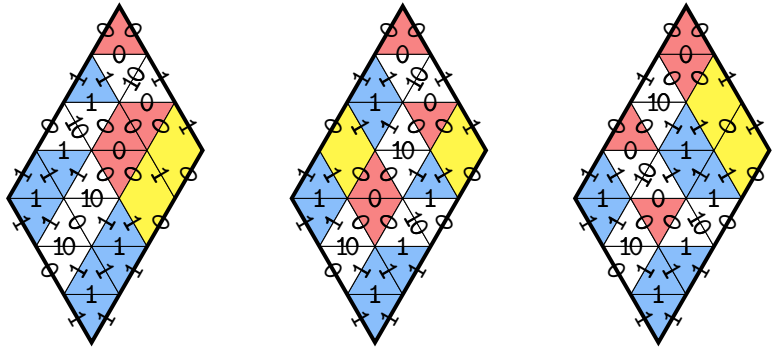

\centering
    \begin{subfigure}[b]{.8\textwidth}
        \centering
        \input{equivariant_parallelogram_example/eqvtpar_v1_row}
        \caption{The set of all $\protect\inlinepar_{\onestr\zerostr\onestr,\zerostr\zerostr\onestr\onestr,\zerostr\onestr\onestr,\onestr\zerostr\onestr\zerostr}$-puzzles using pieces in $\Hfrak \cup \left\{\raisebox{-.4\height}{\input{puzzle_pieces/eqvtpiece}} \right\}$. 
        From left to right, their weights are $(y_4-y_1)(y_4-y_3)$, $(y_4-y_3)(y_6-y_3)$, and $(y_4-y_3)(y_5-y_2)$.}
        \label{fig:eqvtpar_v1}
    \end{subfigure}
    
    \vspace{1\baselineskip}
    \begin{subfigure}[b]{.8\textwidth}
        \centering
        \input{equivariant_parallelogram_example/eqvtpar_v2_row}
        \caption{The set of all $\protect\inlinepar_{\onestr\zerostr\onestr,\onestr\zerostr\onestr\zerostr,\zerostr\onestr\onestr,\onestr\onestr\zerostr\zerostr}$-puzzles using pieces in $\Hfrak \cup \left\{\raisebox{-.4\height}{\input{puzzle_pieces/eqvtpiece}} \right\}$. From left to right, their weights are $(y_7-y_3)(y_6-y_3)$, $(y_5-y_1)(y_7-y_3)$, and $(y_7-y_3)(y_7-y_2)$.}
        \label{fig:eqvtpar_v2}
    \end{subfigure}
    
    \caption{An example of commuting the NW and SE labels for parallelogram-shaped equivariant puzzles}
    \label{fig:equivariantparallelogramsexample}
\end{figure}

The equivariant structure constants in the theorem below are the ones associated to parallelogram-shaped equivariant puzzles, after completing to a triangle as in Figure \ref{fig:parallelogramcomplete}. Since the grey triangular regions each have a unique filling with no equivariant pieces, the structure constant depends only on the fillings of the parallelogram-shaped region and the positions of the equivariant pieces within them. So although the theorem is stated and proved entirely in a geometric manner, it is also a combinatorial statement. 

Now our structure constants are polynomials, and instead of remaining unchanged as we commute the opposite boundary labels on the parallelogram, the indeterminates undergo a permutation dependent on which labels were swapped. This is a stronger statement than the commutative property for parallelogram-shaped puzzles using only pieces in $\Hfrak$, because the equivariant result implies that one, and here we even still find that the \textit{number} of parallelogram-shaped equivariant puzzles is preserved.


\begin{theorem}[Commutative Property of Structure Constants Associated to Parallelogram-Shaped Equivariant Puzzles]
    \label{thm:eqvtparallelogram}

    Define block matrices
$$\Phi_a:=\left[
\begin{array}{c|c}
J_{a} & \mathbf{0} \\ \hline
\mathbf{0} & I_{c}
\end{array}\right]
\quad \text{and} \quad
\Phi_c:=\left[
\begin{array}{c|c}
I_{a} & \mathbf{0} \\ \hline
\mathbf{0} & J_{c}
\end{array}\right]
$$
in $\GL(\C^{a+c})$, where $I_a$ and $J_a$ (resp. $I_c$ and $J_c$) denote the $a\times a$ (resp. $c\times c$) identity and anti-diagonal permutation matrices, respectively.

Let $\alpha,\beta\in\binom{[a]}{a_1}$ and $\gamma,\delta\in\binom{[c]}{c_1}$. Then in $H^*_T(\Gr(a_1+c_1;\C^{a+c}))$, we have 

$$(c_T)_{\sort(\alpha)\gamma,\beta\sort(\delta)}^{(\delta\alpha)^\vee}= \Phi_a\cdot (c_T)_{\sort(\beta)\gamma,\alpha\sort(\delta)}^{(\delta\beta)^\vee} = \Phi_c \cdot (c_T)_{\sort(\alpha)\delta,\beta\sort(\gamma)}^{(\gamma\alpha)^\vee} = \Phi_c\cdot\Phi_a\cdot (c_T)_{\sort(\beta)\delta,\alpha\sort(\gamma)}^{(\gamma\beta)^\vee}.$$
      
In other words commuting the pair $\alpha, \beta$ reverses the $y_1,\ldots,y_a$, and commuting the pair $\gamma,\delta$ reverses the $y_{a+1},\ldots,y_{a+c}$, in the structure constant.
\end{theorem}

\begin{proof}

First we make the general comment that the Weyl group $W:=N(T)/T \cong S_n$ acts on $H^*_T(\Gr(k;\C^n))$ by permutation of the indeterminates (i.e. $w\cdot y_i=y_{w(i)}$) along with $w\cdot[X_\lambda]=[w\cdot X_\lambda]$ for $w\in W$, and this action also commutes with integration. This will be relevant later in the proof as we act with $\Phi_a,\Phi_c\in W$.

Now we begin by showing that $\Phi_a\cdot X_{\sort(\alpha)\gamma}= X_{\sort(\alpha)\gamma}$. To do this we consider the conditions encoded by the string $\sort(\alpha)\gamma= \zerostr^{a_0}\onestr^{a_1}\gamma$. There is no $\onestr$ appearing directly in front of a $\zerostr$ anywhere in the first $a-1$ places, and therefore the string does not encode an essential condition involving any of the subspaces $F_1,F_2, \ldots, F_{a-1}$ of the standard complete flag. Thus we just need to check that the last $c+1$ conditions are satisfied. Let $V\in X_{\sort(\alpha)\gamma}$. For all $i\geq a$ we clearly have $\Phi_a\cdot F_i = F_i$, so $\Phi_a\cdot (V\cap F_i)\subseteq F_i$. Then we have $\Phi_a\cdot (V\cap F_i) = (\Phi_a\cdot (V\cap F_i))\cap F_i \subseteq (\Phi_a\cdot V)\cap F_i$, so $\dim((\Phi_a\cdot V)\cap F_i)\geq \dim(\Phi_a\cdot (V\cap F_i))=\dim(V\cap F_i)$. This implies $\Phi_a\cdot X_{\sort(\alpha)\gamma}\subseteq X_{\sort(\alpha)\gamma}$. Since $\Phi_a$ is an involution, we also have that $\Phi_a\cdot (\Phi_a\cdot V)=V$, so $ X_{\sort(\alpha)\gamma}\subseteq \Phi_a\cdot X_{\sort(\alpha)\gamma}$.

This implies $\Phi_a\cdot [X_{\sort(\alpha)\gamma}]= X_{\sort(\alpha)\gamma}$. We also use Lemma \ref{lemma:commutesplit} to write $\Phi_a\cdot ([X_{\alpha\sort(\delta)}][X^{(\delta\beta)^\vee}]) = [X_{\beta\sort(\delta)}][X^{(\delta\alpha)^\vee}]$.
Then we have
\begin{multline*}
    \Phi_a\cdot (c_T)_{\sort(\beta)\gamma,\alpha\sort(\delta)}^{(\delta\beta)^\vee} = \Phi_a\cdot (c_T)_{\sort(\alpha)\gamma,\alpha\sort(\delta)}^{(\delta\beta)^\vee}= \Phi_a \cdot \underset{\Gr_{a+c}}{\int} [X_{\sort(\alpha)\gamma}] [X_{\alpha\sort(\delta)}][X^{(\delta\beta)^\vee}]\\
    =  \underset{\Gr_{a+c}}{\int} \Phi_a\cdot\left([X_{\sort(\alpha)\gamma}][X_{\alpha\sort(\delta)}][X^{(\delta\beta)^\vee}]\right) = \underset{\Gr_{a+c}}{\int} \left(\Phi_a\cdot[X_{\sort(\alpha)\gamma}]\right) \left(\Phi_a\cdot([X_{\alpha\sort(\delta)}][X^{(\delta\beta)^\vee}])\right)\\
    = \underset{\Gr_{a+c}}{\int} [X_{\sort(\alpha)\gamma}] [X_{\beta\sort(\delta)}][X^{(\delta\alpha)^\vee}] =(c_T)_{\sort(\alpha)\gamma,\beta\sort(\delta)}^{(\delta\alpha)^\vee}.
\end{multline*}

Now we show that $\Phi_c\cdot X_{\beta\sort(\delta)}=  X_{\beta\sort(\delta)}$. In the string $\beta\sort(\delta)$, there is no $\onestr$ appearing directly in front of a $\zerostr$ anywhere in the last $c$ places, so it does not encode an essential condition involving any of the subspaces $F_a,F_{a+1},\ldots,F_{a+c}$ of the standard complete flag. Let $V\in X_{\beta\sort(\delta)}$. Thus we just need to check that the first $a$ conditions are satisfied. For all $i\leq a$, we have $\Phi_c\cdot (V\cap F_i)= V\cap F_i$ since $\Phi_c$ acts trivially on $F_i$. Then we have $V\cap F_i = (V\cap F_i)\cap F_i= (\Phi_c\cdot (V\cap F_i))\cap F_i \subseteq (\Phi_c\cdot V)\cap F_i$, so $\dim((\Phi_c\cdot V)\cap F_i)\geq \dim(V\cap F_i)$. This implies $\Phi_c\cdot X_{\beta\sort(\delta)}\subseteq X_{\beta\sort(\delta)}$. Also, $\Phi_c\cdot (\Phi_c\cdot V)= V$, so $X_{\beta\sort(\delta)}\subseteq \Phi_c \cdot X_{\beta\sort(\delta)}$. 

This implies $\Phi_c\cdot [X_{\beta\sort(\delta)}]=  [X_{\beta\sort(\delta)}]$. We also use Lemma \ref{lemma:commutesplit} to write $\Phi_c\cdot ([X_{\sort(\alpha)\delta}][X^{(\gamma\alpha)^\vee}]) = [X_{\sort(\alpha)\gamma}][X^{(\delta\alpha)^\vee}]$. Then we have
\begin{multline*}
    \Phi_c\cdot (c_T)_{\sort(\alpha)\delta,\beta\sort(\gamma)}^{(\gamma\alpha)^\vee} = \Phi_c\cdot (c_T)_{\sort(\alpha)\delta,\beta\sort(\delta)}^{(\gamma\alpha)^\vee} = \Phi_c\cdot \underset{\Gr_{a+c}}{\int}  [X_{\beta\sort(\delta)}][X_{\sort(\alpha)\delta}][X^{(\gamma\alpha)^\vee}]\\
    = \underset{\Gr_{a+c}}{\int} \Phi_c\cdot\left( [X_{\beta\sort(\delta)}][X_{\sort(\alpha)\delta}][X^{(\gamma\alpha)^\vee}]\right)= \underset{\Gr_{a+c}}{\int} \left(\Phi_c\cdot [X_{\beta\sort(\delta)}]\right) \left(\Phi_c\cdot [X_{\sort(\alpha)\delta}][X^{(\gamma\alpha)^\vee}]\right)\\
    = \underset{\Gr_{a+c}}{\int} [X_{\beta\sort(\delta)}][X_{\sort(\alpha)\gamma}][X^{(\delta\alpha)^\vee}]=(c_T)_{\sort(\alpha)\gamma,\beta\sort(\delta)}^{(\delta\alpha)^\vee}.
\end{multline*}

Now we have taken care of the first two equalities in the theorem statement. Finally, to get $(c_T)_{\sort(\alpha)\gamma,\beta\sort(\delta)}^{(\delta\alpha)^\vee}= \Phi_c\cdot\Phi_a\cdot (c_T)_{\sort(\beta)\delta,\alpha\sort(\gamma)}^{(\gamma\beta)^\vee}$, we apply the previous two arguments one after the other to act with $\Phi_a$ and $\Phi_c$.
\end{proof}

\begin{remark}
    To be more explicit about the combinatorial interpretation of Theorem \ref{thm:eqvtparallelogram}, we recall Section \ref{section:equivariantpuzzles} about equivariant puzzles and Proposition \ref{prop:parallelogramcomplete}. We can write
\begin{align*}    (c_T)_{\sort(\alpha)\gamma,\beta\sort(\delta)}^{(\delta\alpha)^\vee} &= \sum_{\inlinetri_{\sort(\alpha)\gamma,\beta\sort(\delta),\delta\alpha}\text{-puzzles } P} \wt(P)  = \sum_{\inlinetri_{\sort(\alpha)\gamma,\beta\sort(\delta),\delta\alpha}\text{-puzzles } P} \left( \prod_{\substack{\text{equivariant} \\ \text{pieces }p \text{ in } P }} \wt(p) \right) \\
&= \sum_{\inlinepar_{\alpha,\gamma,\beta,\delta}\text{-puzzles } P} \wt(P) = \sum_{\inlinepar_{\alpha,\gamma,\beta,\delta}\text{-puzzles } P} \left( \prod_{\substack{\text{equivariant} \\ \text{pieces }p \text{ in } P }} \wt(p) \right),
\end{align*}

where the weight of an equivariant piece in a parallelogram-shaped puzzle can just be defined to be the same as its weight in the puzzle's completion to a triangle (this is called the MS-weight in \cite{KTpuzzles}). The point is that the weight identifies its position in the puzzle. Then we can consider the combinatorial implications of acting with the elements $\Phi_a$ and $\Phi_c$ on the above expression. In this context the statement of the theorem translates to
\begin{multline*}
    \sum_{\inlinepar_{\alpha,\gamma,\beta,\delta}\text{-puzzles } P} \left( \prod_{\substack{\text{equivariant} \\ \text{pieces }p \text{ in } P }} \wt(p) \right) = \sum_{\inlinepar_{\beta,\gamma,\alpha,\delta}\text{-puzzles } P} \left( \prod_{\substack{\text{equivariant} \\ \text{pieces }p \text{ in } P }} \Phi_a\cdot \wt(p) \right) \\
    = \sum_{\inlinepar_{\alpha,\delta,\beta,\gamma}\text{-puzzles } P} \left( \prod_{\substack{\text{equivariant} \\ \text{pieces }p \text{ in } P }} \Phi_c\cdot \wt(p) \right) = \sum_{\inlinepar_{\beta,\delta,\alpha,\gamma}\text{-puzzles } P} \left( \prod_{\substack{\text{equivariant} \\ \text{pieces }p \text{ in } P }} \Phi_c\cdot\Phi_a\cdot \wt(p) \right).
\end{multline*}
\end{remark}

\begin{remark}
    Theorem \ref{thm:eqvtparallelogram} also proves the commutative property of parallelogram-shaped puzzles for ordinary cohomology (Theorem \ref{thm:Hparallelogramcommute}), because the actions of $\Phi_a$ and $\Phi_c$ are trivial in degree 0 of $H^*_T$.
\end{remark}

\begin{corollary}
    \label{cor:eqvtparallelogramnumber}
    Let $\alpha,\beta\in\binom{[a]}{a_1}$ and $\gamma,\delta\in\binom{[c]}{c_1}$. For equivariant puzzles, i.e. using puzzle pieces in the set $\Hfrak \cup \left\{\raisebox{-.4\height}{
\begin{tikzpicture}[x=1pt,y=1pt,yscale=1,xscale=1,scale=.68]
\draw [fill={rgb, 255:red, 255; green, 246; blue, 74 }] (0.0,0.0) -- (15.0,25.98) -- (30.0,0.0) -- (15.0,-25.98) -- cycle;
\draw (7.5,12.99) node  [rotate=60]  {$\mathtt{0}$};
\draw (22.5,12.99) node  [rotate=-60]  {$\mathtt{1}$};
\draw (22.5,-12.99) node  [rotate=60]  {$\mathtt{0}$};
\draw (7.5,-12.99) node  [rotate=-60]  {$\mathtt{1}$};
\end{tikzpicture}
} \right\}$, we have
    $$ \#\{\inlinepar_{\alpha,\gamma,\beta,\delta} \text{-puzzles}\} =  \#\{\inlinepar_{\beta,\gamma,\alpha,\delta} \text{-puzzles}\} =  \#\{\inlinepar_{\alpha,\delta,\beta,\gamma} \text{-puzzles}\} =  \#\{\inlinepar_{\beta,\delta,\alpha,\gamma} \text{-puzzles}\}.$$
\end{corollary}
\begin{proof}
    This does not automatically follow from Theorem \ref{thm:eqvtparallelogram}, but requires a further simple proof. 
    In the completion of a parallelogram to a triangle (see Figure \ref{fig:parallelogramcomplete}), since the unique identity puzzles on the lower left and lower right do not have equivariant pieces, any equivariant pieces must appear within the parallelogram-shaped region. 
    This implies that the weight of any such piece would be $y_j-y_i$ where $1\leq i\leq a$ and $a+1\leq j\leq a+c$, due to the positioning of the pieces. Then the weight of a $\inlinetri_{\sort(\alpha)\gamma,\beta\sort(\delta),\delta\alpha}$-puzzle will be a product $(y_{j_1}-y_{i_1})(y_{j_2}-y_{i_2})\cdots (y_{j_r}-y_{i_r})$ of such weights. Consider the ring homomorphism $\psi: \Z[y_1,\ldots,y_{a+c}]\ra \Z$ defined by $y_1,\ldots,y_a\mapsto 0$ and $y_{a+1},\ldots,y_{a+c}\mapsto 1$.
    We have $\psi\left((y_{j_1}-y_{i_1})(y_{j_2}-y_{i_2})\cdots (y_{j_r}-y_{i_r})\right)= (1-0)(1-0)\cdots (1-0)= 1$, so $\psi$ maps the weight of any $\inlinetri_{\sort(\alpha)\gamma,\beta\sort(\delta),\delta\alpha}$-puzzle to $1$. Thus $\psi$ maps the structure constant $(c_T)_{\sort(\alpha)\gamma,\beta\sort(\delta)}^{(\delta\alpha)^\vee}$, which is the sum of all the puzzle weights, to a count of the puzzles with the associated boundary labels, i.e. 
    \begin{multline*}               \psi\left(c_{\sort(\alpha)\gamma,\beta\sort(\delta)}^{(\delta\alpha)^\vee}\right) = \psi\left(\sum_{\inlinetri_{\sort(\alpha)\gamma,\beta\sort(\delta),\delta\alpha} \text{-puzzles } P} \wt(P)\right) =  \sum_{\inlinetri_{\sort(\alpha)\gamma,\beta\sort(\delta),\delta\alpha} \text{-puzzles } P} 1 \\
    = \#\{\inlinetri_{\sort(\alpha)\gamma,\beta\sort(\delta),\delta\alpha} \text{-puzzles}\} = \#\{\inlinepar_{\gamma,\beta,\delta,\alpha} \text{-puzzles}\}.   
    \end{multline*} 
    
    But the map $\psi$ is invariant if we permute the indeterminates with $\Phi_a$ or $\Phi_c$. 
 So we have 
    $$\psi\left(c_{\sort(\alpha)\gamma,\beta\sort(\delta)}^{(\delta\alpha)^\vee}\right) = \psi\left(\Phi_a \cdot c_{\sort(\alpha)\gamma,\beta\sort(\delta)}^{(\delta\alpha)^\vee}\right)= \psi\left(\Phi_c \cdot c_{\sort(\alpha)\gamma,\beta\sort(\delta)}^{(\delta\alpha)^\vee}\right) = \psi\left(\Phi_c\cdot \Phi_a \cdot c_{\sort(\alpha)\gamma,\beta\sort(\delta)}^{(\delta\alpha)^\vee}\right).$$ 
    Then we use Theorem \ref{thm:eqvtparallelogram}, to rewrite this as 
    $$\psi\left(c_{\sort(\alpha)\gamma,\beta\sort(\delta)}^{(\delta\alpha)^\vee}\right) = \psi\left(c_{\sort(\beta)\gamma,\alpha\sort(\delta)}^{(\delta\beta)^\vee}\right) = \psi\left(c_{\sort(\alpha)\delta,\beta\sort(\gamma)}^{(\gamma\alpha)^\vee}\right) = \psi\left(c_{\sort(\beta)\delta,\alpha\sort(\gamma)}^{(\gamma\beta)^\vee}\right).$$
    Since $\psi$ sends each structure constant to a count of the puzzles with the associated boundary labels, the claim immediately follows.  
\end{proof}

\textit{Example.} Figure \ref{fig:equivariantparallelogramsexample} at the beginning of this section gives an example of commuting the NW and SE labels (called $\gamma$ and $\delta$ in the theorems) for parallelogram-shaped equivariant puzzles. From left to right, the weights associated to the puzzles in the top row are $(y_4-y_1)(y_4-y_3)$, $(y_4-y_3)(y_6-y_3)$, and $(y_4-y_3)(y_5-y_2)$. This gives a structure constant of $(y_4-y_1)(y_4-y_3)+(y_4-y_3)(y_6-y_3)+(y_4-y_3)(y_5-y_2)$. From left to right, the weights associated to the puzzles in the bottom row are $(y_7-y_3)(y_6-y_3)$, $(y_5-y_1)(y_7-y_3)$, and $(y_7-y_3)(y_7-y_2)$. This gives a structure constant of $(y_7-y_3)(y_6-y_3)+(y_5-y_1)(y_7-y_3)+(y_7-y_3)(y_7-y_2)$, which can also be obtained by reversing the order of $y_4,y_5,y_6,y_7$ in the structure constant for the top row, as in Theorem \ref{thm:eqvtparallelogram}. The number of puzzles before and after commuting are also equal, as in Corollary \ref{cor:eqvtparallelogramnumber}.

\begin{remark}
    We do not have a nice commutativity result like Theorem \ref{thm:eqvtparallelogram} for equivariant trapezoidal puzzles, and it is not the case that commuting the labels on the two equal-length sides of a trapezoidal boundary should give the same effect on the structure constants as with parallelograms (i.e. permuting the $y_i$'s in this way). This is because the crucial factors of being able to write $\Phi_a\cdot [X_{\sort(\alpha)\gamma}]= [X_{\sort(\alpha)\gamma}]$ and  $\Phi_c\cdot [X_{\beta\sort(\delta)}]=  [X_{\beta\sort(\delta)}]$ in the proof of the theorem are particular to the special case of parallelograms. 
\end{remark}

\subsubsection{$d$-step version of Theorem \ref{thm:eqvtparallelogram}}

Here we will generalize Theorem \ref{thm:eqvtparallelogram} to the case of $d$-step flag manifolds. For this section, let $a_1\leq a_2\leq\cdots\leq a_d\leq a$ and $c_1\leq c_2\leq\cdots\leq c_d\leq c$ be non-negative integers, and let $\flag{F}$ and $\flag{\Ftilde}$ be the standard and anti-standard complete flags in $\C^{a+c}$, respectively.  

\begin{customthm}{5'}
    \label{thm:eqvtparallelogramdstep}

    Define $\Phi_a$ and $\Phi_c$ as in Theorem \ref{thm:eqvtparallelogram}.
Let $\alpha,\beta$ be strings with the content of \\ $\zerostr^{a-a_d}\onestr^{a_d-a_{d-1}}\twostr^{a_{d-1}-a_{d-2}}\cdots\dstr^{a_1}$, and let $\gamma,\delta$ be strings with the content of $\zerostr^{c-c_d}\onestr^{c_d-c_{d-1}}\twostr^{c_{d-1}-c_{d-2}}\cdots\dstr^{c_1}$. Then in $H_T^*(\Fl(a_1+c_1,a_2+c_2,\ldots,a_d+c_d; \C^{a+c}))$, we have 

$$(c_T)_{\sort(\alpha)\gamma,\beta\sort(\delta)}^{(\delta\alpha)^\vee}= \Phi_a\cdot (c_T)_{\sort(\beta)\gamma,\alpha\sort(\delta)}^{(\delta\beta)^\vee} = \Phi_c \cdot (c_T)_{\sort(\alpha)\delta,\beta\sort(\gamma)}^{(\gamma\alpha)^\vee} = \Phi_c\cdot\Phi_a\cdot (c_T)_{\sort(\beta)\delta,\alpha\sort(\gamma)}^{(\gamma\beta)^\vee}.$$
      
In other words commuting the pair $\alpha, \beta$ reverses the $y_1,\ldots,y_a$, and commuting the pair $\gamma,\delta$ reverses the $y_{a+1},\ldots,y_{a+c}$ in the structure constant.
\end{customthm} 

\begin{proof}

    The proof will follow the same outline as that of Theorem \ref{thm:eqvtparallelogram}, and the argument is mostly identical, with the only significant difference being in showing that $\Phi_a\cdot X_{\sort(\alpha)\gamma}= X_{\sort(\alpha)\gamma}$ and $\Phi_c\cdot X_{\beta\sort(\delta)}=  X_{\beta\sort(\delta)}$. To do that, we need the additional facts given in Section \ref{section:dstepschubertcalculus}, namely that $\flag{V}\in X_\lambda$ if and only if $p^i(\flag{V})\in p^i(X_\lambda)$ for all $1\leq i\leq d$, and $p^i(X_\lambda)=X_{\lambda^i}\subseteq \Gr(k_i;\C^n)$ in general. We note that here $(\sort(\alpha)\gamma)^i=\sort(\alpha^i)\gamma^i$ and $(\beta\sort(\delta))^i=\beta^i\sort(\delta^i)$ for all $1\leq i\leq d$. Then we can apply the same arguments written in the proof of Theorem \ref{thm:eqvtparallelogram} to prove that $\Phi_a\cdot X_{\sort(\alpha^i)\gamma^i}=X_{\sort(\alpha^i)\gamma^i}$ and $\Phi_c\cdot X_{\beta^i\sort(\delta^i)}=X_{\beta^i\sort(\delta^i)}$ for each $1\leq i\leq d$ individually. Putting the $d$ components together, this then implies that overall $\Phi_a\cdot X_{\sort(\alpha)\gamma}= X_{\sort(\alpha)\gamma}$ and $\Phi_c\cdot X_{\beta\sort(\delta)}=  X_{\beta\sort(\delta)}$.

    From there, we invoke Lemma \ref{lemma:commutesplitdstep} (the $d$-step generalization of Lemma \ref{lemma:commutesplit}) and obtain the same equations as are displayed in the proof of Theorem \ref{thm:eqvtparallelogram} to show that $\Phi_a\cdot (c_T)_{\sort(\beta)\gamma,\alpha\sort(\delta)}^{(\delta\beta)^\vee}= (c_T)_{\sort(\alpha)\gamma,\beta\sort(\delta)}^{(\delta\alpha)^\vee}$ and $\Phi_c \cdot (c_T)_{\sort(\alpha)\delta,\beta\sort(\gamma)}^{(\gamma\alpha)^\vee}= (c_T)_{\sort(\alpha)\gamma,\beta\sort(\delta)}^{(\delta\alpha)^\vee}$ (except that now of course the integrals are over $\Fl(a_1+c_1,a_2+c_2,\ldots,a_d+c_d; \C^{a+c})$), and then it is similarly straightforward to obtain $(c_T)_{\sort(\alpha)\gamma,\beta\sort(\delta)}^{(\delta\alpha)^\vee}= \Phi_c\cdot\Phi_a\cdot (c_T)_{\sort(\beta)\delta,\alpha\sort(\gamma)}^{(\gamma\beta)^\vee}$.
\end{proof}

\begin{remark}
    We obtain the same statement as in Corollary \ref{cor:eqvtparallelogramnumber} for 2-step and 3-step puzzles computing the $H^*_T$ structure constants. The proof strategy is the same.
\end{remark}

\subsection{Symmetric rhombus-shaped puzzles}

What we call a symmetric rhombus-shaped puzzle is a special case of a parallelogram-shaped puzzle where the labels on all four boundary sides have the same content. We will find that we can permute all the boundary labels in this case.

\begin{figure}[h]
    \centering
    \input{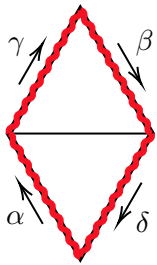}
    \caption{$\protect\inlinerhombus_{\alpha,\gamma,\beta,\delta}$, where the boundary labels all have matching content.}
    \label{fig:rhombusboundarycommutesolo}
\end{figure}

\subsubsection{Commutative property of symmetric rhombus-shaped puzzles}

\begin{theorem}[Commutative Property of Symmetric Rhombus-Shaped Puzzles]
    \label{thm:rhombuscommute}
    For $\alpha,\beta,\gamma,\delta\in\binom{[a]}{a_1}$, and using puzzle pieces in $\Hfrak$, we have
    $$ \#\{\inlinerhombus_{\alpha,\gamma,\beta,\delta} \text{-puzzles}\} =  \#\{\inlinerhombus_{f(\alpha),f(\gamma),f(\beta),f(\delta)} \text{-puzzles}\}$$
    for any bijection $f:\{\alpha,\beta,\gamma,\delta\} \rightarrow \{\alpha,\beta,\gamma,\delta\}$. In other words, we can permute all the boundary labels freely and preserve the number of puzzles.
\end{theorem}

\begin{proof}
We give two different proofs, one geometric and one combinatorial.

    \begin{enumerate}
        \item \textbf{Geometric proof.} If we apply Theorem \ref{thm:parallelogramLRnumbers} to this special case, we simply get 
        $$\#\{\inlinerhombus_{\alpha,\gamma,\beta,\delta} \text{-puzzles}\} = \sum_{\eta} c_{\alpha,\beta}^{\eta^\vee}\cdot c_{\gamma,\delta}^\eta.$$
        
        Then we write 
        \begin{multline*}
        \sum_{\eta} c_{\alpha,\beta}^{\eta^\vee}\cdot c_{\gamma,\delta}^\eta = \sum_\eta \underset{\Gr_a}{\int} [X_\alpha][X_\beta][X^{\eta^\vee}] \underset{\Gr_a}{\int} [X_\gamma][X_\delta][X^{\eta}] = \sum_\eta \underset{\Gr_a\times\Gr_a}{\int}[X_\alpha][X_\beta][X^{\eta^\vee}] \otimes [X_\gamma][X_\delta][X^{\eta}]\\
        = \underset{\Gr_a\times\Gr_a}{\int}([X_\alpha][X_\beta] \otimes [X_\gamma][X_\delta]) \left(\sum_\eta [X^{\eta^\vee}]\otimes [X^{\eta}]\right) = \underset{\Gr_a\times\Gr_a}{\int}([X_\alpha][X_\beta] \otimes [X_\gamma][X_\delta]) \Delta_*([\Gr_a])\\
        =\underset{\Gr_a}{\int} \Delta^*([X_\alpha][X_\beta] \otimes [X_\gamma][X_\delta])[\Gr_a] =\underset{\Gr_a}{\int} ([X_\alpha][X_\beta] [X_\gamma][X_\delta])[\Gr_a] =\underset{\Gr_a}{\int} [X_\alpha][X_\beta] [X_\gamma][X_\delta],
        \end{multline*}
        
        where $\Delta:\Gr_a\hookrightarrow \Gr_a\times \Gr_a$ is the diagonal inclusion map $V\mapsto (V,V)$, hence $\Delta^*$ is the cup product. We can clearly permute $\alpha,\beta,\gamma,\delta$ in the last expression since the cohomology ring is commutative. 

        That finishes the proof, but if the fact that $\Delta_*([\Gr_a])= \sum\limits_\eta [X^{\eta^\vee}]\otimes [X^{\eta}]$ is unclear, we explain it as follows. If we expand $\Delta_*([\Gr_a])$ in the basis $\left\{[X^\xi]\otimes [X^\eta]: \xi,\eta\in \binom{[a]}{a_1}\right\}$, the coefficient on the $[X^\xi]\otimes [X^\eta]$ term is
        \begin{multline*}
            \underset{\Gr_a\times\Gr_a}{\int} \Delta_*([\Gr_a])([X_\xi]\otimes [X_\eta]) = \underset{\Gr_a}{\int} [\Gr_a]\Delta^*([X_\xi]\otimes [X_\eta]) = \underset{\Gr_a}{\int} [\Gr_a][X_\xi][X_\eta] \\
            = \underset{\Gr_a}{\int} 1\cdot [X_\xi] [X_\eta] = \underset{\Gr_a}{\int} [X_\xi] [X_\eta]= \underset{\Gr_a}{\int} [X_\xi][X^{\eta^\vee}] = \langle \xi,\eta^\vee\rangle.
        \end{multline*}

        \item \textbf{Proof by commutative property of parallelogram-shaped puzzles.} A rhombus-shaped puzzle with symmetric content can be seen as two triangular puzzles glued together. We can use the commutativity of triangular puzzles to commute the pairs $\gamma,\beta$ and $\alpha,\delta$, and we can use the commutative property of parallelogram-shaped puzzles (Theorem \ref{thm:Hparallelogramcommute}) to commute the pairs $\alpha,\beta$ and $\gamma,\delta$. Combining these two operations gives all possible permutations of the boundary labels.
    \end{enumerate}
\end{proof}

\begin{remark}
    Theorem \ref{thm:rhombuscommute} should already be deducible from the basic commutative and associative properties of triangular puzzles, but it has not been stated in the literature as far as the author is aware. The associativity of the ring $H^*(\Gr(k;\C^n))$ can be expressed by the condition that
    $$\sum_\nu c_{\gamma,\beta}^{\nu}\cdot c_{\nu,\delta}^{\alpha^\vee} = \sum_\nu c_{\beta,\delta}^{\nu}\cdot c_{\gamma,\nu}^{\alpha^\vee},$$
    which, when seen as summing pairs of triangular puzzles glued together, implies
    $$ \#\{\inlinerhombus_{\alpha,\gamma,\beta,\delta} \text{-puzzles}\} =  \#\{\inlinerhombus_{\gamma,\beta,\delta,\alpha} \text{-puzzles}\}.$$
    (In fact the associativity of puzzles is proved combinatorially in terms of rhombus-shaped puzzles in \cite[\S 3.4]{Purbhoo}.) This combined with the ability to commute the pairs $\gamma,\beta$ and $\alpha,\delta$ from the commutative property of triangular puzzles gives us all possible permutations of the labels.

\end{remark}

\section{Hexagons}
\label{section:hexagons}

For the entirety of this section, let $a_0,a_1,b_0,b_1,c_0,c_1,d_0,d_1,e_0,e_1,z_0,z_1$ be nonnegative integers such that $\zerostr^{b_0+c_0}\onestr^{b_1+c_1}=\zerostr^{e_0+z_0}\onestr^{e_1+z_1}$, $\zerostr^{a_0+b_0}\onestr^{a_1+b_1}=\zerostr^{d_0+e_0}\onestr^{d_1+e_1}$, and $\zerostr^{c_0+d_0}\onestr^{c_1+d_1}=\zerostr^{z_0+a_0}\onestr^{z_1+a_1}$, and define $a:=a_0+a_1$, $b:=b_0+b_1$, $c:=c_0+c_1$, $d:=d_0+d_1$, $e:=e_0+e_1$, and $z:=z_0+z_1$.

We will be looking at hexagonal puzzles with boundary $\inlinehex_{\alpha,\beta,\gamma,\delta,\epsilon,\zeta}$, where $\alpha\in\binom{[a]}{a_1}$, $\beta\in\binom{[b]}{b_1}$, $\gamma\in\binom{[c]}{c_1}$, $\delta\in\binom{[d]}{d_1}$, $\epsilon\in\binom{[e]}{e_1}$, and $\zeta\in\binom{[z]}{z_1}$. The symmetries in the content of these strings described above are necessary and sufficient to say that the boundary is in fact a valid hexagon and can possibly be filled with puzzle pieces. (This is both a matter of hexagon side lengths and the requirements of the discrete Green's theorem discussed in Section \ref{section:otherpuzzleproperties}. The relevant statement under the Green's theorem property is (c), combined with Proposition \ref{prop:hexagoncomplete} below.)

In this section we will give a geometric proof of $180^\circ$ rotational symmetry of hexagonal puzzles, and we will prove commutative properties for the four different types of symmetries that a hexagonal puzzle can possess. We specialize one of these to prove a commutative property for symmetric pentagonal puzzles. In the cases where it is useful or interesting to do so, we also provide a formula for the number of puzzles in terms of simpler Littlewood-Richardson numbers.


\subsection{Geometric interpretation of hexagonal puzzles}

In this section we give a geometric interpretation to the most general kind of hexagonal puzzle by completing to a triangular puzzle as in previous sections, and we use this to derive a formula which will be specialized further as we move on to hexagonal puzzles with special symmetries.

\subsubsection{Completion of a hexagon to a triangle}

\begin{proposition}
    \label{prop:hexagoncomplete}
    For $\alpha\in\binom{[a_0+a_1]}{a_1}$, $\beta\in\binom{[b_0+b_1]}{b_1}$, $\gamma\in\binom{[c_0+c_1]}{c_1}$, $\delta\in\binom{[d_0+d_1]}{d_1}$, $\epsilon\in\binom{[e_0+e_1]}{e_1}$, and $\zeta\in\binom{[z_0+z_1]}{z_1}$, and using puzzle pieces in $\Hfrak \cup \left\{\raisebox{-.4\height}{
\begin{tikzpicture}[x=1pt,y=1pt,yscale=1,xscale=1,scale=.68]
\draw [fill={rgb, 255:red, 255; green, 246; blue, 74 }] (0.0,0.0) -- (15.0,25.98) -- (30.0,0.0) -- (15.0,-25.98) -- cycle;
\draw (7.5,12.99) node  [rotate=60]  {$\mathtt{0}$};
\draw (22.5,12.99) node  [rotate=-60]  {$\mathtt{1}$};
\draw (22.5,-12.99) node  [rotate=60]  {$\mathtt{0}$};
\draw (7.5,-12.99) node  [rotate=-60]  {$\mathtt{1}$};
\end{tikzpicture}
},\raisebox{-.4\height}{
\begin{tikzpicture}[x=1pt,y=1pt,yscale=1,xscale=1,scale=.68]
\draw [fill={rgb, 255:red, 255; green, 255; blue, 255 }] (0.0,0.0) -- (15.0,25.98) -- (30.0,0.0) -- cycle;
\draw (7.5,12.99) node  [rotate=60]  {$\mathtt{1\! 0}$};
\draw (22.5,12.99) node  [rotate=-60]  {$\mathtt{1\! 0}$};
\draw (15.0,0.0) node  [rotate=0]  {$\mathtt{1\! 0}$};
\end{tikzpicture}
} \right\}$, we have a bijection
    $$\{\inlinehex_{\alpha,\beta,\gamma,\delta,\epsilon,\zeta} \text{-puzzles}\} \leftrightarrow \{\inlinetri_{\sort(\alpha)\beta\gamma,\sort(\gamma)\delta\sort(\epsilon),\epsilon\zeta\alpha} \text{-puzzles}\}$$
\end{proposition}

\begin{figure*}[h]
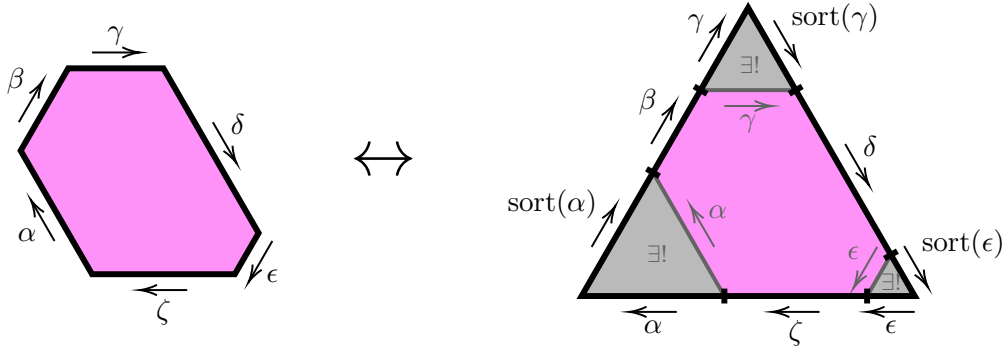

    \centering
    \raisebox{-.5\height}{\input{proof_boundary_diagrams/complete_to_triangle/hexagonfilling}} 
    \Huge $\quad \leftrightarrow \quad$ \normalsize \raisebox{-.5\height}{\input{proof_boundary_diagrams/complete_to_triangle/hexagoncomplete}}
    
    \caption{Puzzles with the boundary $\protect\inlinehex_{\alpha,\beta,\gamma,\delta,\epsilon,\zeta}$ are in bijection with puzzles with the boundary $\protect\inlinetri_{\sort(\alpha)\beta\gamma,\sort(\gamma)\delta\sort(\epsilon),\epsilon\zeta\alpha}$. On the righthand side, given the outer boundary labels, there exists a unique filling of each grey triangular region, and these force the labels around the inner pink region to replicate those for the hexagon itself (see Property 1 of Section \ref{section:otherpuzzleproperties}). The bijection is to glue on the uniquely filled grey regions, and the inverse is to cut them off.}
    \label{fig:hexagoncomplete}
\end{figure*}

\begin{proof}
    See Figure \ref{fig:hexagoncomplete}.
\end{proof}

\subsubsection{A formula in terms of Littlewood-Richardson numbers}



\begin{proposition}
    \label{prop:hexagonstructureconstantformula}
    For $\alpha\in\binom{[a_0+a_1]}{a_1}$, $\beta\in\binom{[b_0+b_1]}{b_1}$, $\gamma\in\binom{[c_0+c_1]}{c_1}$, $\delta\in\binom{[d_0+d_1]}{d_1}$, $\epsilon\in\binom{[e_0+e_1]}{e_1}$, and $\zeta\in\binom{[z_0+z_1]}{z_1}$, and using puzzle pieces in $\Hfrak$, we have
    $$\#\{\inlinehex_{\alpha,\beta,\gamma,\delta,\epsilon,\zeta} \text{-puzzles}\}  =\sum_{\substack{\mu,\theta \text{ such that}\\ \zerostr^{a_0}\mu\onestr^{a_1} = \zerostr^{c_0}\theta\onestr^{e_1}}} c_{\zerostr^{b_0}\alpha^\vee\onestr^{z_1},\theta}^{\onestr^{c_1}\delta\zerostr^{e_0}}\cdot c_{\beta\gamma,\mu}^{(\epsilon\zeta)^\vee}$$
\end{proposition}

\begin{proof}
Once completed to a triangular puzzle as in Proposition \ref{prop:hexagoncomplete}, a hexagonal puzzle can be thought of as a special case of a trapezoidal puzzle (compare with Figure \ref{fig:trapezoidcomplete}). We can apply Proposition \ref{prop:trapezoidgeometric}, letting $\alpha$, $\beta\gamma$, $\sort(\gamma)\delta\sort(\epsilon)$, and $\epsilon\zeta$ replace $\beta$, $\gamma$, $\nu$, and $\delta$, respectively.
This gives us
$$\#\{\inlinehex_{\alpha,\beta,\gamma,\delta,\epsilon,\zeta} \text{-puzzles}\} = \#\{\inlinetrapside_{\alpha,\beta\gamma,\sort(\gamma)\delta\sort(\epsilon),\epsilon\zeta}\text{-puzzles}\} = \sum_{\mu}c_{(\alpha^\vee)\pad,\mu\pad}^{\sort(\gamma)\delta\sort(\epsilon)} \cdot c_{\beta\gamma,\mu}^{(\epsilon\zeta)^\vee}.$$

Now let us look more closely at $c_{(\alpha^\vee)\pad,\mu\pad}^{\sort(\gamma)\delta\sort(\epsilon)}$. We have $(\alpha^\vee)\pad=\zerostr^{b_0+c_0}\alpha^\vee\onestr^{e_1+z_1}= \zerostr^{c_0}\zerostr^{b_0}\alpha^\vee\onestr^{z_1}\onestr^{e_1}$, $\mu\pad=\zerostr^{a_0}\mu\onestr^{a_1}$, and $\sort(\gamma)\delta\sort(\epsilon) =\zerostr^{c_0}\onestr^{c_1}\delta\zerostr^{e_0}\onestr^{e_1}$. Then,
$$c_{(\alpha^\vee)\pad,\mu\pad}^{\sort(\gamma)\delta\sort(\epsilon)} = \underset{\Gr(a_1+z_1+e_1,\C^{a+z+e})}{\int} [X_{\zerostr^{c_0}\zerostr^{b_0}\alpha^\vee\onestr^{z_1}\onestr^{e_1}}][X_{\zerostr^{a_0}\mu\onestr^{a_1}}][X^{\zerostr^{c_0}\onestr^{c_1}\delta\zerostr^{e_0}\onestr^{e_1}}].$$

Now we will use a similar argument as in the geometric proof of Theorem \ref{thm:parallelogramLRnumbers} and so will include less detail this time. We note that for the above integral to be nonzero, $\zerostr^{a_0}\mu\onestr^{a_1}$ must begin with $\zerostr^{c_0}$ and end with $\onestr^{e_1}$. So let us assume this and define $\theta$ so that $\zerostr^{a_0}\mu\onestr^{a_1} = \zerostr^{c_0}\theta\onestr^{e_1}$.

Let $\flag{F}$ and $\flag{\Ftilde}$ be the standard and anti-standard complete flags. Now define the maps 
$$ 
\Psi: \Gr(0,F_{c_0})\times \Gr(a_1+z_1,\langle \bfe_{c_0+1},\ldots,\bfe_{c_0+b_0+a+z_1}\rangle) \times \Gr(e_1,\Ftilde_{e_1}) \hookrightarrow \Gr(a_1+z_1,F_{c_0+b_0+a+z_1}) \times \Gr(e_1,\Ftilde_{e_1}), $$ $$
(U,V,W)=(\{\mathbf{0}\},V,\Ftilde_{e_1})\mapsto (U\oplus V,W) =(\{\mathbf{0}\} \oplus V, \Ftilde_{e_1})
$$
and 
\begin{gather*}
    \Omega:\Gr(a_1+z_1,F_{c_0+b_0+a+z_1}) \times \Gr(e_1,\Ftilde_{e_1}) \hookrightarrow \Gr(a_1+z_1+e_1,\C^{a+z+e}),\\
(V,W)= (V,\Ftilde_{e_1}) \mapsto V\oplus W= V\oplus \Ftilde_{e_1}.
\end{gather*}

Going forward, let us shorten the notation for these spaces to $\Gr_{c_0}:=\Gr(0,F_{c_0})$, $\Gr_{b_0+a+z_1}:=\Gr(a_1+z_1,\langle \bfe_{c_0+1},\ldots,\bfe_{c_0+b_0+a+z_1}\rangle)$, $\Gr_{e_1}:=\Gr(e_1,\Ftilde_{e_1})$, and $\Gr_{a+z+e}:=\Gr(a_1+z_1+e_1,\C^{a+z+e})$.

Now we have
$$[X_{\zerostr^{c_0}\zerostr^{b_0}\alpha^\vee\onestr^{z_1}\onestr^{e_1}}][X^{\zerostr^{c_0}\onestr^{c_1}\delta\zerostr^{e_0}\onestr^{e_1}}] = (\Omega\circ \Psi)_*\left([X_{\zerostr^{c_0}}]\otimes [X_{\zerostr^{b_0}\alpha^\vee\onestr^{z_1}}][X^{\onestr^{c_1}\delta\zerostr^{e_0}}] \otimes [X_{\onestr^{e_1}}]\right).$$ 

Then we continue with
\begin{align*}
c_{(\alpha^\vee)\pad,\zerostr^{c_0}\theta\onestr^{e_1}}^{\sort(\gamma)\delta\sort(\epsilon)} &=\underset{\Gr_{a+z+e}}{\int} [X_{\zerostr^{c_0}\zerostr^{b_0}\alpha^\vee\onestr^{z_1}\onestr^{e_1}}][X_{\zerostr^{c_0}\theta\onestr^{e_1}}][X^{\zerostr^{c_0}\onestr^{c_1}\delta\zerostr^{e_0}\onestr^{e_1}}]\\
&=\underset{\Gr_{a+z+e}}{\int} (\Omega\circ \Psi)_*\left([X_{\zerostr^{c_0}}]\otimes [X_{\zerostr^{b_0}\alpha^\vee\onestr^{z_1}}][X^{\onestr^{c_1}\delta\zerostr^{e_0}}] \otimes [X_{\onestr^{e_1}}]\right)[X_{\zerostr^{c_0}\theta\onestr^{e_1}}]\\
&=\underset{\Gr_{c_0}\times \Gr_{b_0+a+z_1}\times \Gr_{e_1}}{\int} \left([X_{\zerostr^{c_0}}]\otimes [X_{\zerostr^{b_0}\alpha^\vee\onestr^{z_1}}][X^{\onestr^{c_1}\delta\zerostr^{e_0}}] \otimes [X_{\onestr^{e_1}}]\right)(\Omega\circ \Psi)^*\left([X_{\zerostr^{c_0}\theta\onestr^{e_1}}]\right)\\
&=\underset{\Gr_{c_0}\times \Gr_{b_0+a+z_1}\times \Gr_{e_1}}{\int} \left([X_{\zerostr^{c_0}}]\otimes [X_{\zerostr^{b_0}\alpha^\vee\onestr^{z_1}}][X^{\onestr^{c_1}\delta\zerostr^{e_0}}] \otimes [X_{\onestr^{e_1}}]\right)\left([X_{\zerostr^{c_0}}]\otimes [X_\theta]\otimes [X_{\onestr^{e_1}}]\right)\\
&=\underset{\Gr_{c_0}\times \Gr_{b_0+a+z_1}\times \Gr_{e_1}}{\int} [X_{\zerostr^{c_0}}][X_{\zerostr^{c_0}}]\otimes [X_{\zerostr^{b_0}\alpha^\vee\onestr^{z_1}}][X^{\onestr^{c_1}\delta\zerostr^{e_0}}][X_\theta] \otimes [X_{\onestr^{e_1}}][X_{\onestr^{e_1}}]\\
&=\underset{\Gr_{c_0}}{\int} [X_{\zerostr^{c_0}}][X_{\zerostr^{c_0}}] \underset{\Gr_{b_0+a+z_1}}{\int}  [X_{\zerostr^{b_0}\alpha^\vee\onestr^{z_1}}][X^{\onestr^{c_1}\delta\zerostr^{e_0}}][X_\theta] \underset{\Gr_{e_1}}{\int} [X_{\onestr^{e_1}}][X_{\onestr^{e_1}}]\\
&=\underset{\Gr_{b_0+a+z_1}}{\int}  [X_{\zerostr^{b_0}\alpha^\vee\onestr^{z_1}}][X^{\onestr^{c_1}\delta\zerostr^{e_0}}][X_\theta]\\
&=c_{\zerostr^{b_0}\alpha^\vee\onestr^{z_1},\theta}^{\onestr^{c_1}\delta\zerostr^{e_0}}.
\end{align*}

To summarize, we now have that $c_{(\alpha^\vee)\pad,\mu\pad}^{\sort(\gamma)\delta\sort(\epsilon)} = c_{\zerostr^{b_0}\alpha^\vee\onestr^{z_1},\theta}^{\onestr^{c_1}\delta\zerostr^{e_0}}$ when $\zerostr^{a_0}\mu\onestr^{a_1} = \zerostr^{c_0}\theta\onestr^{e_1}$.

Finally, we return to where we began, substituting in the above equality. We have 
$$\#\{\inlinehex_{\alpha,\beta,\gamma,\delta,\epsilon,\zeta} \text{-puzzles}\} = \sum_{\mu} c_{(\alpha^\vee)\pad,\mu\pad}^{\sort(\gamma)\delta\sort(\epsilon)} \cdot c_{\beta\gamma,\mu}^{(\epsilon\zeta)^\vee} = \sum_{\substack{\mu,\theta \text{ such that}\\ \zerostr^{a_0}\mu\onestr^{a_1} = \zerostr^{c_0}\theta\onestr^{e_1}}} c_{\zerostr^{b_0}\alpha^\vee\onestr^{z_1},\theta}^{\onestr^{c_1}\delta\zerostr^{e_0}}\cdot  c_{\beta\gamma,\mu}^{(\epsilon\zeta)^\vee}.$$
\end{proof}


\subsection{$180^\circ$ rotational symmetry of hexagonal puzzles}

\begin{figure}[h]
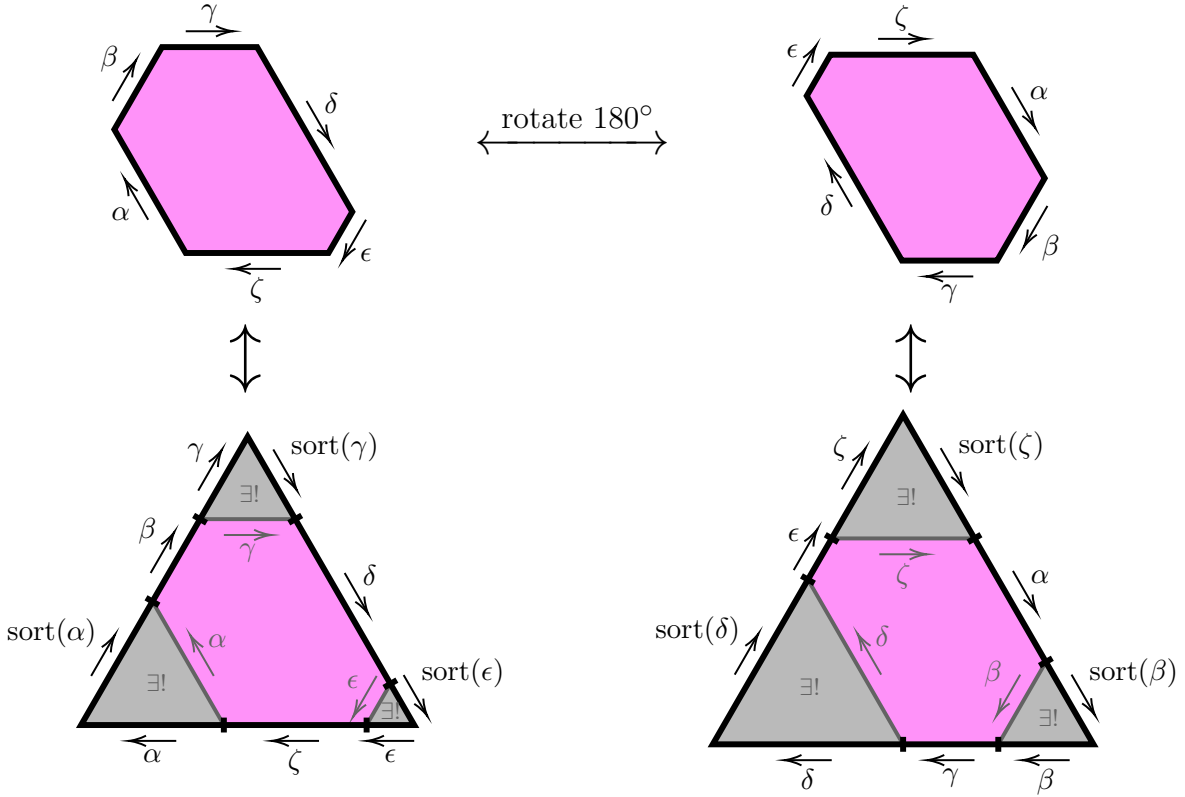

    \centering
    \raisebox{-.5\height}{\input{proof_boundary_diagrams/complete_to_triangle/hexagonfilling}} 
    \LARGE $\quad\quad\, \xleftrightarrow{\text{ rotate }180^\circ} \quad\quad\,$ \normalsize \raisebox{-.5\height}{\input{proof_boundary_diagrams/complete_to_triangle/hexagonrotatefilling}} \;\;

    \vspace{.5\baselineskip}
    \Huge $\updownarrow \qquad\qquad\qquad\qquad\qquad\;\, \updownarrow$  \normalsize
    \vspace{.5\baselineskip}
    
    \raisebox{-.5\height}{\input{proof_boundary_diagrams/complete_to_triangle/hexagoncomplete}} 
    \Huge $\; \quad \, \;$ \normalsize \raisebox{-.5\height}{\input{proof_boundary_diagrams/complete_to_triangle/hexagonrotatecomplete}}

    \caption{Composing the bijections of puzzles shown, when using pieces in $\Hfrak$, implies an equality of Littlewood-Richardson numbers $c_{\sort(\alpha)\beta\gamma,\, \sort(\gamma)\delta\sort(\epsilon)}^{(\epsilon\zeta\alpha)^\vee} = c_{\sort(\delta)\epsilon\zeta,\, \sort(\gamma)\alpha\sort(\epsilon)}^{(\beta\gamma\delta)^\vee}$. }
    \label{fig:hexagonrotate}
\end{figure}

When using any set of puzzle pieces preserved under $180^\circ$ rotation, $180^\circ$ rotation of hexagonal puzzles gives a bijection
$$\{\inlinehex_{\alpha,\beta,\gamma,\delta,\epsilon,\zeta} \text{-puzzles}\} \leftrightarrow \{\inlinehex_{\delta,\epsilon,\zeta,\alpha,\beta,\gamma,} \text{-puzzles}\}.$$ 
We will now proceed by focusing on the case where we use only puzzle pieces in $\Hfrak$. Note that if we complete the hexagonal puzzles to trianglular puzzles on either side of the above bijection, we get a bijection
$$\{\inlinetri_{\sort(\alpha)\beta\gamma,\, \sort(\zeta)\delta\sort(\beta),\,\epsilon\zeta\alpha}\text{-puzzles}\} \leftrightarrow \{\inlinetri_{\sort(\delta)\epsilon\zeta,\, \sort(\gamma)\alpha\sort(\epsilon),\,\beta\gamma\delta}\text{-puzzles}\}.$$
(See Figure \ref{fig:hexagonrotate}.)
This yields an equality of Littlewood-Richardson numbers which has a geometric interpretation as well, and in the following proposition we give a geometric proof of it. It is worth noting that, as shown in Figure \ref{fig:hexagonrotate}, the triangular puzzles are different sizes, which means the Schubert calculus occurs in different Grassmannians. This adds some additional difficulty to the geometric analysis compared to the case of parallelograms.


\begin{proposition}
    \label{hexagonupsidedown}
    For $\alpha\in\binom{[a]}{a_1}$, $\beta\in\binom{[b]}{b_1}$, $\gamma\in\binom{[c]}{c_1}$, $\delta\in\binom{[d]}{d_1}$, $\epsilon\in\binom{[e]}{e_1}$, and $\zeta\in\binom{[z]}{z_1}$, we have
    $$c_{\sort(\alpha)\beta\gamma,\, \sort(\gamma)\delta\sort(\epsilon)}^{(\epsilon\zeta\alpha)^\vee} = c_{\sort(\delta)\epsilon\zeta,\, \sort(\gamma)\alpha\sort(\epsilon)}^{(\beta\gamma\delta)^\vee}.$$
\end{proposition}

\begin{proof}
We give two different arguments.

    \begin{enumerate}
        \item \textbf{Proof by $180^\circ$ rotational symmetry of hexagonal puzzles.} This is trivial to see from Figure \ref{fig:hexagonrotate}.

        \item \textbf{Geometric proof.} 

    By completing to a triangle as in Proposition \ref{prop:hexagoncomplete}, we get that proving this equality of Littlewood-Richardson numbers is equivalent to proving that 
    $$\#\{\inlinehex_{\alpha,\beta,\gamma,\delta,\epsilon,\zeta} \text{-puzzles}\} = \#\{\inlinehex_{\delta,\epsilon,\zeta\alpha,\beta,\gamma} \text{-puzzles}\},$$

using pieces in $\Hfrak$. We apply Proposition \ref{prop:hexagonstructureconstantformula} to the situation of the $180^\circ$ rotated puzzle, and we have 
 $$\#\{\inlinehex_{\delta,\epsilon,\zeta\alpha,\beta,\gamma} \text{-puzzles}\} =\sum_{\substack{\mu,\theta \text{ such that}\\ \zerostr^{d_0}\mu\onestr^{d_1} = \zerostr^{z_0}\theta\onestr^{b_1}}} c_{\zerostr^{e_0}\delta^\vee\onestr^{c_1},\theta}^{\onestr^{z_1}\alpha\zerostr^{b_0}}\cdot c_{\epsilon\zeta,\mu}^{(\beta\gamma)^\vee}.$$

 We use that $c_{\zerostr^{e_0}\delta^\vee\onestr^{c_1},\theta}^{\onestr^{z_1}\alpha\zerostr^{b_0}} = c_{\zerostr^{b_0}\alpha^\vee\onestr^{z_1},\theta}^{\onestr^{c_1}\delta\zerostr^{e_0}}$ and $c_{\epsilon\zeta,\mu}^{(\beta\gamma)^\vee}= c_{\beta\gamma,\mu}^{(\epsilon\zeta)^\vee}$ by basic properties. We also have that $\zerostr^{a_0}\mu\onestr^{a_1} = \zerostr^{c_0}\theta\onestr^{e_1}$ if and only if $\zerostr^{d_0}\mu\onestr^{d_1}= \zerostr^{z_0}\theta\onestr^{b_1}$, since $a_0-c_0=d_0-z_0$ and $a_1-e_1=d_1-b_1$ by definition (see the beginning of Section \ref{section:hexagons}). Making these substitutions in the expression above, we exactly get
 
 $$\#\{\inlinehex_{\delta,\epsilon,\zeta\alpha,\beta,\gamma} \text{-puzzles}\} =\sum_{\substack{\mu,\theta \text{ such that}\\ \zerostr^{a_0}\mu\onestr^{a_1} = \zerostr^{c_0}\theta\onestr^{e_1}}} c_{\zerostr^{b_0}\alpha^\vee\onestr^{z_1},\theta}^{\onestr^{c_1}\delta\zerostr^{e_0}} \cdot c_{\beta\gamma,\mu}^{(\epsilon\zeta)^\vee},$$
 
 which is the same expression as the one for the non-rotated version, $\#\{\inlinehex_{\alpha,\beta,\gamma,\delta,\epsilon,\zeta} \text{-puzzles}\}$, according to Proposition \ref{prop:hexagonstructureconstantformula}.
    \end{enumerate}
\end{proof}

\begin{remark}
    What is interesting about Proposition \ref{fig:hexagonrotate} is the juxtaposition of the two given proofs. One is based on a trivial observation on the puzzle level, and the other provides an understanding of why the symmetry should be true on a geometric level once we correspond it to a Schubert calculus problem, which is obviously much more involved and interesting. In fact, searching for this geometric manifestation of the $180^\circ$ rotational symmetry of hexagonal puzzles was the original question that led to the rest of the work in this paper. We discovered all of these commutativity results as an unexpected product of that exploration.
\end{remark}

\subsection{Hexagonal puzzles with opposite sides symmetry}

\begin{figure}[h]
    \centering
    \input{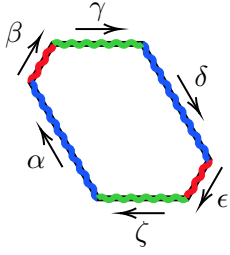}   \caption{$\protect\inlinehex_{\alpha,\beta,\gamma,\delta,\epsilon,\zeta}$, displaying ``opposite sides symmetry.''}
    \label{fig:hexagonboundarycommuteoppositesolo}
\end{figure}

For hexagonal puzzles, what we will refer to as ``opposite sides symmetry'' is a situation where each pair of opposite (parallel) boundary sides has labels with matching content. (See Figure \ref{fig:hexagonboundarycommuteoppositesolo}.)
We will find that we can commute the labels on these pairs independently.

\subsubsection{Commutative property of hexagonal puzzles with opposite sides symmetry}

\begin{theorem}[Commutative Property of Hexagonal Puzzles with Opposite Sides Symmetry]
    \label{thm:hexagonopposite}
    Let $\alpha\in\binom{[a]}{a_1}$, $\beta\in\binom{[b]}{b_1}$, $\gamma\in\binom{[c]}{c_1}$, $\delta\in\binom{[d]}{d_1}$, $\epsilon\in\binom{[e]}{e_1}$, and $\zeta\in\binom{[z]}{z_1}$. If $\sort(\alpha)=\sort(\delta)$ (and hence $\sort(\beta)=\sort(\epsilon)$ and $\sort(\gamma)=\sort(\zeta)$ also), then, using puzzle pieces in $\Hfrak$, we have
    \begin{multline*}
        \#\{\inlinehex_{\alpha,\beta,\gamma,\delta,\epsilon,\zeta} \text{-puzzles}\} = \#\{\inlinehex_{\delta,\beta,\gamma,\alpha,\epsilon,\zeta} \text{-puzzles}\} = \#\{\inlinehex_{\alpha,\epsilon,\gamma,\delta,\beta,\zeta} \text{-puzzles}\} = \#\{\inlinehex_{\alpha,\beta,\zeta,\delta,\epsilon,\gamma} \text{-puzzles}\} \\
        = \#\{\inlinehex_{\delta,\epsilon,\gamma,\alpha,\beta,\zeta} \text{-puzzles}\} = \#\{\inlinehex_{\delta,\beta,\zeta,\alpha,\epsilon,\gamma} \text{-puzzles}\} = \#\{\inlinehex_{\alpha,\epsilon,\zeta,\delta,\beta,\gamma} \text{-puzzles}\} = \#\{\inlinehex_{\delta,\epsilon,\zeta,\alpha,\beta,\gamma} \text{-puzzles}\}
    \end{multline*}
    In other words, we can commute each pair $(\alpha,\delta)$, $(\beta,\epsilon)$, and $(\gamma,\zeta)$ independently and preserve the number of puzzles.
\end{theorem}

\begin{proof}

   We will begin by proving the first equality
   $$\#\{\inlinehex_{\alpha,\beta,\gamma,\delta,\epsilon,\zeta} \text{-puzzles}\} = \#\{\inlinehex_{\delta,\beta,\gamma,\alpha,\epsilon,\zeta} \text{-puzzles}\}.$$
    
    Using Proposition \ref{prop:hexagonstructureconstantformula} and the basic fact that $c_{\zerostr^{b_0}\alpha^\vee\onestr^{z_1},\theta}^{\onestr^{c_1}\delta\zerostr^{e_0}} =c_{\zerostr^{e_0}\delta^\vee\onestr^{c_1},\theta}^{\onestr^{z_1}\alpha\zerostr^{b_0}}$, we write 

    $$\#\{\inlinehex_{\alpha,\beta,\gamma,\delta,\epsilon,\zeta} \text{-puzzles}\} =\sum_{\substack{\mu,\theta \text{ such that}\\ \zerostr^{a_0}\mu\onestr^{a_1} = \zerostr^{c_0}\theta\onestr^{e_1}}} c_{\zerostr^{b_0}\alpha^\vee\onestr^{z_1},\theta}^{\onestr^{c_1}\delta\zerostr^{e_0}}\cdot c_{\beta\gamma,\mu}^{(\epsilon\zeta)^\vee} = \sum_{\substack{\mu,\theta \text{ such that}\\ \zerostr^{a_0}\mu\onestr^{a_1} = \zerostr^{c_0}\theta\onestr^{e_1}}} c_{\zerostr^{e_0}\delta^\vee\onestr^{c_1},\theta}^{\onestr^{z_1}\alpha\zerostr^{b_0}}\cdot c_{\beta\gamma,\mu}^{(\epsilon\zeta)^\vee}.$$
    
    Given the symmetries in the content of the strings, we make some substitutions of $a_0=d_0$, $a_1=d_1$, $b_0=e_0$, and $c_1=z_1$ in the above expression, and we end up with 

    $$\#\{\inlinehex_{\alpha,\beta,\gamma,\delta,\epsilon,\zeta} \text{-puzzles}\} = \sum_{\substack{\mu,\theta \text{ such that}\\ \zerostr^{d_0}\mu\onestr^{d_1} = \zerostr^{c_0}\theta\onestr^{e_1}}} c_{\zerostr^{b_0}\delta^\vee\onestr^{z_1},\theta}^{\onestr^{c_1}\alpha\zerostr^{e_0}}\cdot c_{\beta\gamma,\mu}^{(\epsilon\zeta)^\vee},$$

    which is exactly the expression for $\#\{\inlinehex_{\delta,\beta,\gamma,\alpha,\epsilon,\zeta} \text{-puzzles}\}$, according to Proposition \ref{prop:hexagonstructureconstantformula}. This completes the first step of proving that $\alpha$ and $\delta$ commute.

    The rest of the claim of the theorem, which is that the pairs $\beta,\epsilon$ and $\gamma,\zeta$ also commute independently, just follows from $120^\circ$ rotational symmetry of puzzles. We can rotate our labeled hexagonal boundary so that another of these pairs takes the place of $\alpha$ and $\delta$, complete to a triangle as in Proposition \ref{prop:hexagoncomplete}, use the same argument as above, and then undo the rotation. Conjugating the process we used to commute the pair $\alpha,\delta$ with a $120^\circ$ or $240^\circ$ rotation allows us to commute the pairs $\beta,\epsilon$ or $\gamma,\zeta$, respectively. Combining these operations produces the $2^3=8$ different configurations written in the theorem statement.    
\end{proof}

\subsection{Hexagonal puzzles with two-way symmetry}

\begin{figure}[h]
    \centering
    \input{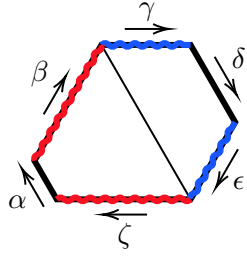}   \caption{$\protect\inlinehex_{\alpha,\beta,\gamma,\delta,\epsilon,\zeta}$, displaying ``two-way symmetry.''}
    \label{fig:hexagonboundarycommutetwowaysolo}
\end{figure}

For hexagonal puzzles, what we will refer to as ``two-way symmetry'' is a situation where there are two pairs of boundary sides whose labels have matching content, and where there is an angular difference of $60^\circ$ between the sides in each pair. (See Figure \ref{fig:hexagonboundarycommutetwowaysolo}.) (In fact, if one such pair of boundary sides has labels with matching content, then the other pair must also, by the definitions given in the first paragraph of Section \ref{section:hexagons}.)
We will find that we can commute the labels on these pairs independently.

\subsubsection{Commutative property of hexagonal puzzles with two-way symmetry}

\begin{theorem}[Commutative Property of Hexagonal Puzzles with Two-Way Symmetry]
    \label{thm:hexagontwoway}
    Let $\alpha\in\binom{[a]}{a_1}$, $\beta\in\binom{[b]}{b_1}$, $\gamma\in\binom{[c]}{c_1}$, $\delta\in\binom{[d]}{d_1}$, $\epsilon\in\binom{[e]}{e_1}$, and $\zeta\in\binom{[z]}{z_1}$. If $\sort(\beta)=\sort(\zeta)$ (and hence $\sort(\gamma)=\sort(\epsilon)$ also), then, using puzzle pieces in $\Hfrak$, then, using puzzle pieces in $\Hfrak$, we have  
    $$\#\{\inlinehex_{\alpha,\beta,\gamma,\delta,\epsilon,\zeta} \text{-puzzles}\} = \#\{\inlinehex_{\alpha,\zeta,\gamma,\delta,\epsilon,\beta} \text{-puzzles}\} = \#\{\inlinehex_{\alpha,\beta,\epsilon,\delta,\gamma,\zeta} \text{-puzzles}\} = \#\{\inlinehex_{\alpha,\zeta,\epsilon,\delta,\gamma,\beta} \text{-puzzles}\}.$$
    In other words, we can commute the pairs $\beta,\zeta$ and $\gamma,\epsilon$ independently and preserve the number of puzzles.
\end{theorem}

\begin{proof}
We give two different arguments. 

\begin{enumerate}
    \item \textbf{Geometric proof}. Given the symmetries in the content of the pairs $\beta,\zeta$ and $\gamma,\epsilon$, we use Proposition \ref{prop:splitLRnumbers} to write 
    $$c_{\beta\gamma,\mu}^{(\epsilon\zeta)^\vee} = c_{\zeta\gamma,\mu}^{(\epsilon\beta)^\vee} = c_{\beta\epsilon,\mu}^{(\gamma\zeta)^\vee} = c_{\zeta\epsilon,\mu}^{(\gamma\beta)^\vee}.$$
    If we make substitutions of these into the expression in the statement of Proposition \ref{prop:hexagonstructureconstantformula}, the sequence of equalities in the theorem statement immediately follows.

    \item \textbf{Proof by commutative property of trapezoidal puzzles.} A puzzle with these boundary conditions can be seen as two trapezoidal puzzles glued together. 
    In Figure \ref{fig:hexagonboundarycommutetwowaysolo}, the line cutting across hexagonal boundary is a line where $\tenstr$ puzzle piece edge labels cannot appear, due to the discrete Green's theorem argument mentioned in Section \ref{section:otherpuzzleproperties}. 
    
    We can sum over each fixed choice of label $\rho$ along that line, and take the product of the number of puzzles in each of the two trapezoidal regions to get the expression
    $$\#\{\inlinehex_{\alpha,\beta,\gamma,\delta,\epsilon,\zeta} \text{-puzzles}\} = \sum_\rho \left(\#\{\inlinetrapside_{\alpha,\beta,\rho,\zeta} \text{-puzzles}\}\right)\left(\#\{\inlinetrapsiderotate_{\rho^\vee,\gamma,\delta,\epsilon} \text{-puzzles}\}\right).$$ 
    
    Then we use the commutative property of trapezoidal puzzles, Theorem \ref{thm:trapezoidcommute}, to write
    \begin{align*}
        &\sum_\rho \left(\#\{\inlinetrapside_{\alpha,\beta,\rho,\zeta} \text{-puzzles}\}\right)\left(\#\{\inlinetrapsiderotate_{\rho^\vee,\gamma,\delta,\epsilon} \text{-puzzles}\}\right) \\
        =& \sum_\rho \left(\#\{\inlinetrapside_{\alpha,\zeta,\rho,\beta} \text{-puzzles}\}\right)\left(\#\{\inlinetrapsiderotate_{\rho^\vee,\gamma,\delta,\epsilon} \text{-puzzles}\}\right) \\
        =&\sum_\rho \left(\#\{\inlinetrapside_{\alpha,\beta,\rho,\zeta} \text{-puzzles}\}\right)\left(\#\{\inlinetrapsiderotate_{\rho^\vee,\epsilon,\delta,\gamma} \text{-puzzles}\}\right) \\
        =&\sum_\rho \left(\#\{\inlinetrapside_{\alpha,\zeta,\rho,\beta} \text{-puzzles}\}\right)\left(\#\{\inlinetrapsiderotate_{\rho^\vee,\epsilon,\delta,\gamma} \text{-puzzles}\}\right).
    \end{align*}
    
This sequence of equalities translates directly to the theorem statement.
\end{enumerate}
\end{proof}

\subsubsection{Commutative property of symmetric pentagonal puzzles}

\begin{figure}[h]
    \centering
    \input{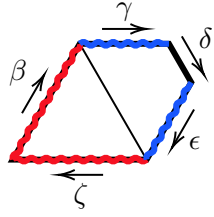}   \caption{ $\protect\inlinepentside_{\beta,\gamma,\delta,\epsilon,\zeta}$, displaying pentagonal symmetry.}
    \label{fig:pentagonboundarycommutesolo}
\end{figure}

A symmetric pentagonal puzzle is one where there are two pairs of boundary sides that have labels with matching content. (See Figure \ref{fig:pentagonboundarycommutesolo}.) This can just be seen as a degenerate hexagonal puzzle with two-way symmetry, and hence the same statements hold, which is why it is included in this section. 

\begin{corollary}[Commutative Property of Symmetric Pentagonal Puzzles]
\label{cor:pentagoncommute}
Let $\beta\in\binom{[b]}{b_1}$, $\gamma\in\binom{[c]}{c_1}$, $\delta\in\binom{[d]}{d_1}$, $\epsilon\in\binom{[e]}{e_1}$, and $\zeta\in\binom{[z]}{z_1}$. If $\sort(\beta)=\sort(\zeta)$ (and hence $\sort(\gamma)=\sort(\epsilon)$ also), then, using puzzle pieces in $\Hfrak$, we have 
$$\#\{\inlinepentside_{\beta,\gamma,\delta,\epsilon,\zeta} \text{-puzzles}\} = \#\{\inlinepentside_{\zeta,\gamma,\delta,\epsilon,\beta} \text{-puzzles}\} = \#\{\inlinepentside_{\beta,\epsilon,\delta,\gamma,\zeta} \text{-puzzles}\} = \#\{\inlinepentside_{\zeta,\epsilon,\delta,\gamma,\beta} \text{-puzzles}\}.$$
    In other words, we can commute the pairs $\beta,\zeta$ and $\gamma,\epsilon$ independently and preserve the number of puzzles.
\end{corollary}

\begin{proof}
    This is just Theorem \ref{thm:hexagontwoway} where $\alpha$ is empty. 
\end{proof}

\subsection{Hexagonal puzzles with three-way symmetry}

\begin{figure}[h]
    \centering
    \input{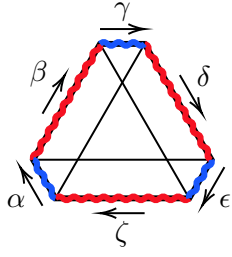}   \caption{$\protect\inlinehex_{\alpha,\beta,\gamma,\delta,\epsilon,\zeta}$, displaying ``three-way symmetry.''}
    \label{fig:hexagonboundarycommutethreewaysolo}
\end{figure}

For hexagonal puzzles, what we will refer to as ``three-way symmetry'' is a situation where there are two triples of boundary sides whose labels have matching content, and where there is a $60^\circ$ angular difference between the sides in each triple. (See Figure \ref{fig:hexagonboundarycommutethreewaysolo}.) 
In Theorem \ref{thm:hexagoncommutethreeway}, we will find that we can permute the labels within each of the triples independently. In Theorem \ref{thm:hexagonthreewayLRnumbers}, we will also present a formula for the number of puzzles with three-way symmetry in terms of Littlewood-Richardson numbers.

\subsubsection{Commutative property of hexagonal puzzles with three-way symmetry}

\begin{theorem}[Commutative Property of Hexagonal Puzzles with Three-Way Symmetry]
    \label{thm:hexagoncommutethreeway}
    Let $\alpha\in\binom{[a]}{a_1}$, $\beta\in\binom{[b]}{b_1}$, $\gamma\in\binom{[c]}{c_1}$, $\delta\in\binom{[d]}{d_1}$, $\epsilon\in\binom{[e]}{e_1}$, and $\zeta\in\binom{[z]}{z_1}$. If $\sort(\alpha)=\sort(\gamma)=\sort(\epsilon)$ (and hence $\sort(\beta)=\sort(\delta)=\sort(\zeta)$ also), then, using puzzle pieces in $\Hfrak$, we have 
$$\#\{\inlinehex_{\alpha,\beta,\gamma,\delta,\epsilon,\zeta} \text{-puzzles}\} = \#\{\inlinehex_{f(\alpha),g(\beta),f(\gamma),g(\delta),f(\epsilon),g(\zeta)} \text{-puzzles}\}$$
    for any pair of bijections $f:\{\alpha,\gamma,\epsilon\}\rightarrow \{\alpha,\gamma,\epsilon\}$ and $g:\{\beta,\delta,\zeta\}\rightarrow \{\beta,\delta,\zeta\}$.
    In other words, we can permute the triples $\alpha,\gamma,\epsilon$ and $\beta,\delta,\zeta$ independently and preserve the number of puzzles. 
\end{theorem}

\begin{proof}

Note that a hexagonal puzzle with three-way symmetry is a special case of a hexagonal puzzle with two-way symmetry, where the two-way symmetry is 3-fold symmetric. We simply apply Theorem \ref{thm:hexagontwoway}, potentially conjugated with $120^\circ$ or $240^\circ$ rotation, as many times as needed in order to get the stated equality for a given pair $f,g$.
For example, to get the equality $\#\{\inlinehex_{\alpha,\beta,\gamma,\delta,\epsilon,\zeta} \text{-puzzles}\}= \#\{\inlinehex_{\gamma,\beta,\alpha,\zeta,\epsilon,\delta} \text{-puzzles}\}$, we rotate the hexagonal boundary $-120^\circ$ and let the pairs $\alpha,\gamma$ and $\delta,\zeta$ replace $\gamma,\epsilon$ and $\zeta,\beta$ in Theorem \ref{thm:hexagontwoway}, respectively. We commute them there, and then we rotate back $120^\circ$. 
\end{proof}

\subsubsection{A formula in terms of Littlewood-Richardson numbers}

\begin{theorem}
\label{thm:hexagonthreewayLRnumbers}
Let $\beta\in\binom{[b]}{b_1}$, $\gamma\in\binom{[c]}{c_1}$, $\delta\in\binom{[d]}{d_1}$, $\epsilon\in\binom{[e]}{e_1}$, and $\zeta\in\binom{[z]}{z_1}$. If $\sort(\alpha)=\sort(\gamma)=\sort(\epsilon)$ (and hence $\sort(\beta)=\sort(\delta)=\sort(\zeta)$ also), then, using puzzle pieces in $\Hfrak$, we have 
$$\#\{\inlinehex_{\alpha,\beta,\gamma,\delta,\epsilon,\zeta} \text{-puzzles}\} =\sum_{\lambda,\mu,\rho,\sigma} c_{\lambda\pad,\mu\pad}^{\beta\gamma} \cdot c_{\rho\pad,\sigma\pad}^{\epsilon\zeta} \cdot c_{\rho,\lambda}^{\delta^\vee}\cdot c_{\mu,\sigma}^{\alpha^\vee}.$$    
\end{theorem}


\begin{proof}
We give two different arguments, one geometric and one combinatorial. 

\begin{enumerate}
    \item \textbf{Geometric proof.} We will begin with the statement of Proposition \ref{prop:hexagonstructureconstantformula}, and given the symmetric content of the strings, we will make the substitutions $c_0=a_0$, $c_1=a_1$, $e_0=a_0$, $e_1=a_1$, and $z_1=b_1$. We will get that $\theta=\mu$. We proceed with basic steps, writing
\begin{gather*}
\#\{\inlinehex_{\alpha,\beta,\gamma,\delta,\epsilon,\zeta} \text{-puzzles}\} =\sum_{\substack{\mu,\theta \text{ such that}\\ \zerostr^{a_0}\mu\onestr^{a_1} = \zerostr^{c_0}\theta\onestr^{e_1}}} c_{\zerostr^{b_0}\alpha^\vee\onestr^{z_1},\theta}^{\onestr^{c_1}\delta\zerostr^{e_0}}\cdot c_{\beta\gamma,\mu}^{(\epsilon\zeta)^\vee} =\sum_{\substack{\mu,\theta \text{ such that}\\ \zerostr^{a_0}\mu\onestr^{a_1} = \zerostr^{a_0}\theta\onestr^{a_1}}} c_{\zerostr^{b_0}\alpha^\vee\onestr^{b_1},\theta}^{\onestr^{a_1}\delta\zerostr^{a_0}} \cdot c_{\beta\gamma,\mu}^{(\epsilon\zeta)^\vee}\\
=\sum_{\mu} c_{\zerostr^{b_0}\alpha^\vee\onestr^{b_1},\mu}^{\onestr^{a_1}\delta\zerostr^{a_0}} \cdot c_{\beta\gamma,\mu}^{(\epsilon\zeta)^\vee} = \sum_{\mu} c_{\zerostr^{b_0}\alpha^\vee\onestr^{b_1},\zerostr^{a_0}\delta^\vee\onestr^{a_1}}^{\mu^\vee} \cdot c_{\beta\gamma,\mu}^{(\epsilon\zeta)^\vee} = \sum_{\mu} c_{(\alpha^\vee)\pad,(\delta^\vee)\pad}^{\mu^\vee} \cdot c_{\beta\gamma,\mu}^{(\epsilon\zeta)^\vee} = \sum_{\mu} d_{\alpha^\vee,\delta^\vee}^{\mu^\vee} \cdot c_{\beta\gamma,\mu}^{(\epsilon\zeta)^\vee}.
\end{gather*}
In the last expression, $d_{\alpha^\vee,\delta^\vee}^{\mu^\vee}$ is the coefficient from Corollary \ref{cor:pushforwardcoeffs}.
We make the substitutions $d_{\alpha^\vee,\delta^\vee}^{\mu^\vee}=d_{\delta^\vee,\alpha^\vee}^{\mu^\vee}$ and $c_{\beta\gamma,\mu}^{(\epsilon\zeta)^\vee} = \underset{\Gr_{b+c}}{\int} [X_{\beta\gamma}][X^{(\epsilon\zeta)^\vee}][X_\mu]$, where $\Gr_{b+c}:= \Gr(b_1+c_1,\langle \bfe_{c+1},\bfe_{c+2},\ldots,\bfe_{2c+b}\rangle)$, into our last expression, and we arrive at
$$\#\{\inlinehex_{\alpha,\beta,\gamma,\delta,\epsilon,\zeta} \text{-puzzles}\} = \sum_{\mu} d_{\delta^\vee,\alpha^\vee}^{\mu^\vee} \underset{\Gr_{b+c}}{\int} [X_{\beta\gamma}][X^{(\epsilon\zeta)^\vee}][X_\mu] =  \underset{\Gr_{b+c}}{\int} [X_{\beta\gamma}][X^{(\epsilon\zeta)^\vee}] \left(\sum_{\mu} d_{\delta^\vee,\alpha^\vee}^{\mu^\vee} [X_\mu]\right).$$

We have $\sum_{\mu} d_{\delta^\vee,\alpha^\vee}^{\mu^\vee} [X_\mu] = \sum_{\mu} d_{\delta^\vee,\alpha^\vee}^{\mu^\vee}[X^{\mu^\vee}] = \Omega_*([X^{\delta^\vee}]\otimes[X^{\alpha^\vee}])$, where $\Omega$ is the map
\begin{gather*}
\Omega: \Gr(b_1,\langle \bfe_{c+1},\bfe_{c+2},\ldots,\bfe_{c+b}\rangle)\times \Gr(c_1,\langle \bfe_{c+b+1},\bfe_{c+b+2},\ldots,\bfe_{2c+b}\rangle) \hookrightarrow  \Gr_{b+c},\\
 (V,W)\mapsto V\oplus W.
\end{gather*}
Going forward, let us shorten our notation for these spaces to $\Gr_b:=\Gr(b_1,\langle \bfe_{c+1},\bfe_{c+2},\ldots,\bfe_{c+b}\rangle)$ and $\Gr_c:=\Gr(c_1,\langle \bfe_{c+b+1},\bfe_{c+b+2},\ldots,\bfe_{2c+b}\rangle)$.

    Now note that by  Lemma \ref{lemma:directsum}, we also have $[X_{\beta\gamma}][X^{(\epsilon\zeta)^\vee}] = \Omega_*([X_\beta][X^{\zeta^\vee}]\otimes [X_\gamma][X^{\epsilon^\vee}])$.
     We make these substitutions where we left off, writing
    \begin{align*}
        \underset{\Gr_{b+c}}{\int} [X_{\beta\gamma}][X^{(\epsilon\zeta)^\vee}] \left(\sum_{\mu} d_{\delta^\vee,\alpha^\vee}^{\mu^\vee} [X_\mu]\right) 
        &=\underset{\Gr_{b+c}}{\int} [X_{\beta\gamma}][X^{(\epsilon\zeta)^\vee}]\Omega_*([X^{\delta^\vee}]\otimes[X^{\alpha^\vee}])\\
        &=\underset{\Gr_{b+c}}{\int} \Omega_*([X_\beta][X^{\zeta^\vee}]\otimes [X_\gamma][X^{\epsilon^\vee}])\Omega_*([X^{\delta^\vee}]\otimes[X^{\alpha^\vee}])\\
        &=\underset{\Gr_{b}\times\Gr_{c}}{\int} \Omega^*(\Omega_*([X_\beta][X^{\zeta^\vee}]\otimes [X_\gamma][X^{\epsilon^\vee}]))([X^{\delta^\vee}]\otimes[X^{\alpha^\vee}])\\
        &=\underset{\Gr_{b}\times\Gr_{c}}{\int} \Omega^*([X_{\beta\gamma}][X^{(\epsilon\zeta)^\vee}])([X^{\delta^\vee}]\otimes[X^{\alpha^\vee}])\\
        &=\underset{\Gr_{b}\times\Gr_{c}}{\int} \Omega^*([X_{\beta\gamma}])\Omega^*([X^{(\epsilon\zeta)^\vee}])([X^{\delta^\vee}]\otimes[X^{\alpha^\vee}]).
    \end{align*}

    Above we included some somewhat redundant intermediate steps, but this is done just to illustrate some of the different forms the expression can take.

Using Corollary \ref{cor:pushforwardcoeffs}, we write $\Omega^*([X_{\beta\gamma}]) = \sum\limits_{\lambda,\mu} c_{\lambda\pad,\mu\pad}^{\beta\gamma} [X_{\lambda}]\otimes [X_{\mu}]$ and $\Omega^*([X^{(\epsilon\zeta)^\vee}])=  \Omega^*([X_{\epsilon\zeta}]) =\sum\limits_{\rho,\sigma} c_{\rho\pad,\sigma\pad}^{\epsilon\zeta} [X_{\rho}]\otimes [X_{\sigma}].$
Then we continue with 
\begin{multline*}
    \underset{\Gr_{b}\times\Gr_{c}}{\int} \Omega^*([X_{\beta\gamma}])\Omega^*([X^{(\epsilon\zeta)^\vee}])([X^{\delta^\vee}]\otimes[X^{\alpha^\vee}])\\
    = \underset{\Gr_{b}\times\Gr_{c}}{\int} \left(\sum_{\lambda,\mu} c_{\lambda\pad,\mu\pad}^{\beta\gamma} [X_{\lambda}]\otimes [X_{\mu}]\right) \left(\sum_{\rho,\sigma} c_{\rho\pad,\sigma\pad}^{\epsilon\zeta} [X_{\rho}]\otimes [X_{\sigma}]\right) ([X^{\delta^\vee}]\otimes[X^{\alpha^\vee}]) \\
    =\sum_{\lambda,\mu,\rho,\sigma} c_{\lambda\pad,\mu\pad}^{\beta\gamma} \cdot c_{\rho\pad,\sigma\pad}^{\epsilon\zeta} \underset{\Gr_{b}}{\int} [X_\lambda][X_\rho][X^{\delta^\vee}] \underset{\Gr_{c}}{\int} [X_\mu][X_\sigma][X^{\alpha^\vee}]\\
    =\sum_{\lambda,\mu,\rho,\sigma} c_{\lambda\pad,\mu\pad}^{\beta\gamma} \cdot c_{\rho\pad,\sigma\pad}^{\epsilon\zeta} \cdot c_{\rho,\lambda}^{\delta^\vee}\cdot c_{\mu,\sigma}^{\alpha^\vee}.
\end{multline*}

\begin{figure}[h]
    \centering
    \begin{subfigure}[t]{0.41\textwidth}
        \centering
        \input{proof_boundary_diagrams/hexagonthreewaysplitup}
        \caption{A hexagonal puzzle with three-way symmetry can be split up into four subpuzzles, two with parallelogram-shaped boundaries and two triangular boundaries.}
        \label{fig:hexagonthreewaysplitup}
    \end{subfigure}
    ~\quad
    \begin{subfigure}[t]{0.41\textwidth}
        \centering
        \input{proof_boundary_diagrams/hexagonthreewaybetagammatriangle}
        \caption{$\#\protect\{\inlinetri_{\sort(\mu)\sort(\lambda),\beta\gamma,(\mu\lambda)^\vee} \text{-puzzles}\}= \#\{\protect\inlinepar_{\beta,\gamma,\lambda^\vee,\mu^\vee} \text{-puzzles}\}$ as a consequence of the unique fillings of the grey regions.}
        \label{fig:hexagonthreewaybetagammatriangle}
    \end{subfigure}
    \caption{}
    \label{fig:hexagonthreewaysplitupfigures}
\end{figure}

\item \textbf{Combinatorial proof.} We can split a $\inlinehex_{\alpha,\beta,\gamma,\delta,\epsilon,\zeta}$-puzzle with three-way symmetry into two parallelogram-shaped subpuzzles and two triangular subpuzzles, as shown in Figure \ref{fig:hexagonthreewaysplitup}. This is because, as a consequence of the discrete Green's theorem statements given in Section \ref{section:otherpuzzleproperties}, it is not possible for any of the drawn lines to have any $\tenstr$ puzzle piece edge labels. By statement (b) of that section, each of the three trapezoidal regions (comprised of a parallelogram and a triangle glued together) cannot have any $\tenstr$s on along its base, which implies that no parallelogram or triangle in the picture can have $\tenstr$s along its boundary either. 

We can fix choices of labels for the internal lines as depicted and sum over them, which gives
\begin{multline*}\#\{\inlinehex_{\alpha,\beta,\gamma,\delta,\epsilon,\zeta} \text{-puzzles}\} \\
=\sum_{\lambda,\mu,\rho,\sigma} \left(\#\{\inlinepar_{\beta,\gamma,\lambda^\vee,\mu^\vee} \text{-puzzles}\}\right) \left(\#\{\inlinepar_{\epsilon,\zeta,\sigma^\vee,\rho^\vee} \text{-puzzles}\}\right) \left(\#\{\inlinetri_{\rho,\lambda,\delta} \text{-puzzles}\}\right) \left(\#\{\inlinetri_{\mu,\sigma,\alpha} \text{-puzzles}\}\right).
\end{multline*}

Then, by completing the first parallelogram to a triangle as in Figure \ref{fig:hexagonthreewaybetagammatriangle} and applying Proposition \ref{prop:equalityofLR}, we can write 
$$ \#\{\inlinepar_{\beta,\gamma,\lambda^\vee,\mu^\vee} \text{-puzzles}\} =  \#\{\inlinetri_{\sort(\mu)\sort(\lambda),\beta\gamma,(\mu\lambda)^\vee} \text{-puzzles}\} = c_{\sort(\mu)\sort(\lambda),\beta\gamma}^{\mu\lambda}  =c_{\mu\pad,\lambda\pad}^{\beta\gamma} =c_{\lambda\pad,\mu\pad}^{\beta\gamma}.$$

Very similarly, we can write 
$$\#\{\inlinepar_{\epsilon,\zeta,\sigma^\vee,\rho^\vee} \text{-puzzles}\} = \#\{\inlinetri_{\sort(\rho)\sort(\sigma),\epsilon\zeta,(\rho\sigma)^\vee} \text{-puzzles}\}= c_{\sort(\rho)\sort(\sigma),\epsilon\zeta}^{\rho\sigma} = c_{\rho\pad,\sigma\pad}^{\epsilon\zeta}.$$

Making the appropriate substitutions, we now have
\begin{multline*}
    \sum_{\lambda,\mu,\rho,\sigma} \left(\#\{\inlinepar_{\beta,\gamma,\lambda^\vee,\mu^\vee} \text{-puzzles}\}\right) \left(\#\{\inlinepar_{\epsilon,\zeta,\sigma^\vee,\rho^\vee} \text{-puzzles}\}\right) \left(\#\{\inlinetri_{\rho,\lambda,\delta} \text{-puzzles}\}\right) \left(\#\{\inlinetri_{\mu,\sigma,\alpha} \text{-puzzles}\}\right) \\
    = \sum_{\lambda,\mu,\rho,\sigma} c_{\lambda\pad,\mu\pad}^{\beta\gamma} \cdot c_{\rho\pad,\sigma\pad}^{\epsilon\zeta} \cdot c_{\rho,\lambda}^{\delta^\vee}\cdot c_{\mu,\sigma}^{\alpha^\vee}.
\end{multline*}
\end{enumerate}
\end{proof}

\subsection{Hexagonal puzzles with all-way symmetry}

\begin{figure}[h]
    \centering
    \input{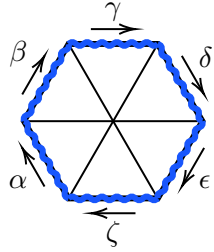}   \caption{$\protect\inlinehex_{\alpha,\beta,\gamma,\delta,\epsilon,\zeta}$, displaying ``all-way symmetry.''}
    \label{fig:hexagonboundarycommuteallwaysolo}
\end{figure}

For hexagonal puzzles, what we will refer to as ``all-way symmetry'' is a situation where the labels on all the boundary sides have matching content. (See Figure \ref{fig:hexagonboundarycommuteallwaysolo}.) 
In Theorem \ref{thm:hexagoncommuteallway}, we will find that we can freely permute all the boundary labels in this case. In Theorem \ref{thm:hexagonallwayLRnumbers}, we will also present a formula for the number of puzzles with three-way symmetry in terms of Littlewood-Richardson numbers. Finally, in Theorem \ref{thm:hexagonallwaynumberofpuzzles}, we give an explicit count of the number of hexagonal puzzles with all-way symmetry. This result in fact implies Theorem \ref{thm:hexagoncommuteallway} and lessens its impact, however we still give Theorem \ref{thm:hexagoncommuteallway} the same treatment and prominence as the many similar theorems in the paper, for the sake of consistency. 

\subsubsection{Commutative property of hexagonal puzzles with all-way symmetry}

\begin{theorem}[Commutative Property of Hexagonal Puzzles with All-Way Symmetry]
    \label{thm:hexagoncommuteallway}
    Let $\alpha\in\binom{[a]}{a_1}$, $\beta\in\binom{[b]}{b_1}$, $\gamma\in\binom{[c]}{c_1}$, $\delta\in\binom{[d]}{d_1}$, $\epsilon\in\binom{[e]}{e_1}$, and $\zeta\in\binom{[z]}{z_1}$. If $\sort(\alpha)=\sort(\beta)=\sort(\gamma)=\sort(\delta)=\sort(\epsilon)=\sort(\zeta)$, then, using puzzle pieces in $\Hfrak$, we have  
    $$\#\{\inlinehex_{\alpha,\beta,\gamma,\delta,\epsilon,\zeta} \text{-puzzles}\} = \#\{\inlinehex_{f(\alpha),f(\beta),f(\gamma),f(\delta),f(\epsilon),f(\zeta)} \text{-puzzles}\}$$
    for any bijection $f:\{\alpha,\beta,\gamma,\delta,\epsilon,\zeta\}\rightarrow \{\alpha,\beta,\gamma,\delta,\epsilon,\zeta\}$. In other words, we can permute all the boundary labels freely and preserve the number of puzzles.
    
\end{theorem}

\begin{proof}
    This type of symmetry is just a combination of opposite sides symmetry and three-way symmetry for hexagonal puzzles. Thus, we can just apply Theorem \ref{thm:hexagonopposite} and Theorem \ref{thm:hexagoncommutethreeway} to prove this result.
\end{proof}

\subsubsection{A formula in terms of Littlewood-Richardson numbers}

\begin{theorem}
\label{thm:hexagonallwayLRnumbers}
Let $\alpha\in\binom{[a]}{a_1}$, $\beta\in\binom{[b]}{b_1}$, $\gamma\in\binom{[c]}{c_1}$, $\delta\in\binom{[d]}{d_1}$, $\epsilon\in\binom{[e]}{e_1}$, and $\zeta\in\binom{[z]}{z_1}$. If $\sort(\alpha)=\sort(\beta)=\sort(\gamma)=\sort(\delta)=\sort(\epsilon)=\sort(\zeta)$, then, using puzzle pieces in $\Hfrak$, we have
$$\#\{\inlinehex_{\alpha,\beta,\gamma,\delta,\epsilon,\zeta} \text{-puzzles}\}  =\sum_{\xi,\tau,\lambda,\mu,\rho,\sigma} c_{\beta,\xi}^\mu\cdot c_{\lambda^\vee,\xi^\vee}^{\gamma^\vee} \cdot c_{\rho,\lambda}^{\delta^\vee}\cdot c_{\epsilon,\tau}^\rho \cdot c_{\sigma^\vee,\tau^\vee}^{\zeta^\vee}  \cdot  c_{\mu,\sigma}^{\alpha^\vee}.$$
\end{theorem}

\begin{proof}
We give two different arguments, one geometric and one combinatorial.
    \begin{enumerate}
        \item \textbf{Geometric proof.} We start with the statement of Theorem \ref{thm:hexagonthreewayLRnumbers}, since we are dealing with a special case of three-way symmetry. We have
        $$\#\{\inlinehex_{\alpha,\beta,\gamma,\delta,\epsilon,\zeta} \text{-puzzles}\} =\sum_{\lambda,\mu,\rho,\sigma} c_{\lambda\pad,\mu\pad}^{\beta\gamma} \cdot c_{\rho\pad,\sigma\pad}^{\epsilon\zeta} \cdot c_{\rho,\lambda}^{\delta^\vee}\cdot c_{\mu,\sigma}^{\alpha^\vee}.$$    

        Using Proposition \ref{prop:equalityofLR}, we have $$c_{\lambda\pad,\mu\pad}^{\beta\gamma} = c_{\sort(\mu)\sort(\lambda),\beta\gamma}^{\mu\lambda}= c_{\mu^\vee\lambda^\vee,\sort(\mu)\sort(\lambda)}^{(\beta\gamma)^\vee}.$$ 
        Now, given the all-way symmetry in string content, we have $\sort(\gamma)=\sort(\lambda)$ and $\sort(\beta)=\sort(\mu)$, hence the above Littlewood-Richardson number corresponds to puzzles with split symmetry. So we use Proposition \ref{prop:splitsymmetryformula} to write 
        $$c_{\lambda^\vee\mu^\vee,\sort(\mu)\sort(\lambda)}^{(\beta\gamma)^\vee} = \sum_{\xi,o}c_{\xi\pad,o\pad}^{\sort(\mu)\sort(\lambda)} \cdot c_{\mu^\vee,\xi}^{\beta^\vee}\cdot c_{\lambda^\vee,o}^{\gamma^\vee}.$$
        Now we look at $c_{\xi\pad,o\pad}^{\sort(\mu)\sort(\lambda)}$, and we note that we are dealing with something similar to what we did in the geometric proof of Theorem \ref{thm:parallelogramLRnumbers} (where we analyzed $c_{\lambda\pad,\mu\pad}^{\sort(\alpha)\sort(\delta)}$), except this time it is much easier since all our strings have the same content. Following the same procedure as in that proof, now we simply get that $c_{\xi\pad,o\pad}^{\sort(\mu)\sort(\lambda)} = \langle \xi,o^\vee \rangle$, so the summands in the above sum are only nonzero when $o=\xi^\vee$. Putting all this together, we end up with
        $$c_{\lambda\pad,\mu\pad}^{\beta\gamma} =  \sum_{\xi} c_{\mu^\vee,\xi}^{\beta^\vee}\cdot c_{\lambda^\vee,\xi^\vee}^{\gamma^\vee}.$$
        By a similar argument, we have that 
        $$c_{\rho\pad,\sigma\pad}^{\epsilon\zeta} =  c_{\sort(\rho)\sort(\sigma),\epsilon\zeta}^{\rho\sigma} = \sum_{\tau} c_{\sigma^\vee,\tau}^{\epsilon^\vee}\cdot c_{\rho^\vee,\tau^\vee}^{\zeta^\vee}.$$
        Now we make these substitutions in our original expression and perform basic steps to get it into the form shown in the theorem.
\begin{multline*}
    \#\{\inlinehex_{\alpha,\beta,\gamma,\delta,\epsilon,\zeta} \text{-puzzles}\}  =\sum_{\lambda,\mu,\rho,\sigma} c_{\lambda\pad,\mu\pad}^{\beta\gamma} \cdot c_{\rho\pad,\sigma\pad}^{\epsilon\zeta} \cdot c_{\rho,\lambda}^{\delta^\vee}\cdot c_{\mu,\sigma}^{\alpha^\vee} =\sum_{\xi,\tau,\lambda,\mu,\rho,\sigma} c_{\mu^\vee,\xi}^{\beta^\vee}\cdot c_{\lambda^\vee,\xi^\vee}^{\gamma^\vee} \cdot c_{\rho^\vee,\tau}^{\epsilon^\vee}\cdot c_{\sigma^\vee,\tau^\vee}^{\zeta^\vee}  \cdot c_{\rho,\lambda}^{\delta^\vee}\cdot c_{\mu,\sigma}^{\alpha^\vee} \\
    =\sum_{\xi,\tau,\lambda,\mu,\rho,\sigma} c_{\beta,\xi}^\mu\cdot c_{\lambda^\vee,\xi^\vee}^{\gamma^\vee} \cdot c_{\epsilon,\tau}^\rho \cdot c_{\sigma^\vee,\tau^\vee}^{\zeta^\vee}  \cdot c_{\rho,\lambda}^{\delta^\vee}\cdot c_{\mu,\sigma}^{\alpha^\vee} =\sum_{\xi,\tau,\lambda,\mu,\rho,\sigma} c_{\beta,\xi}^\mu\cdot c_{\lambda^\vee,\xi^\vee}^{\gamma^\vee} \cdot c_{\rho,\lambda}^{\delta^\vee}\cdot c_{\epsilon,\tau}^\rho \cdot c_{\sigma^\vee,\tau^\vee}^{\zeta^\vee}  \cdot  c_{\mu,\sigma}^{\alpha^\vee}.
\end{multline*}

        \item \textbf{Combinatorial proof.} We start where the combinatorial proof of Theorem \ref{thm:hexagonthreewayLRnumbers} leaves off. Now the parallelogram-shaped regions are rhombi, and each rhombi-shaped subpuzzle can be further split into two triangular subpuzzles. This is again a consequence of the discrete Green's theorem, where now each of the six different trapezoidal regions shown in Figure \ref{fig:hexagonallwaysplitup} cannot have any $\tenstr$s along their boundary, so neither can any of the six triangular regions. 

\begin{figure}[h]
    \centering
    \input{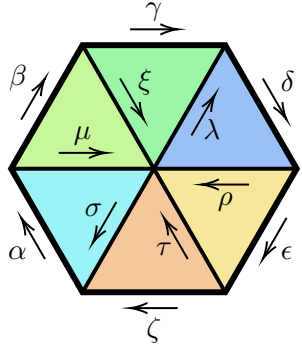}
    \caption{A hexagonal puzzle with all-way symmetry can be split up into six triangular subpuzzles.}
    \label{fig:hexagonallwaysplitup}
\end{figure}

        We fix choices of labels for the internal lines as depicted and sum over them, which gives
        \begin{multline*}
             \#\{\inlinehex_{\alpha,\beta,\gamma,\delta,\epsilon,\zeta} \text{-puzzles}\} \\
            =\sum_{\lambda,\mu,\rho,\sigma} \left(\#\{\inlinetri_{\beta,\xi,\mu^\vee} \text{-puzzles}\}\right) \left(\#\{\inlinetri_{\lambda^\vee,\xi^\vee,\gamma} \text{-puzzles}\}\right) \left(\#\{\inlinetri_{\rho,\lambda,\delta} \text{-puzzles}\}\right)  \\
            \left(\#\{\inlinetri_{\epsilon,\tau,\rho^\vee} \text{-puzzles}\}\right) \left(\#\{\inlinetri_{\sigma^\vee,\tau^\vee,\zeta} \text{-puzzles}\}\right) \left(\#\{\inlinetri_{\mu,\sigma,\alpha} \text{-puzzles}\}\right) \\
            =\sum_{\xi,\tau,\lambda,\mu,\rho,\sigma} c_{\beta,\xi}^\mu\cdot c_{\lambda^\vee,\xi^\vee}^{\gamma^\vee} \cdot c_{\rho,\lambda}^{\delta^\vee}\cdot c_{\epsilon,\tau}^\rho \cdot c_{\sigma^\vee,\tau^\vee}^{\zeta^\vee}  \cdot  c_{\mu,\sigma}^{\alpha^\vee}.
        \end{multline*}
    \end{enumerate}
\end{proof}

\subsubsection{An explicit count of hexagonal puzzles with all-way symmetry}

\begin{theorem}
\label{thm:hexagonallwaynumberofpuzzles}
Let $\alpha\in\binom{[a]}{a_1}$, $\beta\in\binom{[b]}{b_1}$, $\gamma\in\binom{[c]}{c_1}$, $\delta\in\binom{[d]}{d_1}$, $\epsilon\in\binom{[e]}{e_1}$, and $\zeta\in\binom{[z]}{z_1}$. If $\sort(\alpha)=\sort(\beta)=\sort(\gamma)=\sort(\delta)=\sort(\epsilon)=\sort(\zeta)$, then, using puzzle pieces in $\Hfrak$, we have
$$\#\{\inlinehex_{\alpha,\beta,\gamma,\delta,\epsilon,\zeta} \text{-puzzles}\} =  \begin{cases}
        \binom{a}{a_1} &\text{if } \alpha=\beta=\gamma=\delta=\epsilon=\zeta=\zerostr^{a_0}\onestr^{a_1} \\
        0 &\text{otherwise}
        \end{cases}.$$
\end{theorem}

\begin{proof}
Note that we have $\alpha,\beta,\gamma,\delta,\epsilon,\zeta\in \binom{[a]}{a_1}$, since these strings all have the same content. Also let $\mu\in \binom{[a]}{a_1}$. Now consider the classes $[X_{\mu^\vee}],[X_\alpha],[X_\beta],[X_\gamma],[X_\delta],[X_\epsilon],[X_\zeta] \in H^*(\Gr(a_1,\C^a))$. We will take their product and expand it in the Schubert basis as follows.
\begin{align*}
   & [X_{\mu^\vee}][X_{\beta}][X_{\gamma}][X_{\delta}][X_{\epsilon}][X_{\zeta}][X_{\alpha}] \\
   =& ([X_{\mu^\vee}][X_{\beta}])[X_{\gamma}][X_{\delta}][X_{\epsilon}][X_{\zeta}][X_{\alpha}] = \left(\sum_{\xi} c_{\mu^\vee,\beta}^{\xi^\vee} [X_{\xi^\vee}]\right)[X_{\gamma}][X_{\delta}][X_{\epsilon}][X_{\zeta}][X_{\alpha}] \\
   =& \sum_{\xi} c_{\mu^\vee,\beta}^{\xi^\vee} ([X_{\xi^\vee}][X_{\gamma}])[X_{\delta}][X_{\epsilon}][X_{\zeta}][X_{\alpha}] = \sum_{\xi} c_{\mu^\vee,\beta}^{\xi^\vee} \left(\sum_{\lambda} c_{\xi^\vee,\gamma}^{\lambda} [X_\lambda] \right) [X_{\delta}][X_{\epsilon}][X_{\zeta}][X_{\alpha}] \\
   =& \sum_{\xi,\lambda} c_{\mu^\vee,\beta}^{\xi^\vee} \cdot c_{\xi^\vee,\gamma}^{\lambda} 
 ([X_\lambda] [X_{\delta}])[X_{\epsilon}][X_{\zeta}][X_{\alpha}] = \sum_{\xi,\lambda} c_{\mu^\vee,\beta}^{\xi^\vee} \cdot c_{\xi^\vee,\gamma}^{\lambda} 
 \left(\sum_{\rho} c_{\lambda,\delta}^{\rho^\vee} [X_{\rho^\vee}] \right)[X_{\epsilon}][X_{\zeta}][X_{\alpha}] \\
 =& \sum_{\xi,\lambda,\rho} c_{\mu^\vee,\beta}^{\xi^\vee} \cdot c_{\xi^\vee,\gamma}^{\lambda} \cdot c_{\lambda,\delta}^{\rho^\vee}
 ([X_{\rho^\vee}][X_{\epsilon}])[X_{\zeta}][X_{\alpha}] = \sum_{\xi,\lambda,\rho} c_{\mu^\vee,\beta}^{\xi^\vee} \cdot c_{\xi^\vee,\gamma}^{\lambda} \cdot c_{\lambda,\delta}^{\rho^\vee}
 \left(\sum_{\tau} c_{\rho^\vee,\epsilon}^{\tau^\vee} [X_{\tau^\vee}]\right)[X_{\zeta}][X_{\alpha}] \\
 =& \sum_{\xi,\lambda,\rho,\tau} c_{\mu^\vee,\beta}^{\xi^\vee} \cdot c_{\xi^\vee,\gamma}^{\lambda} \cdot c_{\lambda,\delta}^{\rho^\vee}\cdot c_{\rho^\vee,\epsilon}^{\tau^\vee}
 ([X_{\tau^\vee}][X_{\zeta}])[X_{\alpha}] = \sum_{\xi,\lambda,\rho,\tau} c_{\mu^\vee,\beta}^{\xi^\vee} \cdot c_{\xi^\vee,\gamma}^{\lambda} \cdot c_{\lambda,\delta}^{\rho^\vee}\cdot c_{\rho^\vee,\epsilon}^{\tau^\vee}
 \left(\sum_{\sigma} c_{\tau^\vee,\zeta}^{\sigma} [X_{\sigma}]\right) [X_{\alpha}] \\
 =& \sum_{\xi,\lambda,\rho,\tau,\sigma} c_{\mu^\vee,\beta}^{\xi^\vee} \cdot c_{\xi^\vee,\gamma}^{\lambda} \cdot c_{\lambda,\delta}^{\rho^\vee}\cdot c_{\rho^\vee,\epsilon}^{\tau^\vee}\cdot c_{\tau^\vee,\zeta}^{\sigma}
 ([X_{\sigma}][X_{\alpha}]) = \sum_{\xi,\lambda,\rho,\tau,\sigma} c_{\mu^\vee,\beta}^{\xi^\vee} \cdot c_{\xi^\vee,\gamma}^{\lambda} \cdot c_{\lambda,\delta}^{\rho^\vee}\cdot c_{\rho^\vee,\epsilon}^{\tau^\vee}\cdot c_{\tau^\vee,\zeta}^{\sigma}
 \left(\sum_{\upsilon} c_{\sigma,\alpha}^{\upsilon^\vee}[X_{\upsilon^\vee}]\right) \\
 =& \sum_{\xi,\lambda,\rho,\tau,\sigma,\upsilon} c_{\mu^\vee,\beta}^{\xi^\vee} \cdot c_{\xi^\vee,\gamma}^{\lambda} \cdot c_{\lambda,\delta}^{\rho^\vee}\cdot c_{\rho^\vee,\epsilon}^{\tau^\vee}\cdot c_{\tau^\vee,\zeta}^{\sigma}\cdot c_{\sigma,\alpha}^{\upsilon^\vee} [X_{\upsilon^\vee}]
\end{align*}

Now we use the equality derived above to write
\begin{multline*}
     \underset{\Gr(a_1,\C^a)}{\int} [X_{\mu^\vee}][X_{\beta}][X_{\gamma}][X_{\delta}][X_{\epsilon}][X_{\zeta}][X_{\alpha}][X^{\mu^\vee}] = \underset{\Gr(a_1,\C^a)}{\int} \sum_{\xi,\lambda,\rho,\tau,\sigma,\upsilon} c_{\mu^\vee,\beta}^{\xi^\vee} \cdot c_{\xi^\vee,\gamma}^{\lambda} \cdot c_{\lambda,\delta}^{\rho^\vee}\cdot c_{\rho^\vee,\epsilon}^{\tau^\vee}\cdot c_{\tau^\vee,\zeta}^{\sigma}\cdot c_{\sigma,\alpha}^{\upsilon^\vee} [X_{\upsilon^\vee}][X^{\mu^\vee}] \\
    = \sum_{\xi,\lambda,\rho,\tau,\sigma,\upsilon} c_{\mu^\vee,\beta}^{\xi^\vee} \cdot c_{\xi^\vee,\gamma}^{\lambda} \cdot c_{\lambda,\delta}^{\rho^\vee}\cdot c_{\rho^\vee,\epsilon}^{\tau^\vee}\cdot c_{\tau^\vee,\zeta}^{\sigma}\cdot c_{\sigma,\alpha}^{\upsilon^\vee} \underset{\Gr(a_1,\C^a)}{\int} [X_{\upsilon^\vee}][X^{\mu^\vee}] = \sum_{\xi,\lambda,\rho,\tau,\sigma,\upsilon} c_{\mu^\vee,\beta}^{\xi^\vee} \cdot c_{\xi^\vee,\gamma}^{\lambda} \cdot c_{\lambda,\delta}^{\rho^\vee}\cdot c_{\rho^\vee,\epsilon}^{\tau^\vee}\cdot c_{\tau^\vee,\zeta}^{\sigma}\cdot c_{\sigma,\alpha}^{\upsilon^\vee} \langle \upsilon^\vee,\mu^\vee\rangle \\
    = \sum_{\xi,\lambda,\rho,\tau,\sigma} c_{\mu^\vee,\beta}^{\xi^\vee} \cdot c_{\xi^\vee,\gamma}^{\lambda} \cdot c_{\lambda,\delta}^{\rho^\vee}\cdot c_{\rho^\vee,\epsilon}^{\tau^\vee}\cdot c_{\tau^\vee,\zeta}^{\sigma}\cdot c_{\sigma,\alpha}^{\mu^\vee},
\end{multline*}
where the last equality in the sequence comes from noting that the summand equals 0 whenever $\upsilon\neq \mu$.

We will now also argue that 
\begin{multline*}
    \underset{\Gr(a_1,\C^a)}{\int} [X_{\mu^\vee}][X_{\beta}][X_{\gamma}][X_{\delta}][X_{\epsilon}][X_{\zeta}][X_{\alpha}][X^{\mu^\vee}] = \begin{cases}
        1 &\text{if } \alpha=\beta=\gamma=\delta=\epsilon=\zeta=\zerostr^{a_0}\onestr^{a_1} \\
        0 &\text{otherwise}
    \end{cases}.
\end{multline*}

First suppose that $\alpha=\beta=\gamma=\delta=\epsilon=\zeta=\zerostr^{a_0}\onestr^{a_1}$. Then $[X_\alpha]=[X_\beta]=[X_\gamma]=[X_\delta]=[X_\epsilon]=[X_\zeta]=[X_{\zerostr^{a_0}\onestr^{a_1}}]=[\Gr(a_1,\C^a)]=1$, and we have 
$$\underset{\Gr(a_1,\C^a)}{\int} [X_{\mu^\vee}][X_{\beta}][X_{\gamma}][X_{\delta}][X_{\epsilon}][X_{\zeta}][X_{\alpha}][X^{\mu^\vee}] = \underset{\Gr(a_1,\C^a)}{\int} [X_{\mu^\vee}][X^{\mu^\vee}]= \langle \mu^\vee,\mu^\vee \rangle =1.$$
On the other hand, suppose that one of $\alpha,\beta,\gamma,\delta,\epsilon,\zeta$ is not equal to $\zerostr^{a_0}\onestr^{a_1}$. Without loss of generality, suppose that $\alpha\neq \zerostr^{a_0}\onestr^{a_1}$. Recall that $\deg([X_\alpha])=\codim(X_\alpha)=\ell(\alpha)$, where $\ell(\alpha)$ is the length of $\alpha$. Then we have $\deg([X_\alpha])>0$, because the string $\alpha$ must have at least one inversion. Since we already necessarily have that $\deg([X_{\mu^\vee}])+\deg([X^{\mu^\vee}]) = \dim(\Gr(a_1,\C^a))$, this then implies that $\deg([X_{\mu^\vee}])+\deg([X^{\mu^\vee}]) +  \deg([X_\alpha])>\dim(\Gr(a_1,\C^a))$. Thus $[X_{\mu^\vee}][X^{\mu^\vee}][X_\alpha]= 0$, and likewise $[X_{\mu^\vee}][X_{\beta}][X_{\gamma}][X_{\delta}][X_{\epsilon}][X_{\zeta}][X_{\alpha}][X^{\mu^\vee}]=0$. So in this case we have
$$\underset{\Gr(a_1,\C^a)}{\int} [X_{\mu^\vee}][X_{\beta}][X_{\gamma}][X_{\delta}][X_{\epsilon}][X_{\zeta}][X_{\alpha}][X^{\mu^\vee}] = \underset{\Gr(a_1,\C^a)}{\int} 0=0.$$
This finishes the argument for the above claim.

At this point, we have derived two different expressions for 
$\underset{\Gr(a_1,\C^a)}{\int} [X_{\mu^\vee}][X_{\beta}][X_{\gamma}][X_{\delta}][X_{\epsilon}][X_{\zeta}][X_{\alpha}][X^{\mu^\vee}]$, and we combine them to obtain
$$\sum_{\xi,\lambda,\rho,\tau,\sigma} c_{\mu^\vee,\beta}^{\xi^\vee} \cdot c_{\xi^\vee,\gamma}^{\lambda} \cdot c_{\lambda,\delta}^{\rho^\vee}\cdot c_{\rho^\vee,\epsilon}^{\tau^\vee}\cdot c_{\tau^\vee,\zeta}^{\sigma}\cdot c_{\sigma,\alpha}^{\mu^\vee}= \begin{cases}
        1 &\text{if } \alpha=\beta=\gamma=\delta=\epsilon=\zeta=\zerostr^{a_0}\onestr^{a_1} \\
        0 &\text{otherwise}
    \end{cases} . $$

Now if we additionally sum over all $\mu\in \binom{[a]}{a_1}$, we obtain
\begin{align*}
    \sum_{\mu}\sum_{\xi,\lambda,\rho,\tau,\sigma} c_{\mu^\vee,\beta}^{\xi^\vee} \cdot c_{\xi^\vee,\gamma}^{\lambda} \cdot c_{\lambda,\delta}^{\rho^\vee}\cdot c_{\rho^\vee,\epsilon}^{\tau^\vee}\cdot c_{\tau^\vee,\zeta}^{\sigma}\cdot c_{\sigma,\alpha}^{\mu^\vee} &= \sum_{\mu} \begin{cases}
        1 &\text{if } \alpha=\beta=\gamma=\delta=\epsilon=\zeta=\zerostr^{a_0}\onestr^{a_1} \\
        0 &\text{otherwise}
    \end{cases} \\
    &= \begin{cases}
        \binom{a}{a_1} &\text{if } \alpha=\beta=\gamma=\delta=\epsilon=\zeta=\zerostr^{a_0}\onestr^{a_1} \\
        0 &\text{otherwise}
        \end{cases}.
\end{align*}

Finally, we apply Theorem \ref{thm:hexagonallwayLRnumbers} along with the above to write
\begin{multline*}
    \#\{\inlinehex_{\alpha,\beta,\gamma,\delta,\epsilon,\zeta} \text{-puzzles}\}  =\sum_{\xi,\tau,\lambda,\mu,\rho,\sigma} c_{\beta,\xi}^\mu\cdot c_{\lambda^\vee,\xi^\vee}^{\gamma^\vee} \cdot c_{\rho,\lambda}^{\delta^\vee}\cdot c_{\epsilon,\tau}^\rho \cdot c_{\sigma^\vee,\tau^\vee}^{\zeta^\vee}  \cdot  c_{\mu,\sigma}^{\alpha^\vee} \\
    = \sum_{\mu}\sum_{\xi,\lambda,\rho,\tau,\sigma} c_{\mu^\vee,\beta}^{\xi^\vee} \cdot c_{\xi^\vee,\gamma}^{\lambda} \cdot c_{\lambda,\delta}^{\rho^\vee}\cdot c_{\rho^\vee,\epsilon}^{\tau^\vee}\cdot c_{\tau^\vee,\zeta}^{\sigma}\cdot c_{\sigma,\alpha}^{\mu^\vee}  =  \begin{cases}
        \binom{a}{a_1} &\text{if } \alpha=\beta=\gamma=\delta=\epsilon=\zeta=\zerostr^{a_0}\onestr^{a_1} \\
        0 &\text{otherwise}
        \end{cases},
\end{multline*}
where the second equality in the sequence simply comes from rearranging variables and from basic symmetries of Littlewood-Richardson numbers.
\end{proof}

\begin{remark}
    Theorem \ref{thm:hexagonallwaynumberofpuzzles} implies Theorem \ref{thm:hexagoncommuteallway} and renders it less interesting, since Theorem \ref{thm:hexagonallwaynumberofpuzzles} says that the only case where there exist any puzzles is when all the boundary labels are equal to $\zerostr^{a_0}\onestr^{a_1}$, and so permuting them is trivial. In any other case, the number of puzzles always equals 0 regardless of how we permute the boundary labels.
\end{remark}

\section{2-step and 3-step puzzles}
\label{section:2step3step}
In Section \ref{section:propertiesofpuzzles2step3step}, we discussed how all the puzzle properties we use in this paper also work for 2-step and 3-step puzzles. In Section \ref{Section:geometriclemmasdstep}, we gave $d$-step analogues of the geometric lemmas that underpin all of our geometric proofs in this paper. Then, essentially, any result or proof we gave for Grassmannian puzzles which invoked only these properties or lemmas can be straightforwardly replicated for 2-step and 3-step puzzles. 

Specifically, for 2-step and 3-step puzzles we obtain the analogous statements of:
\begin{itemize}
    \item all of the ``complete to a triangle'' operations (i.e. Propositions \ref{prop:trapezoidcomplete}, \ref{prop:parallelogramcomplete}, and \ref{prop:hexagoncomplete}) to give geometric interpretations to polygonal puzzles.
    \item Theorems \ref{thm:puzzlecommutesplit}, \ref{thm:trapezoidcommute}, \ref{thm:Hparallelogramcommute}, \ref{thm:rhombuscommute}, and \ref{thm:hexagontwoway}, and Corollary \ref{cor:pentagoncommute} (i.e. the commutative properties of puzzles with split symmetry, trapezoidal puzzles, parallelogram-shaped puzzles, symmetric rhombus-shaped puzzles, hexagonal puzzles with two-way symmetry, and symmetric pentagonal puzzles).
    \item Theorems \ref{thm:hexagonthreewayLRnumbers} and \ref{thm:hexagonallwayLRnumbers} (i.e. the formulas in terms of Littlewood-Richardson numbers for hexagonal puzzles with three-way symmetry and hexagonal puzzles with all-way symmetry), via the combinatorial arguments, and Theorem \ref{thm:hexagonallwaynumberofpuzzles} (i.e. the explicit count of hexagonal puzzles with all-way symmetry).
    \item Remarks \ref{remark:splitsymmetryKtheory}, \ref{remark:Ktheorytrapezoidpuzzles}, \ref{remark:Ktheoryparallelogrampuzzles} (i.e. the K-theory analogues of Proposition \ref{prop:splitLRnumbers} and Theorems \ref{thm:puzzlecommutesplit}, \ref{thm:trapezoidcommute}, and \ref{thm:Hparallelogramcommute}.
    
    \explainalot{
    \item Theorem \ref{thm:puzzlecommutesplit} (commutative property of puzzles with split symmetry).
    \begin{itemize}
        \item A geometric argument via an analogous version of Proposition \ref{prop:splitLRnumbers} (invoking Lemma \ref{lemma:commutesplitdstep}) for $d$-step structure constants.
        \item A combinatorial argument using the commutative properties of triangular and trapezoidal puzzles (see Remark 
    \end{itemize}
    \item Theorem \ref{thm:trapezoidcommute} (commutative property of trapezoidal puzzles).
    \item Theorem \ref{thm:Hparallelogramcommute} (commutative property of parallelogram-shaped puzzles).
    \item Theorem \ref{thm:rhombuscommute} (commutative property of symmetric rhombus-shaped puzzles).
    \item Theorem \ref{thm:hexagontwoway} (commutative property of hexagonal puzzles with two-way symmetry).
    \item Corollary \ref{cor:pentagoncommute} (commutative property of symmetric pentagonal puzzles).
    \item Theorem \ref{thm:hexagontwoway} (commutative property of hexagonal puzzles with three-way symmetry).
    \item Theorem \ref{thm:hexagonthreewayLRnumbers} (formula in terms of Littlewood-Richardson numbers for hexagonal puzzles with three-way symmetry). The combinatorial proof.
    \item Theorem \ref{thm:hexagonallwayLRnumbers} (formula in terms of Littlewood-Richardson numbers for hexagonal puzzles with all-way symmetry). The combinatorial proof.
    }
\end{itemize}

On the other hand, what cannot be straightforwardly replicated is any result or proof which relied on Proposition \ref{prop:equalityofLR} (which was an equality of Littlewood-Richardson numbers proved via skew tableaux rules). This includes Corollary \ref{cor:pushforwardcoeffs} about the pushforward and pullback coefficients for the direct sum map, along with anything it was used for (namely the geometric proofs of the formulas in terms of Littlewood-Richardson numbers given in each section). We also have not attempted a 2-step or 3-step version of Lemma \ref{lemma:uniquefilling} or a puzzle-based strategy to prove a statement analogous to Theorem \ref{thm:parallelogramLRnumbers}.

\justthesis{
Anywhere we have a figure with a line cutting across a boundary shape, there also can be no composite labels appearing along it for a 2-step or 3-step puzzle. Anywhere we were able to split a puzzle into subpuzzles, we can do the same for 2-step and 3-step puzzles.}

\explainalot{
All of the ``complete to a triangle'' operations to give geometric interpretations to polygonal puzzles work the same for 2-step and 3-step puzzles.

All the commutative properties can be proved for 2-step and 3-step.

\begin{itemize}
    \item \textbf{Split symmetry.} There is an analogous version of Proposition \ref{prop:splitLRnumbers} (invoking Lemma \ref{lemma:commutesplitdstep}) for $d$-step structure constants. Theorem \ref{thm:puzzlecommutesplit} (the commutative property of puzzles with split symmetry) can be proved in that way, or via the commutative properties of triangular and trapezoidal puzzles. 

    \item Theorem \ref{thm:puzzlecommutesplit} (the commutative property of puzzles with split symmetry).
    \begin{itemize}
        \item A geometric argument via an analogous version of Proposition \ref{prop:splitLRnumbers} (invoking Lemma \ref{lemma:commutesplitdstep}) for $d$-step structure constants.
        \item A combinatorial argument using the commutative properties of triangular and trapezoidal puzzles. 
    \end{itemize}
    \item Theorem \ref{thm:trapezoidcommute} (commutative property of trapezoidal puzzles).
\end{itemize}
}

\section{Further questions}
\label{section:furtherquestions}

The following are open questions that could be pursued in future work. 

\begin{enumerate}

\item Find a commutativity statement for hexagonal equivariant puzzles. This may be significantly more difficult, and little work has been done on it so far.  Commuting the boundary labels would at least have to correspond to some permutation of the $y_i$ indeterminates, as it did in the parallelogram case. We would need to determine what the expected permutation would be for various cases of special hexagonal symmetry, if there is one at all. It is uncertain whether there even exists a nice statement to be made, as we did not seem to find one for equivariant trapezoidal puzzles.

\item 
Find a commutativity statement for parallelogram-shaped equivariant K-theory puzzles, i.e. using pieces in $\Hfrak \cup \left\{\raisebox{-.4\height}{
\begin{tikzpicture}[x=1pt,y=1pt,yscale=1,xscale=1,scale=.68]
\draw [fill={rgb, 255:red, 255; green, 246; blue, 74 }] (0.0,0.0) -- (15.0,25.98) -- (30.0,0.0) -- (15.0,-25.98) -- cycle;
\draw (7.5,12.99) node  [rotate=60]  {$\mathtt{0}$};
\draw (22.5,12.99) node  [rotate=-60]  {$\mathtt{1}$};
\draw (22.5,-12.99) node  [rotate=60]  {$\mathtt{0}$};
\draw (7.5,-12.99) node  [rotate=-60]  {$\mathtt{1}$};
\end{tikzpicture}
}, \raisebox{-.4\height}{
\begin{tikzpicture}[x=1pt,y=1pt,yscale=1,xscale=1,scale=.68]
\draw [fill={rgb, 255:red, 255; green, 255; blue, 255 }] (0.0,0.0) -- (15.0,25.98) -- (30.0,0.0) -- cycle;
\draw (7.5,12.99) node  [rotate=60]  {$\mathtt{1\! 0}$};
\draw (22.5,12.99) node  [rotate=-60]  {$\mathtt{1\! 0}$};
\draw (15.0,0.0) node  [rotate=0]  {$\mathtt{1\! 0}$};
\end{tikzpicture}
}\right\}$. All that is known is that the statement and proof method of Theorem \ref{thm:eqvtparallelogram} do not directly extend to this case, since e.g. it is not true that $[\Ocal_{\alpha\gamma}][\Ical^{(\delta\beta)^\vee}]=\Phi_a\cdot ([\Ocal_{\beta\gamma}][\Ical^{(\delta\alpha)^\vee}])$ or $[\Ocal_{\alpha\gamma}][\Ical^{(\delta\beta)^\vee}]=\Phi_c\cdot ([\Ocal_{\alpha\delta}][\Ical^{(\gamma\beta)^\vee}])$. The K-theory analogue of Lemma \ref{lemma:commutesplit} (i.e. Remark \ref{remark:lemmatwoktheory}) is not useful in the equivariant case.


\item Extend our commutativity results to the types of puzzles that include the piece \raisebox{-.4\height}{
\begin{tikzpicture}[x=1pt,y=1pt,yscale=1,xscale=1,scale=.68]
\draw [fill={rgb, 255:red, 255; green, 255; blue, 255 }] (30.0,25.98) -- (15.0,0.0) -- (0.0,25.98) -- cycle;
\draw (22.5,12.99) node  [rotate=60]  {$\mathtt{1\! 0}$};
\draw (7.5,12.99) node  [rotate=-60]  {$\mathtt{1\! 0}$};
\draw (15.0,25.98) node  [rotate=0]  {$\mathtt{1\! 0}$};
\end{tikzpicture}
}, namely those that compute structure constants for ideal sheaf K-theory classes and SSM classes. When using this piece, we cannot perform the operation of uniquely completing a polygonal puzzle to a triangular puzzle and get a bijection as we do for the cases in this paper, since there is no unique identity puzzle. Quite little work has been done on this, but so far in computer experiments, we found that the sum of the weights of $\inlinepar_{\alpha,\gamma,\beta,\delta}$-puzzles using pieces in $\Hfrak \cup \left\{\raisebox{-.4\height}{
\begin{tikzpicture}[x=1pt,y=1pt,yscale=1,xscale=1,scale=.68]
\draw [fill={rgb, 255:red, 255; green, 255; blue, 255 }] (0.0,0.0) -- (15.0,25.98) -- (30.0,0.0) -- cycle;
\draw (7.5,12.99) node  [rotate=60]  {$\mathtt{1\! 0}$};
\draw (22.5,12.99) node  [rotate=-60]  {$\mathtt{1\! 0}$};
\draw (15.0,0.0) node  [rotate=0]  {$\mathtt{1\! 0}$};
\end{tikzpicture}
},\raisebox{-.4\height}{
\begin{tikzpicture}[x=1pt,y=1pt,yscale=1,xscale=1,scale=.68]
\draw [fill={rgb, 255:red, 255; green, 255; blue, 255 }] (30.0,25.98) -- (15.0,0.0) -- (0.0,25.98) -- cycle;
\draw (22.5,12.99) node  [rotate=60]  {$\mathtt{1\! 0}$};
\draw (7.5,12.99) node  [rotate=-60]  {$\mathtt{1\! 0}$};
\draw (15.0,25.98) node  [rotate=0]  {$\mathtt{1\! 0}$};
\end{tikzpicture}
} \right\}$ is not preserved when the opposite boundary labels are swapped. It is possible that the ``correct'' thing to do would be to just see if we can commute the pairs $\alpha,\beta$ and $\gamma,\delta$ for puzzles with the boundary $\inlinetri_{\sort(\alpha)\gamma,\beta\sort(\delta),\delta\alpha}$. This would no longer be about the parallelogram per se; rather, it would just be triangular boundary conditions of a parallelogram-inspired form.

\item Extend our commutativity results to hexagonal K-theory puzzles using pieces in $\Hfrak \cup \left\{\raisebox{-.4\height}{
\begin{tikzpicture}[x=1pt,y=1pt,yscale=1,xscale=1,scale=.68]
\draw [fill={rgb, 255:red, 255; green, 255; blue, 255 }] (0.0,0.0) -- (15.0,25.98) -- (30.0,0.0) -- cycle;
\draw (7.5,12.99) node  [rotate=60]  {$\mathtt{1\! 0}$};
\draw (22.5,12.99) node  [rotate=-60]  {$\mathtt{1\! 0}$};
\draw (15.0,0.0) node  [rotate=0]  {$\mathtt{1\! 0}$};
\end{tikzpicture}
}\right\}$. Our proofs of the hexagonal commutativity properties in $H^*$ used Corollary \ref{cor:pushforwardcoeffs}, but the analogous coefficients for K-theory are not currently understood by the author, which has hindered an attempt to carry out a direct translation of those proof methods.

\item Find more direct combinatorial proofs, perhaps on the level of actual puzzle pieces, for all the commutativity results in this paper. The proofs given are either geometric, or they are combinatorial proofs that still first rely on the operation of completing a polygonal puzzle to a triangular puzzle and working in that context. We have no more native combinatorial understanding of these results. Trying to apply the strategy of migration as in \cite{Purbhoo} might be fruitful for this question, and in that paper it is already used to commute boundary labels for rhombus-shaped puzzles.

\item Prove 2-step and 3-step versions of any of our results that cannot already be straightforwardly replicated based on our 1-step arguments (see Section \ref{section:2step3step}). This includes commutative properties of hexagonal puzzles with opposite sides symmetry and hexagonal puzzles with all-way symmetry (i.e. prove analogues of Theorem \ref{thm:hexagonopposite} and Theorem \ref{thm:hexagoncommuteallway}). (The barrier to replicating the proof method used for Theorem \ref{thm:hexagonopposite} is that it traces back to Proposition \ref{prop:equalityofLR}, for which we do not have a 2-step or 3-step analogue currently.)

\item Extend our results to 4-step puzzles, and perhaps even hypothetical higher step puzzles. Assuming that 4-step puzzles, and any higher step puzzle rule that might be found in the future, share (approximately) the same properties as given in Section \ref{section:propertiesofpuzzles}, it is likely that many our results would extend to them, in a similar manner as for 2-step and 3-step puzzles, as described in Section \ref{section:2step3step}.


\end{enumerate}

\section*{Acknowledgements}
We thank Allen Knutson for the questions that first led to this work and for his help. We thank Pasha Pylyavskyy for inspiring conversations and questions about trapezoidal and parallelogram-shaped puzzles. We thank Paul Zinn-Justin for his encouragement and support.

\printbibliography

\end{document}